\documentclass[reqno]{amsart}
\usepackage{amssymb}
\usepackage{amsfonts}
\usepackage{amsmath}
\usepackage{amsthm, amsmath, amssymb, enumerate}
\usepackage{amssymb, enumitem}
\usepackage{bbm}
\usepackage{amscd}
\usepackage[all]{xy}
\usepackage{hyperref}
\usepackage{cleveref}
\usepackage{amsmath}
\usepackage{amssymb}
\usepackage{amsthm}
\usepackage[sort&compress,numbers]{natbib}
\usepackage[left=2cm,right=2cm,bottom=3cm,top=3cm]{geometry}
\usepackage{tikz}
\usepackage{tikz-cd}
\usepackage{wrapfig}
\usepackage{refstyle}

\newtheorem{theorem}{Theorem}[section]

\newtheorem{proposition}[theorem]{Proposition}
\newtheorem{lemma}[theorem]{Lemma}
\newtheorem{corollary}[theorem]{Corollary}
\theoremstyle{definition}
\newtheorem*{claim*}{Claim}
\newtheorem{definition}[theorem]{Definition}

\newtheorem{remark}[theorem]{Remark}

\AtBeginDocument{\def\MR#1{}}

\begin{document}
\title{Derivation length and automorphism length of unital C*-algebras}
\author{Martino Lupini}
\address{Dipartimento di Matematica, Universit\`{a} di Bologna, Piazza di
Porta S. Donato, 5, 40126 Bologna,\ Italy}
\email{martino.lupini@unibo.it}
\urladdr{http://www.lupini.org/}
\thanks{The author was partially supported by the Starting Grant 101077154
\textquotedblleft Definable Algebraic Topology\textquotedblright\ from the
European Research Council, the Universit\'{e} Paris Cit\'{e}, the Gruppo
Nazionale per le Strutture Algebriche, Geometriche e le loro Applicazioni
(GNSAGA) of the Istituto Nazionale di Alta Matematica (INDAM), and the
University of Bologna. Part of this research was conducted during a visit of
the author to Chalmers University and the University of Gothenburg. The
hospitality of these institutions is gratefully acknowledged.}
\date{\today }

\begin{abstract}
This paper is a contribution to the study of the ordinal-valued invariants
of derivation and automorphism length. Akemann--Pedersen and Elliott proved
that a separable unital C*-algebra has derivation length $0$ if and only if
it has automorphism length at most $1/2$ if and only if it is a finite
direct sum of homogeneous C*-algebras and simple C*-algebras.
Kadison--Lance--Ringrose and Somerset proved that a separable unital
C*-algebra has derivation length at most $1$ if and only if it has
automorphism length at most $1$ if and only if its primitive spectrum has
finite connecting order.

In this paper, we prove a complete comparison between the two lengths:
either 
\begin{equation*}
\ell _{\mathrm{Aut}}\left( A\right) =\ell _{\mathrm{aut}}\left( A\right) 
\end{equation*}%
or, for some countable ordinal $\alpha $ other than $1$ or a limit ordinal, 
\begin{equation*}
\ell _{\mathrm{aut}}\left( A\right) =\alpha \quad \text{and}\quad \ell _{%
\mathrm{Aut}}\left( A\right) =\alpha +1/2\text{.}
\end{equation*}
\end{abstract}

\maketitle
\tableofcontents

\setcounter{tocdepth}{1}

\section{Introduction}

Building on previous works of\ Saint-Raymond \cite%
{saint-raymond_espaces_1975,saint-raymond_espaces_1976}, Solecki \cite%
{solecki_polish_1999} and Farah--Solecki \cite{farah_borel_2006}, Bergfalk,
Caputi, Casarosa, Codenotti, Panagiotopoulos, and the author \cite%
{bergfalk_definable_2024,bergfalk_definable_2024-1,lupini_looking_2024,lupini_complexity_2025,casarosa_homological_2026,bergfalk_definable_2026,lupini_applications_2025,caputi_general_2026}
have been applying tools from descriptive set theory to algebraic invariants 
\emph{explicitly presented }as quotients of Polish groups. Such \emph{groups
with a Polish cover }(whose importance had already been recognized by Moore 
\cite{moore_group_1976}, who called them \emph{pseudo-polonais groups}) are
groups of the form $G/N$ where $G$ is a Polish group, and $N$ is a normal
subgroup of $G$ that is also a Polish group with respect to a topology
(necessarily unique, possibly different from the subspace topology) making
the inclusion $N\rightarrow G$ continuous. Such groups admit a rich
structure of algebraic, topological, and complexity-theoretic\ nature,
including in particular the sequence of \emph{phantom subgroups}. In turn,
phantom subgroups allow one to define the \emph{phantom length}, which is a
countable ordinal (or a \emph{fractional}\ countable ordinal) measuring the
complexity of the group with a Polish cover, and essentially corresponding
to the least index of a phantom subgroup that becomes trivial.

Groups with a Polish cover arise naturally across several mathematical
subjects, including algebra (\textrm{Ext} and cohomology of groups),
topology (homology, cohomology, and $\mathrm{K}$-homology of spaces), and
operator algebras (Picard groups, cohomology groups). Their consideration as
groups with a Polish cover produces ordinal-valued notions of \emph{length},
which are new invariants whose study can be seen as a far-reaching
generalization of several classical lines of research. These include the
problem of determining when all the derivations of a C*-algebra are inner,
corresponding to the value $0$ for the \emph{derivation length }\cite%
{akemann_central_1979,elliott_some_1977}, or when the space of inner
derivations is norm-closed, corresponding to derivation length at most $1$,
as defined below, and hence to a value in $\left\{ 0,1\right\} $ \cite%
{kadison_derivations_1967,kadison_derivations_1967-1}. Likewise, the problem
of whether the group of inner automorphisms of a C*-algebra is \emph{%
norm-closed }within the group of all automorphisms is equivalent to
determining whether the \emph{automorphism length} is at most $1$ (Lemma \ref%
{Lemma:phantom-automorphism-1}).

Separable unital C*-algebras of derivation length $0$ were classified in 
\cite{akemann_central_1979,elliott_some_1977} as the finite direct sums of
homogeneous C*-algebras and simple C*-algebras. These are also precisely the
separable unital C*-algebras of automorphism length at most $1/2$. It was
shown in \cite%
{kadison_derivations_1967,kadison_derivations_1967-1,somerset_inner_1993}
that a separable unital C*-algebra has derivation length at most $1$ if and
only if it has automorphism length at most $1$ if and only if its primitive
spectrum has finite \emph{connecting order}.

In this paper, we establish a general \emph{local comparison} result\emph{\ }%
for TSI groups with a Polish cover in terms of \emph{pseudomorphisms},
allowing one to relate their phantom lengths. We apply this result to the
pseudomorphism from derivations to automorphisms given by the exponential
map to prove a complete comparison for every separable unital C*-algebra $A$%
: either 
\begin{equation*}
\ell _{\mathrm{Aut}}\left( A\right) =\ell _{\mathrm{aut}}\left( A\right) 
\end{equation*}%
or, for some countable ordinal $\alpha $ other than $1$ or a limit ordinal, 
\begin{equation*}
\ell _{\mathrm{aut}}\left( A\right) =\alpha \quad \text{and}\quad \ell _{%
\mathrm{Aut}}\left( A\right) =\alpha +1/2
\end{equation*}%
(Theorem \ref{Theorem:equality}). Thus the two lengths agree whenever the
derivation length is not an ordinal. Consequently, for every countable
ordinal $\alpha $, 
\begin{equation*}
\ell _{\mathrm{aut}}\left( A\right) \leq \alpha +1/2\Leftrightarrow \ell _{%
\mathrm{Aut}}\left( A\right) \leq \alpha +1/2\text{.}
\end{equation*}%
We also have the blockwise relation 
\begin{equation*}
\ell _{\mathrm{aut}}\left( A\right) \in \lbrack \alpha ,\alpha
+1)\Leftrightarrow \ell _{\mathrm{Aut}}\left( A\right) \in \lbrack \alpha
,\alpha +1)\text{.}
\end{equation*}%
Furthermore, neither the derivation length nor the automorphism length can
be a limit ordinal or $1/2+\varepsilon $, although these values can arise
for the automorphism length in the nonunital case, as described in work in
preparation \cite{bergfalk_definable_2026}.

This paper is divided into 9 sections, including this introduction. In
Section \ref{Section:phantom} we will recall the results from the theory of
groups (or, more generally, homogeneous spaces) with a Polish cover, to be
used in the rest of the paper. In Section \ref{Section:shani} we recall
Shani's Dichotomy Theorem, and we deduce from it the criterion for lower
bounds on complexity that will be applied in Section \ref{Section:dichotomy}%
. Section \ref{Section:TSI} introduces TSI groups with a Polish cover and
pseudomorphisms between them, and establishes the local comparison result
for their phantom subgroups (Proposition \ref{Proposition:local-comparison}%
), which is the main abstract tool of this paper. The theory of Banach--Lie
groups and Banach--Lie algebras is recalled in Section \ref%
{Section:Banach-Lie}. We will be exclusively concerned with the Banach--Lie
group of automorphisms of a C*-algebra, and its unitary group, together with
their Banach--Lie algebras of derivations and self-adjoint elements.
Derivation length and automorphism length of a separable unital C*-algebra
are introduced in Section \ref{Section:derivation} and Section \ref%
{Section:automorphism}, respectively, where they are also related to each
other. Section \ref{Section:dichotomy} presents the proof of the Dichotomy
Theorem for derivations and automorphisms of C*-algebras. Section \ref%
{Section:comparison} concludes with the comparison between derivation and
automorphism length.

\subsection*{Acknowledgments}

We would like to thank Shaun Allison, Jeffrey Bergfalk, Matteo Casarosa,
Luigi Caputi, Alessandro Codenotti, Eusebio Gardella, Michael Hartz,
Aristotelis Panagiotopoulos, and Marvin for many inspiring conversations. We
also thank Ilja Gogi\'{c} for his comments on a preliminary version of this
article, and for suggesting the cocycle argument in the proof of Proposition %
\ref{Proposition:from-derivation-to-automorphism}. ChatGPT by OpenAI
(version GPT-5.6 Sol) and Claude by Anthropic (version Opus 5) have been
used in the preparation of this manuscript.

\section{Phantom spaces\label{Section:phantom}}

In this section we recall some notions from the theory of \emph{homogeneous
spaces with a Polish cover} and \emph{Borel complexity theory}; see also 
\cite%
{bergfalk_definable_2024,bergfalk_definable_2024-1,lupini_looking_2024,casarosa_homological_2026,lupini_complexity_2025}%
.

\subsection{Ordinal notation}

In this section we introduce some notation about ordinals and
\textquotedblleft fractional ordinals\textquotedblright. We let $\omega $ be
the first \emph{infinite }ordinal, and $\omega _{1}$ be the first \emph{%
uncountable }ordinal. We let $\omega _{1}[1/2]$ be the set that contains,
for every countable ordinal $\alpha $, both $\alpha $ and $\alpha +1/2$
subject to the relations%
\begin{equation*}
\alpha <\alpha +\frac{1}{2}<\alpha +1
\end{equation*}%
and%
\begin{equation*}
\beta +(\alpha +\frac{1}{2})=\left( \beta +\alpha \right) +\frac{1}{2}\text{.%
}
\end{equation*}%
We also let $\omega _{1}^{\mathrm{Ph}}$ be the set that contains $\omega
_{1}[1/2]$ and, for every $\alpha <\omega _{1}$, the element $\alpha
+1/2+\varepsilon $ subject to the relations%
\begin{equation*}
\alpha +\frac{1}{2}<\alpha +\frac{1}{2}+\varepsilon <\alpha +1
\end{equation*}%
and%
\begin{equation*}
\beta +(\alpha +\frac{1}{2}+\varepsilon )=\left( \beta +\alpha \right) +%
\frac{1}{2}+\varepsilon \text{\quad and\quad }(\alpha +\frac{1}{2}%
)+\varepsilon =\alpha +\frac{1}{2}+\varepsilon \text{.}
\end{equation*}%
If $\alpha $ is a successor ordinal, we let $\alpha -1$ be its immediate
predecessor, and%
\begin{equation*}
\alpha -\frac{1}{2}=\left( \alpha -1\right) +\frac{1}{2}
\end{equation*}%
and%
\begin{equation*}
\alpha -\frac{1}{2}+\varepsilon =(\alpha -\frac{1}{2})+\varepsilon \text{.}
\end{equation*}%
We also set $\omega _{1}^{\mathbf{\Delta }}\subseteq \omega _{1}^{\mathrm{Ph}%
}$ to be the subset containing $\omega _{1}[1/2]$ together with $\lambda
+1/2+\varepsilon $ whenever $\lambda $ is either zero or a limit ordinal.

\subsection{Complexity classes}

Recall that a \emph{Polish space} is a second countable topological space
that admits a compatible complete metric. A \emph{Polish group} is a
topological group whose topology is Polish. A \emph{complexity class }$%
\Gamma $ is an assignment $X\mapsto \Gamma \left( X\right) $ of a collection 
$\Gamma \left( X\right) $ of Borel subsets of $X$ to each Polish space $X$,
in such a way that a continuous function $f:X\rightarrow Y$ induces a
function $\Gamma \left( Y\right) \rightarrow \Gamma \left( X\right) $ by
taking preimages. We will consider the complexity classes $\boldsymbol{\Pi }%
_{\alpha }^{0}$ and $\boldsymbol{\Sigma }_{\alpha }^{0}$ for $\alpha <\omega
_{1}$ as defined in \cite{kechris_classical_1995}; see also \cite%
{gao_invariant_2009}. For $\alpha <\omega _{1}$ we will also consider the
complexity class $D(\boldsymbol{\Pi }_{\alpha }^{0})$ comprising the sets
that can be written as the intersection of an element of $\boldsymbol{\Pi }%
_{\alpha }^{0}$ and an element of $\boldsymbol{\Sigma }_{\alpha }^{0}$.

The \emph{dual class }$\check{\Gamma}$ of a complexity class $\Gamma $
contains the \emph{complements }of the elements of $\Gamma $. A complexity
class is not \emph{self-dual} if it is different from its dual class. In
this case, for a Polish space $X$ and a Borel subset $A$, one says that $%
\Gamma $ is the complexity class of $A$ if $A\in \Gamma \left( X\right) $
and $A\notin \check{\Gamma}\left( X\right) $.

We define an increasing sequence $\left( \Gamma _{\alpha }\right) _{\alpha
\in \omega _{1}^{\mathrm{Ph}}}$ of complexity classes as follows: For $%
\lambda <\omega _{1}$ that is either zero or a limit ordinal, $n<\omega $,
and $i\in \left\{ 0,1/2,1/2+\varepsilon \right\} $, define:%
\begin{equation*}
\Gamma _{\lambda +n+i}:=\left\{ 
\begin{array}{ll}
\boldsymbol{\Pi }_{1+\lambda }^{0} & \text{for }n=0\text{ and }i=0\text{;}
\\ 
\boldsymbol{\Sigma }_{1+\lambda +1}^{0} & \text{for }n=0\text{ and }i=1/2%
\text{;} \\ 
D(\boldsymbol{\Pi }_{1+\lambda +1}^{0}) & \text{for }n=0\text{ and }%
i=1/2+\varepsilon \text{;} \\ 
D(\boldsymbol{\Pi }_{1+\lambda +n+1}^{0}) & \text{for }n\geq 1\text{ and }%
i\in \left\{ 1/2,1/2+\varepsilon \right\} \text{;} \\ 
\boldsymbol{\Pi }_{1+\lambda +n+1}^{0} & \text{for }n\geq 1\text{ and }i=0%
\text{;}%
\end{array}%
\right.
\end{equation*}%
For a Borel subset $A$ of a Polish space $X$, we define the \emph{Borel
class }of $A$ in $X$ to be the least $\alpha \in \omega _{1}^{\mathrm{Ph}}$
such that $A\in \Gamma _{\alpha }\left( X\right) $. Notice that for every
successor ordinal $\sigma $, the classes $\Gamma _{\sigma +1/2}$ and $\Gamma
_{\sigma +1/2+\varepsilon }$ are equal. Thus, by definition the Borel class
of $A$ in $X$ must be in $\omega _{1}^{\mathbf{\Delta }}$.

\subsection{Homogeneous spaces with a Polish cover}

Let $G$ be a Polish group. We call a subgroup $N$ of $G$ a \emph{Polish
subgroup} if it is a Polish group with respect to a (necessarily unique)
Polish group topology that makes the inclusion $N\rightarrow G$ continuous.
This is the notion traditionally called a \emph{Polishable subgroup} in the
literature; in particular, the topology on $N$ need not be the subspace
topology inherited from $G$. A \emph{homogeneous space with a Polish cover}
is a pointed space of the form $G/N$ where $G$ is a Polish group and $N$ is
a Polish subgroup of $G$, with distinguished point $\ast $ corresponding to
the trivial $N$-coset.\ A \emph{homogeneous subspace} with a Polish cover of 
$G/N$ is a subspace of the form $H/N$ where $H$ is a Polish subgroup of $G$
containing $N$.\ In particular, $\left\{ \ast \right\} =N/N$ is a
homogeneous subspace with a Polish cover. We define the Borel class of $H/N$
in $G/N$ to be the Borel class of $H$ in $G$. It is proved in \cite%
{lupini_complexity_2025} building on \cite%
{farah_borel_2006,solecki_polish_1999} that for every $\alpha \in \omega
_{1} $ there exists a \emph{least }homogeneous subspace with a Polish cover $%
\mathrm{Ph}^{\alpha }\left( G/N\right) =G_{\alpha }/N$ of $G/N$ of Borel
class at most $\alpha $.

We define the \emph{Borel length }$\ell _{\boldsymbol{\Delta }}\left(
G/N\right) \in \omega _{1}^{\mathbf{\Delta }}$ of $G/N$ to be the Borel
class of $\left\{ \ast \right\} $ in $G/N$ or, equivalently, the Borel class
of $N$ in $G$. One defines the \emph{phantom length }$\ell _{\mathrm{Ph}%
}\left( G/N\right) \in \omega _{1}^{\mathrm{Ph}}$ of $G/N$ as follows:

\begin{definition}
Let $G/N$ be a homogeneous space with a Polish cover. Set $\mathrm{Ph}%
^{\alpha }\left( G/N\right) :=G_{\alpha }/N$. For $\alpha <\omega _{1}$
define:

\begin{itemize}
\item $\ell _{\mathrm{Ph}}\left( G/N\right) \leq \alpha \Leftrightarrow 
\mathrm{Ph}^{\alpha }\left( G/N\right) =\left\{ \ast \right\} $;

\item $\ell _{\mathrm{Ph}}\left( G/N\right) \leq \alpha +1/2\Leftrightarrow
\left\{ \ast \right\} \in \boldsymbol{\Sigma }_{2}^{0}(\mathrm{Ph}^{\alpha
}\left( G/N\right) )$;

\item $\ell _{\mathrm{Ph}}\left( G/N\right) \leq \alpha +1/2+\varepsilon
\Leftrightarrow \left\{ \ast \right\} \in D(\boldsymbol{\Pi }_{2}^{0})(%
\mathrm{Ph}^{\alpha }\left( G/N\right) )$.
\end{itemize}
\end{definition}

If $N$ is open in $G$, then it is also closed, and both the phantom length
and the Borel length of $G/N$ are zero. If $N$ is not open, it follows from 
\cite[Theorem 6.1]{lupini_complexity_2025} that the phantom length of $G/N$
completely determines the Borel length of $\left\{ \ast \right\} $ in $G/N$.
Thus, in all cases:

\begin{theorem}
\label{Theorem:lupini-complexity}Let $G/N$ be a homogeneous space with a
Polish cover. Then%
\begin{equation*}
\Gamma _{\ell _{\boldsymbol{\Delta }}\left( G/N\right) }=\Gamma _{\ell _{%
\mathrm{Ph}}\left( G/N\right) }\text{;}
\end{equation*}%
in particular, $N\in \Gamma _{\ell _{\mathrm{Ph}}\left( G/N\right) }\left(
G\right) $. Consequently, since the sequence $\left( \Gamma _{\alpha
}\right) _{\alpha \in \omega _{1}^{\mathrm{Ph}}}$ is increasing, for every $%
\alpha \in \omega _{1}^{\mathrm{Ph}}$,%
\begin{equation*}
\ell _{\mathrm{Ph}}\left( G/N\right) \leq \alpha \Rightarrow \ell _{%
\boldsymbol{\Delta }}\left( G/N\right) \leq \alpha \text{,}
\end{equation*}%
while, for $\alpha \in \omega _{1}^{\mathbf{\Delta }}$,%
\begin{equation*}
\ell _{\boldsymbol{\Delta }}\left( G/N\right) \leq \alpha \Rightarrow \ell _{%
\mathrm{Ph}}\left( G/N\right) \leq \alpha
\end{equation*}%
\emph{provided that }$\alpha $\emph{\ is not of the form }$\sigma +1/2$\emph{%
\ for a successor ordinal }$\sigma $; in that remaining case one obtains
only $\ell _{\mathrm{Ph}}\left( G/N\right) \leq \sigma +1/2+\varepsilon $.
\end{theorem}

\begin{remark}
\label{Remark:Delta-exception}The exceptional case in Theorem \ref%
{Theorem:lupini-complexity} is genuine, and it is the only one. Indeed, $%
\Gamma $ is injective on $\omega _{1}^{\mathbf{\Delta }}$, and the only
coincidences in the sequence $\left( \Gamma _{\alpha }\right) _{\alpha \in
\omega _{1}^{\mathrm{Ph}}}$ are $\Gamma _{\sigma +1/2}=\Gamma _{\sigma
+1/2+\varepsilon }$ for $\sigma $ a successor ordinal. The phantom length
can take the value $\alpha +3/2+\varepsilon $ for a countable ordinal $%
\alpha $, which is of that form with $\sigma =\alpha +1$; the corresponding
Borel length is then $\alpha +3/2$. In particular, the phantom length is a
strictly finer invariant than the Borel length, and it need not belong to $%
\omega _{1}^{\mathbf{\Delta }}$.
\end{remark}

We will use the first assertion of Theorem \ref{Theorem:lupini-complexity}
in the following form: if $\ell _{\mathrm{Ph}}\left( G/N\right) \leq \alpha $
for some $\alpha \in \omega _{1}^{\mathrm{Ph}}$, then $N\in \Gamma _{\alpha
}\left( G\right) $. Notice that this does \emph{not} require $\alpha $ to
belong to $\omega _{1}^{\mathbf{\Delta }}$.

One defines a group with a Polish cover to be a homogeneous space with a
Polish cover $G/N$ where $N$ is furthermore normal in $G$. In this case, a
subgroup with a Polish cover is a homogeneous subspace with a Polish cover $%
H/N$ such that $H$ is normal in $G$.\ When $G/N$ is a group with a Polish
cover, $\mathrm{Ph}^{\alpha }\left( G/N\right) $ is a subgroup with a Polish
cover for every $\alpha <\omega _{1}$.

\subsection{Classes of Polish groups}

Let $G$ be a Polish group. Then $G$ admits a compatible left-invariant
metric \cite[Theorem 2.1.1]{gao_invariant_2009}, and a compatible complete
metric; see also \cite[Corollary 2.2.2]{gao_invariant_2009}. However, $G$
may not admit a compatible metric that is both complete and left-invariant.
Recall that a Polish group $G$ is:

\begin{itemize}
\item non-Archimedean\emph{\ }if every identity neighborhood of $G$ contains
an open subgroup \cite{gao_non-archimedean_2014};

\item CLI if it admits a compatible complete left-invariant metric \cite[%
Definition 2.2.5]{gao_invariant_2009};

\item TSI if it admits a compatible two-sided invariant metric \cite%
{gao_non-archimedean_2014}.
\end{itemize}

For example, the Polish group $S_{\infty }$ of permutations of $\mathbb{N}$
endowed with the topology of pointwise convergence, where $\mathbb{N}$ has
the discrete topology, is not CLI%
%TCIMACRO{\TeXButton{@}{\@}}%
%BeginExpansion
\@%
%EndExpansion
. Since a compatible two-sided invariant metric on a Polish group is
necessarily complete \cite[Corollary 1.2.2]{becker_descriptive_1996}, every 
\emph{abelian} Polish group is TSI, and every TSI\ Polish group is CLI%
%TCIMACRO{\TeXButton{@}{\@}}%
%BeginExpansion
\@%
%EndExpansion
. A Polish group is TSI if and only if it is SIN \cite{rosendal_global_2013}%
, i.e., it admits a basis of \emph{conjugation-invariant }neighborhoods of
the identity \cite{klee_invariant_1952}. A Polish group is both TSI and
non-Archimedean if and only if it admits a basis of neighborhoods of the
identity consisting of \emph{open normal subgroups} \cite%
{gao_non-archimedean_2014}.

Suppose that $X=G/N$ is a homogeneous space with a Polish cover and that $N$
is non-Archimedean. If $N$ is open in $G$, then it is also closed and the
phantom length of $X$ is zero. If $N$ is not open, it follows from \cite[%
Theorem 6.1]{lupini_complexity_2025} that the phantom length of $X$ belongs
to $\omega _{1}[1/2]$. Thus, in all cases the phantom length of $X$ belongs
to $\omega _{1}[1/2]$. In particular, the phantom length of such homogeneous
spaces with a Polish cover is invariant under Borel-definable
basepoint-preserving bijections.

\section{Dichotomy\label{Section:shani}}

In this section we recall Shani's Dichotomy Theorem from \cite[Theorem 1.1]%
{shani_generic_2024}. We use it to obtain a sufficient criterion to
establish a lower bound to the complexity of certain equivalence relations.

\subsection{$c_{0}\left( \mathbb{R}\right) $-orthogonality}

To simplify the statement and proofs of our main results, we introduce a new
notion concerning equivalence relations on Polish spaces. Denote by $c_{0}$
the separable Banach space of vanishing sequences of real numbers, endowed
with the supremum norm.\ We regard $c_{0}$ as a subspace of $\mathbb{R}^{%
\mathbb{N}}$. For a subset $X$ of $\mathbb{R}^{\mathbb{N}}$, we let $X/c_{0}$
be the quotient of $X$ by the restriction of the coset relation of $c_{0}$
on $\mathbb{R}^{\mathbb{N}}$ to $X$. For $\boldsymbol{t}\in \mathbb{R}^{%
\mathbb{N}}$ we define%
\begin{equation*}
\left\Vert \boldsymbol{t}\right\Vert _{\infty }:=\sup_{n}\left\vert
t_{n}\right\vert
\end{equation*}%
to be its supremum norm. A subset $Y$ of $\mathbb{R}^{\mathbb{N}}$ is:

\begin{itemize}
\item $\delta $-separated if, for every distinct $\boldsymbol{t},\boldsymbol{%
s}\in Y$, $\left\Vert \boldsymbol{t}-\boldsymbol{s}\right\Vert _{\infty
}\geq \delta $;

\item well-separated if it is $\delta $-separated for some $\delta >0$.
\end{itemize}

\begin{definition}
\label{Definition:c0-orthogonal}An equivalence relation $E$ on a Polish
space $X$ is $c_{0}$-\emph{orthogonal} if for every Borel function 
\begin{equation*}
f:\left( 0,1\right) ^{\mathbb{N}}\rightarrow X
\end{equation*}%
that induces a function 
\begin{equation*}
f_{0}:\left( 0,1\right) ^{\mathbb{N}}/c_{0}\rightarrow X/E\text{,}
\end{equation*}%
there exists an uncountable well-separated $Y\subseteq \left( 0,1\right) ^{%
\mathbb{N}}$ such that $f_{0}$ is constant on $Y/c_{0}$.
\end{definition}

Clearly, in the definition of $c_{0}$-orthogonal one can replace $\left(
0,1\right) $ with any other open interval.

\subsection{Tail equivalence and dichotomy}

We now establish a relation between the notion of $c_{0}$-orthogonality and
the relation of tail equivalence of countably many binary sequences. Recall
that $E_{0}$ denotes the relation on $2^{\mathbb{N}}$ of \emph{tail
equivalence }of binary sequences. Let $\left\langle \cdot ,\cdot
\right\rangle :\mathbb{N}\times \mathbb{N}\rightarrow \mathbb{N}$ be a
bijection such that $\left\langle n,m\right\rangle \leq \left\langle
n^{\prime },m^{\prime }\right\rangle $ for all $n,m,n^{\prime },m^{\prime
}\in \mathbb{N}$ such that $n\leq n^{\prime }$ and $m\leq m^{\prime }$.
Define $\varphi :\left( 2^{\mathbb{N}}\right) ^{\mathbb{N}}\rightarrow 
\mathbb{R}^{\mathbb{N}}$ by setting%
\begin{equation*}
\varphi \left( x\right) _{\left\langle n,m\right\rangle }=2^{-n}(x_{n})_{m}%
\text{.}
\end{equation*}%
Then $\varphi $ is continuous and, by \cite[Lemma 8.5.3]{gao_invariant_2009}%
, is a reduction of $E_{0}^{\mathbb{N}}$ to the coset relation of $c_{0}$;
explicitly, 
\begin{equation*}
xE_{0}^{\mathbb{N}}y\Leftrightarrow \varphi \left( x\right) -\varphi \left(
y\right) \in c_{0}\text{.}
\end{equation*}%
In particular, 
\begin{equation*}
\varphi \left( x\right) \in c_{0}\Leftrightarrow \forall n\exists m\forall
k\geq m\text{, }x_{n,k}=0\text{;}
\end{equation*}%
thus, $\varphi $ induces an injective function%
\begin{equation*}
\varphi _{0}:\left( 2^{\mathbb{N}}\right) ^{\mathbb{N}}/E_{0}^{\mathbb{N}%
}\rightarrow \lbrack 0,1]^{\mathbb{N}}/c_{0}\text{.}
\end{equation*}%
We recall the statement of Shani's Dichotomy Theorem for the relation $%
E_{0}^{\mathbb{N}}$ on $\left( 2^{\mathbb{N}}\right) ^{\mathbb{N}}$ \cite[%
Theorem 1.1]{shani_generic_2024}.

\begin{theorem}[Shani]
\label{Theorem:Shani}Let $F$ be an analytic equivalence relation on a Polish
space $Y$. Suppose that the relation $E_{0}^{\mathbb{N}}$ of tail
equivalence of countably many binary sequences is not Borel reducible to $F$%
. Let 
\begin{equation*}
f:\left( 2^{\mathbb{N}}\right) ^{\mathbb{N}}\rightarrow Y
\end{equation*}%
be a Borel function that induces a function 
\begin{equation*}
f_{0}:\left( 2^{\mathbb{N}}\right) ^{\mathbb{N}}/E_{0}^{\mathbb{N}%
}\rightarrow Y/F\text{.}
\end{equation*}%
Then there exist $k\in \omega $, a comeager $X\subseteq \left( 2^{\mathbb{N}%
}\right) ^{\mathbb{N}}$, and a Borel function $g:(2^{\mathbb{N}%
})^{k}\rightarrow Y$ that induces a function 
\begin{equation*}
g_{0}:\left( 2^{\mathbb{N}}\right) ^{k}/E_{0}^{k}\rightarrow Y/F
\end{equation*}%
such that 
\begin{equation*}
f_{0}|_{X/E_{0}^{\mathbb{N}}}=(g_{0}\circ \pi _{k})|_{X/E_{0}^{\mathbb{N}}}%
\text{,}
\end{equation*}%
where 
\begin{equation*}
\pi _{k}:\left( 2^{\mathbb{N}}\right) ^{\mathbb{N}}/E_{0}^{\mathbb{N}%
}\rightarrow \left( 2^{\mathbb{N}}\right) ^{k}/E_{0}^{k}
\end{equation*}%
is the projection on the first $k$ coordinates.
\end{theorem}

The formulation above is the totalized form of Shani's theorem. Indeed, \cite%
[Theorem 1.1]{shani_generic_2024} initially provides a Borel homomorphism $%
\widetilde{g}:C\rightarrow Y$ on a comeager Borel subset $C$ of $\left( 2^{%
\mathbb{N}}\right) ^{k}$. Let 
\begin{equation*}
\Gamma _{k}:=\biggl(\bigoplus_{\mathbb{N}}\mathbb{Z}/2\mathbb{Z}\biggr)^{k}
\end{equation*}%
act on $\left( 2^{\mathbb{N}}\right) ^{k}$ by coordinatewise addition. Its
orbit equivalence relation is $E_{0}^{k}$. The set 
\begin{equation*}
C_{\ast}:=\bigcap_{\gamma \in \Gamma _{k}}\gamma C
\end{equation*}%
is Borel, comeager, and $E_{0}^{k}$-invariant. Restrict $\widetilde{g}$ to $%
C_{\ast}$ and extend it to all of $\left( 2^{\mathbb{N}}\right) ^{k}$ by
assigning a fixed value $y_{0}\in Y$ on the complement of $C_{\ast}$. Since $%
C_{\ast}$ is $E_{0}^{k}$-invariant, the resulting map is a total Borel
homomorphism. Finally, intersecting the comeager factorization set supplied
by Shani's theorem with the inverse image of $C_{\ast}$ under the coordinate
projection onto the first $k$ coordinates preserves comeagerness and the
asserted factorization. This yields precisely the formulation stated above.

We isolate the following consequence of Shani's Dichotomy Theorem.

\begin{corollary}
\label{Corollary:Shani}Let $F$ be an analytic equivalence relation on a
Polish space $Y$. If $F$ is not $c_{0}$-orthogonal, then the relation $%
E_{0}^{\mathbb{N}}$ is Borel reducible to $F$.
\end{corollary}

\begin{proof}
We prove the contrapositive. Suppose that the relation $E_{0}^{\mathbb{N}}$
of tail equivalence of countably many binary sequences is \emph{not Borel
reducible} to $F$. We verify that $F$ is $c_{0}$-orthogonal. Let 
\begin{equation*}
h:\left( -2,2\right) ^{\mathbb{N}}\rightarrow Y
\end{equation*}%
be a Borel function that induces a function 
\begin{equation*}
h_{0}:\left( -2,2\right) ^{\mathbb{N}}/c_{0}\rightarrow Y/F\text{.}
\end{equation*}%
Consider the Borel function $h\circ \varphi :\left( 2^{\mathbb{N}}\right) ^{%
\mathbb{N}}\rightarrow Y$. By Shani's Dichotomy Theorem, there exist $k\in
\omega $, a comeager subset $X\subseteq \left( 2^{\mathbb{N}}\right) ^{%
\mathbb{N}}$ and a Borel function 
\begin{equation*}
g:\left( 2^{\mathbb{N}}\right) ^{k}\rightarrow Y
\end{equation*}%
that induces a function 
\begin{equation*}
g_{0}:\left( 2^{\mathbb{N}}\right) ^{k}/E_{0}^{k}\rightarrow Y/F
\end{equation*}%
such that 
\begin{equation*}
(h_{0}\circ \varphi _{0})|_{X/E_{0}^{\mathbb{N}}}=(g_{0}\circ \pi
_{k})|_{X/E_{0}^{\mathbb{N}}}\text{.}
\end{equation*}%
By the Kuratowski--Ulam theorem \cite[Theorem 8.41]{kechris_classical_1995},
there exist $x_{0}\in \left( 2^{\mathbb{N}}\right) ^{k}$ and a comeager
subset $Z\subseteq 2^{\mathbb{N}}$ such that for every $z\in Z$ there exists 
$w_{z}\in X$ such that $w_{z}$ has as initial segment the concatenation $%
x_{0}\smallfrown z$. Consider%
\begin{equation*}
X_{0}:=\left\{ w_{z}:z\in Z\right\} \text{.}
\end{equation*}%
Observe that, for distinct $w,w^{\prime }\in X_{0}$,%
\begin{equation*}
\left\Vert \varphi \left( w\right) -\varphi \left( w^{\prime }\right)
\right\Vert _{\infty }\geq 2^{-(k+1)}\text{.}
\end{equation*}%
Thus, we have that%
\begin{equation*}
\left\{ \varphi \left( w\right) :w\in X_{0}\right\}
\end{equation*}%
is an uncountable well-separated subset of $\left( -2,2\right) ^{\mathbb{N}}$%
. Furthermore, for $x,x^{\prime }\in X_{0}/E_{0}^{\mathbb{N}}$ we have that%
\begin{equation*}
\left( h_{0}\circ \varphi _{0}\right) \left( x\right) =\left( g_{0}\circ \pi
_{k}\right) \left( x\right) =\left( g_{0}\circ \pi _{k}\right) \left(
x^{\prime }\right) =\left( h_{0}\circ \varphi _{0}\right) \left( x^{\prime
}\right) \text{.}
\end{equation*}%
This concludes the proof.
\end{proof}

\subsection{A lemma on potential complexity}

We recall a lemma on the potential complexity of coset equivalence
relations, whose proof is the same as \cite[Theorem 3.3]%
{lupini_complexity_2025}.

\begin{lemma}
\label{Lemma:complexity-class}Suppose that $G$ is a Polish group, and $H$ is
a Polish subgroup of $G$. If $U\subseteq G$ is a nonempty open subset,
define $E_{H}^{U}$ to be the restriction to $U$ of the coset relation $%
E_{H}^{G}$ of $H$ in $G$. Then:

\begin{enumerate}
\item $E_{H}^{U}$ is potentially $\boldsymbol{\Pi }_{2}^{0}$ if and only if $%
H$ is closed in $G$;

\item $E_{H}^{U}$ is potentially $\boldsymbol{\Sigma }_{2}^{0}$ if and only
if $H$ is $D(\boldsymbol{\Pi }_{2}^{0})$ in $G$;

\item for every $\alpha \in \omega _{1}^{\mathrm{Ph}}$ \emph{with }$\alpha
\neq 1/2$, $E_{H}^{U}$ is potentially $\Gamma _{\alpha }$ if and only if $%
G/H $ has phantom length at most $\alpha $;

\item $E_{H}^{U}$ is potentially $\boldsymbol{\Sigma }_{1}^{0}$ if and only
if $H$ is open in $G$.
\end{enumerate}
\end{lemma}

\begin{proof}
(1) If $E_{H}^{U}$ is potentially $\boldsymbol{\Pi }_{2}^{0}$, then $%
E_{H}^{G}$ is potentially $\boldsymbol{\Pi }_{2}^{0}$, and the conclusion
follows from \cite[Theorem 3.3(1)]{lupini_complexity_2025}. Conversely, if $%
H $ is closed in $G$, then $E_{H}^{G}$ is closed in $G\times G$, and hence $%
E_{H}^{U}$ is $\boldsymbol{\Pi }_{1}^{0}$.

(2) The forward implication follows from \cite[Proposition 3.1]%
{lupini_complexity_2025}, and the converse implication follows from \cite[%
Theorem 3.3(2)]{lupini_complexity_2025}.

(3) This is proved as in \cite[Theorem 3.3(3)]{lupini_complexity_2025} using 
\cite[Proposition 3.1]{lupini_complexity_2025}. The hypothesis $\alpha \neq
1/2$ is necessary, and it is the only value of $\omega _{1}^{\mathrm{Ph}}$
that has to be excluded. Indeed, $\Gamma _{1/2}=\boldsymbol{\Sigma }_{2}^{0}$%
, and $\boldsymbol{\Sigma }_{2}^{0}$ is precisely the complexity class that
is excluded in \cite[Theorem 3.3(3)]{lupini_complexity_2025} and treated
separately in \cite[Theorem 3.3(2)]{lupini_complexity_2025}; accordingly,
the correct statement at $\alpha =1/2$ is item (2) above, and not the
assertion that $H$ be $\boldsymbol{\Sigma }_{2}^{0}$ in $G$. For every other 
$\alpha \in \omega _{1}^{\mathrm{Ph}}$, $\Gamma _{\alpha }$ is either $%
\boldsymbol{\Pi }_{\beta }^{0}$, or $D(\boldsymbol{\Pi }_{\beta }^{0})$, or $%
\boldsymbol{\Sigma }_{\beta }^{0}$ for some $\beta \neq 2$.

(4) If $E_{H}^{U}$ is potentially $\boldsymbol{\Sigma }_{1}^{0}$ then $%
E_{H}^{U}$ has countably many classes \cite[Lemma 12.5.2]{gao_invariant_2009}%
. Thus, $E_{H}^{G}$ has countably many classes, and the conclusion follows
from \cite[Theorem 3.3(3)]{lupini_complexity_2025}. Conversely, if $H$ is
open in $G$, then $E_{H}^{G}$ is open in $G\times G$, and hence $E_{H}^{U}$
is open in $U\times U$.
\end{proof}

\subsection{Homomorphisms\label{Subsection:homomorphisms}}

Suppose that $G/N$ and $H/M$ are groups with a Polish cover. Then a \emph{%
Borel-definable }function $f:G/N\rightarrow H/M$ is a function that admits a 
\emph{Borel lift }$\varphi :G\rightarrow H$. By definition, this means that $%
\varphi $ is a Borel function such that, for every $g\in G$, $\varphi
(g)M=f(gN)$. It follows from \cite[Proposition 3.11]%
{bergfalk_definable_2024-1} that if $f$ is a\emph{\ }Borel-definable
bijection, then its inverse is also Borel-definable; see also \cite[Lemma 3.8%
]{kechris_borel_2016}.

It follows from Lemma \ref{Lemma:complexity-class} that if $f:G/N\rightarrow
H/M$ is an \emph{injective }Borel-definable function, and $\alpha <\omega
_{1}$, then:

\begin{itemize}
\item $\ell _{\mathrm{Ph}}\left( H/M\right) \leq \alpha \Rightarrow \ell _{%
\mathrm{Ph}}\left( G/N\right) \leq \alpha $;

\item $\ell _{\mathrm{Ph}}\left( H/M\right) \leq \alpha +\frac{1}{2}%
+\varepsilon \Rightarrow \ell _{\mathrm{Ph}}\left( G/N\right) \leq \alpha +%
\frac{1}{2}+\varepsilon $.
\end{itemize}

Thus, if $f$ is a \emph{bijective }Borel-definable homomorphism, and $\alpha
<\omega _{1}$, then:

\begin{itemize}
\item $\ell _{\mathrm{Ph}}\left( H/M\right) =\alpha \Leftrightarrow \ell _{%
\mathrm{Ph}}\left( G/N\right) =\alpha $;

\item $\ell _{\mathrm{Ph}}\left( H/M\right) \in \left\{ \alpha +1/2,\alpha
+1/2+\varepsilon \right\} \Leftrightarrow \ell _{\mathrm{Ph}}\left(
G/N\right) \in \left\{ \alpha +1/2,\alpha +1/2+\varepsilon \right\} $.
\end{itemize}

We also say that $f:G/N\rightarrow H/M$ is an \emph{airomorphism }if it
admits a lift $\varphi :G\rightarrow H$ that is a continuous group
homomorphism. If there exists an injective airomorphism $f:G/N\rightarrow
H/M $ then by functoriality of the phantom subspaces one has $\ell _{\mathrm{%
Ph}}\left( G/N\right) \leq \ell _{\mathrm{Ph}}\left( H/M\right) $. An
iso-airomorphism is a bijective airomorphism whose inverse is also an
airomorphism. Thus, if there exists an iso-airomorphism $G/N\rightarrow H/M$
then $\ell _{\mathrm{Ph}}\left( G/N\right) =\ell _{\mathrm{Ph}}\left(
H/M\right) $.

\section{TSI groups with a Polish cover\label{Section:TSI}}

In this section we introduce \emph{TSI groups with a Polish cover}. We
consider a natural notion of \emph{pseudomorphism} between such groups,
which we use to obtain a comparison criterion for their phantom lengths.

\subsection{Metrics}

In a Polish group $G$, one can define the canonical \emph{left}, \emph{right}%
, and \emph{two-sided }uniformities \cite{rosendal_global_2013}. A Polish
group $G$ is TSI if and only if these three uniformities coincide, in which
case we refer to it as \emph{the uniformity} of $G$. For a TSI\ Polish group 
$G$, we endow $G^{n}$ with the product uniformity, and a subset $X\subseteq
G^{n}$ with the induced uniformity. A function $f$ between uniform spaces is
a \emph{uniform isomorphism }if it is uniformly continuous with uniformly
continuous inverse.\ For a group $G$ with a distinguished two-sided
invariant metric $d$ we define, for $\varepsilon >0$,%
\begin{equation*}
\mathrm{\mathrm{Ball}}_{\varepsilon }\left( G\right) :=\left\{ g\in
G:d\left( g,1\right) <\varepsilon \right\} \text{.}
\end{equation*}

\begin{definition}
\label{Definition:TSI}A group with a Polish cover $G/N$ is \emph{TSI }if
there exist compatible two-sided invariant metrics $d_{\infty }$ and $d_{0}$
on $N$ and $G$, respectively, such that $d_{0}\leq d_{\infty }\leq 1$ and
the canonical action $G\curvearrowright \left( N,d_{\infty }\right) $ by
conjugation is isometric.
\end{definition}

Suppose that $G/N$ is a group with a TSI Polish cover. Write $\mathrm{Ph}%
^{\alpha }\left( G/N\right) =G_{\alpha }/N$ for $\alpha <\omega _{1}$. Let $%
d_{\infty }$ be a compatible two-sided invariant metric on $N$ and $d_{0}$
be a compatible two-sided invariant metric on $G$, with $d_{0}\leq d_{\infty
}\leq 1$ as in Definition \ref{Definition:TSI}. We define recursively
two-sided invariant metrics $d_{\alpha }$ on $G_{\alpha }$ for $\alpha
<\omega _{1}$ by setting%
\begin{equation*}
d_{\alpha +1}\left( x,1\right) :=\inf R_{\alpha }\left( x\right) \text{,}
\end{equation*}%
where $R_{\alpha }\left( x\right) $ is the set of those $r>0$ for which
there exists a sequence $\left( w_{n}\right) $ in $N$ such that:

\begin{itemize}
\item $d_{\infty }\left( w_{n},1\right) \leq r$ for every $n\in \mathbb{N}$,
and

\item $\lim_{k}\ d_{\alpha }\left( w_{k},x\right) =0$.
\end{itemize}

The set $R_{\alpha }\left( x\right) $ is upward closed, so that $d_{\alpha
+1}\left( x,1\right) \leq r$ whenever $r\in R_{\alpha }\left( x\right) $,
and $r\in R_{\alpha }\left( x\right) $ for every $r>d_{\alpha +1}\left(
x,1\right) $. (We do \emph{not} claim that $R_{\alpha }\left( x\right) $
contains its infimum, and we will only use the two implications just stated.)

For $\lambda $ a limit ordinal, define%
\begin{equation*}
d_{\lambda }=\sum_{n\in \mathbb{N}}2^{-n-1}d_{\lambda _{n}}\text{,}
\end{equation*}%
where $\left( \lambda _{n}\right) $ is an increasing sequence of successor
ordinals that is cofinal in $\lambda $. The weights are normalized so that $%
\sum_{n\in \mathbb{N}}2^{-n-1}=1$, whence $d_{\lambda }\leq d_{\infty }\leq
1 $ on $N$, as required by Definition \ref{Definition:TSI}.

\begin{lemma}
\label{Lemma:metric}Let $G/N$ be a TSI group with a Polish cover, as
witnessed by metrics $d_{\infty }$ and $d_{0}$. Set $\mathrm{Ph}^{\alpha
}\left( G/N\right) =G_{\alpha }/N$ and define metrics $d_{\alpha }$ for $%
\alpha <\omega _{1}$ as above. Then $\mathrm{Ph}^{\alpha }\left( G/N\right) $
is a \emph{normal }subgroup of $G/N$, and $d_{\infty }$ and $d_{\alpha }$
witness that $\mathrm{Ph}^{\alpha }\left( G/N\right) $ is a TSI group with a
Polish cover.
\end{lemma}

\begin{proof}
Consider initially the case when $\alpha =1$. We need to verify that, for $%
x,x^{\prime }\in G_{1}$ and $g\in G$:

\begin{itemize}
\item $gxg^{-1}\in G_{1}$;

\item $d_{1}(gxg^{-1},1)\leq d_{1}\left( x,1\right) $;

\item $d_{1}(x^{-1},1)\leq d_{1}\left( x,1\right) $;

\item $d_{1}\left( xx^{\prime },1\right) \leq d_{1}\left( x,1\right)
+d_{1}\left( x^{\prime },1\right) $.
\end{itemize}

Since $x\in G_{1}$, for every $\varepsilon >0$ there exist $w\in N$ and a
sequence $\left( z_{n}\right) $ in $\mathrm{\mathrm{Ball}}_{\varepsilon
}\left( N\right) $ such that 
\begin{equation*}
\lim_{n}\ d_{0}\left( z_{n},wx\right) =0\text{.}
\end{equation*}%
Then $gwg^{-1}\in N$ and $\left( gz_{n}g^{-1}\right) $ is a sequence in $%
\mathrm{\mathrm{Ball}}_{\varepsilon }\left( N\right) $ with%
\begin{equation*}
\lim_{n}d_{0}\left( gz_{n}g^{-1},gwg^{-1}gxg^{-1}\right) =0\text{.}
\end{equation*}%
As this holds for every $\varepsilon >0$, $gxg^{-1}\in G_{1}$.

Let $r\in R_{0}(x)$, and let $\left( w_{n}\right) $ be a sequence in $N$
witnessing this. Since $d_{\infty }$ and $d_{0}$ are two-sided invariant and 
$N$ is normal in $G$ and the action $G\curvearrowright \left( N,d_{\infty
}\right) $ is isometric, $(gw_{n}g^{-1})_{n\in \mathbb{N}}$ is a sequence in 
$N$ witnessing that $r\in R_{0}(gxg^{-1})$, and $(w_{n}^{-1})$ is a sequence
in $N$ witnessing that $r\in R_{0}(x^{-1})$. Taking infima gives the second
and third assertions above. Applying these inequalities with $g^{-1}$ and $%
x^{-1}$, respectively, gives the reverse inequalities as well.

Similarly, let $r\in R_{0}(x)$ and $r^{\prime }\in R_{0}(x^{\prime })$, and
choose corresponding witnessing sequences $\left( w_{n}\right) $ and $\left(
w_{n}^{\prime }\right) $. Then $\left( w_{n}w_{n}^{\prime }\right) $ is a
sequence in $N$ witnessing that $r+r^{\prime }\in R_{0}(xx^{\prime })$.
Taking the infimum over $r$ and $r^{\prime }$ proves the fourth assertion.

For $r>0$, set%
\begin{equation*}
V_{r}:=\overline{\mathrm{\mathrm{Ball}}_{r}\left( N\right) }^{G_{0}}\text{.}
\end{equation*}%
For $x\in G_{1}$ and $0<s<r$, the definition of $d_{1}$ and the two
implications preceding the lemma give%
\begin{equation*}
\left\{ y\in G_{1}:d_{1}\left( x,y\right) <r\right\} \subseteq \left(
xV_{r}\cap V_{r}x\right) \cap G_{1}
\end{equation*}%
and%
\begin{equation*}
\left( xV_{s}\cap V_{s}x\right) \cap G_{1}\subseteq \left\{ y\in
G_{1}:d_{1}\left( x,y\right) <r\right\} \text{.}
\end{equation*}%
The sets $\left( xV_{r}\cap V_{r}x\right) \cap G_{1}$ form a neighborhood
basis for the canonical Polish group topology of $G_{1}$ \cite[Section 2]%
{solecki_polish_1999}. Therefore the two displayed inclusions show that $%
d_{1}$ is compatible with this topology. The case for arbitrary $\alpha $
follows by induction, noticing that:

\begin{enumerate}
\item $\mathrm{Ph}^{\alpha +1}\left( G/N\right) =\mathrm{Ph}^{1}(\mathrm{Ph}%
^{\alpha }(G/N))$, and

\item for $\lambda $ limit, $G_{\lambda }$ is the intersection of $G_{\beta
} $ for $\beta <\lambda $ endowed with the projective-limit topology,
\end{enumerate}

as shown in \cite[Section 2]{solecki_polish_1999}.
\end{proof}

\subsection{Pseudomorphisms}

We prove that in the context of TSI groups with a Polish cover, locally
defined uniformly continuous functions must map phantom subgroups to phantom
subgroups. In the following definitions, we assume that $G/N$ and $H/M$ are
TSI groups with a Polish cover.

\begin{definition}
\label{Definition:local-uniformly-continuous}A local uniformly continuous
map $G/N\rightarrow H/M$ is given by:

\begin{itemize}
\item an open symmetric identity neighborhood $W\subseteq G$;

\item an open identity neighborhood $U\subseteq N$;

\item a $\left( G,H\right) $-uniformly continuous function $\varphi
:W\rightarrow H$;
\end{itemize}

such that, setting $W_{\infty }:=W\cap N$:

\begin{enumerate}
\item $\varphi \left( 1\right) =1$;

\item $\varphi (x^{-1})=\varphi \left( x\right) ^{-1}$;

\item $\varphi $ maps $W_{\infty }$ to $M$;

\item $U\subseteq W_{\infty }$;

\item $\varphi |_{U}$ is $\left( N,M\right) $-uniformly continuous.
\end{enumerate}
\end{definition}

\begin{definition}
\label{Definition:pseudomorphism}A \emph{pseudomorphism }$G/N\rightarrow H/M$
is given by a local uniformly continuous map $\left( W,U,\varphi \right)
:G/N\rightarrow H/M$ for which there exists a function $\sigma :D\rightarrow
M$ where%
\begin{equation*}
D:=\left\{ \left( z,x\right) \in N\times W:zx\in W\right\}
\end{equation*}%
such that, setting%
\begin{equation*}
D_{U}:=\left\{ \left( z,x\right) \in D:z\in U\right\} \text{:}
\end{equation*}

\begin{itemize}
\item $\varphi \left( zx\right) =\sigma \left( z,x\right) \varphi \left(
x\right) $ for $\left( z,x\right) \in D$;

\item $\sigma |_{D_{U}}$ is $\left( N,M\right) $-uniformly continuous in the
first variable, with modulus independent of the second variable.
\end{itemize}
\end{definition}

\begin{proposition}
\label{Proposition:local-comparison}Let $G/N$ and $H/M$ be\ TSI groups with
a Polish cover. For $\alpha <\omega _{1}$ write $\mathrm{Ph}^{\alpha }\left(
G/N\right) =G_{\alpha }/N$ and $\mathrm{Ph}^{\alpha }\left( H/M\right)
=H_{\alpha }/M$. Suppose that $\left( W,U,\varphi \right) :G/N\rightarrow
H/M $ is a pseudomorphism. For $\alpha <\omega _{1}$, define 
\begin{equation*}
W_{\alpha }:=W\cap G_{\alpha }\text{.}
\end{equation*}%
Then, for every $\alpha <\omega _{1}$:

\begin{enumerate}
\item $\varphi $ maps $W_{\alpha }$ to $H_{\alpha }$;

\item $\varphi |_{W_{\alpha }}$ is $\left( G_{\alpha },H_{\alpha }\right) $%
-uniformly continuous.
\end{enumerate}
\end{proposition}

\begin{proof}
(1) and (2): Let us fix compatible two-sided invariant metrics $d_{0}$ and $%
d_{\infty }$ on $G$ and $N$ witnessing that $G/N$ is TSI, and let $d_{\alpha
}$ be the corresponding metric on $G_{\alpha }$ defined as in Lemma \ref%
{Lemma:metric}. Let us do the same for $H$.

For $\alpha =0$ we have that $G_{0}$ and $H_{0}$ are the closures of $N$ in $%
G$ and of $M$ in $H$, respectively, endowed with the induced metrics.
Suppose that $x\in W_{0}=W\cap G_{0}$. Then $x$ is the limit of a sequence $%
\left( w_{k}\right) $ in $N$ and, $W$ being open, we can assume that $%
w_{k}\in W\cap N=W_{\infty }$ for every $k\in \mathbb{N}$. Hence $\varphi
\left( w_{k}\right) \in M$ for every $k\in \mathbb{N}$ and, by continuity of 
$\varphi $,%
\begin{equation*}
\varphi \left( x\right) =\lim_{k}\varphi \left( w_{k}\right) \in H_{0}\text{,%
}
\end{equation*}%
which proves (1). Assertion (2) for $\alpha =0$ is immediate from the $%
\left( G,H\right) $-uniform continuity of $\varphi $.

We now consider the case $\alpha =1$. Suppose that $x\in W\cap G_{1}$. We
claim that $\varphi \left( x\right) \in H_{1}$. Fix $\varepsilon >0$. By the
uniform continuity of $\varphi |_{U}:U\rightarrow M$, choose $\delta >0$
such that $\varphi \left( \mathrm{\mathrm{Ball}}_{\delta }\left( N\right)
\right) \subseteq \mathrm{\mathrm{Ball}}_{\varepsilon }(M)$, $\mathrm{%
\mathrm{Ball}}_{\delta }\left( N\right) \subseteq U$, and the closed $d_{0}$%
-ball of radius $\delta $ around the identity is contained in $W$. There
exist $z\in N $ and a sequence $\left( w_{k}\right) $ in $\mathrm{\mathrm{%
Ball}}_{\delta }\left( N\right) $ such that $\lim_{k}\ d_{0}\left(
zx,w_{k}\right) =0$. It follows that $zx\in W$, so that $\left( z,x\right)
\in D$. Then:

\begin{itemize}
\item $\varphi \left( zx\right) =\sigma \left( z,x\right) \varphi \left(
x\right) $;

\item $\sigma \left( z,x\right) \in M$;

\item $\varphi \left( w_{k}\right) \in \mathrm{\mathrm{Ball}}_{\varepsilon
}(M)$;

\item $\lim_{k}\ d_{0}(\varphi \left( w_{k}\right) ,\sigma \left( z,x\right)
\varphi \left( x\right) )=0$.
\end{itemize}

This, together with the same argument applied to $x^{-1}$, shows that $%
\varphi \left( x\right) \in H_{1}$.

We now show that $\varphi |_{W_{1}}:W_{1}\rightarrow H_{1}$ is $\left(
G_{1},H_{1}\right) $-uniformly continuous. Suppose that $\varepsilon >0$.
Notice first that $\sigma \left( 1,x\right) =1$ for every $x\in W$. By the $%
\left( N,M\right) $-uniform continuity of $\sigma $ in the first variable,
with modulus independent of the second variable, we can choose $\eta >0$
such that the closed $d_{\infty }$-ball of radius $\eta $ around the
identity is contained in $U$ and%
\begin{equation*}
d_{\infty }\left( z,1\right) \leq \eta \Rightarrow \sigma \left( z,x\right)
\in \mathrm{\mathrm{Ball}}_{\varepsilon }(M)
\end{equation*}%
whenever $x,zx\in W$. Set $\delta :=\eta /2$. Suppose that $x,x^{\prime }\in
W_{1}$ satisfy $d_{1}\left( x,x^{\prime }\right) \leq \delta $. Since $\eta
>d_{1}\left( x,x^{\prime }\right) $, we have $\eta \in R_{0}(x^{\prime
}x^{-1})$. Hence there exists a sequence $\left( w_{n}\right) $ in $N$ such
that $d_{\infty }\left( w_{n},1\right) \leq \eta $ for every $n$ and, after
discarding finitely many terms if necessary:

\begin{itemize}
\item $\lim_{k}\ d_{0}(w_{k}x,x^{\prime })=0$;

\item $w_{n}x\in W$ and hence $\left( w_{n},x\right) \in D_{U}$ for every $%
n\in \mathbb{N}$.
\end{itemize}

Then, for every $n\in \mathbb{N}$,%
\begin{equation*}
\varphi \left( w_{n}x\right) =\sigma \left( w_{n},x\right) \varphi \left(
x\right)
\end{equation*}%
where%
\begin{equation*}
\sigma \left( w_{n},x\right) \in \mathrm{\mathrm{Ball}}_{\varepsilon }(M)
\end{equation*}%
and%
\begin{equation*}
\lim_{k}\ d_{0}(\sigma \left( w_{k},x\right) \varphi \left( x\right)
,\varphi \left( x^{\prime }\right) )=\lim_{k}\ d_{0}(\varphi \left(
w_{k}x\right) ,\varphi \left( x^{\prime }\right) )=0\text{.}
\end{equation*}%
This, together with the same argument for $x^{-1}$, shows that%
\begin{equation*}
d_{1}\left( \varphi \left( x\right) ,\varphi \left( x^{\prime }\right)
\right) \leq \varepsilon \text{.}
\end{equation*}%
This concludes the proof for $\alpha =1$.

We now proceed by transfinite induction. Suppose that (1) and (2) hold at
some $\alpha <\omega _{1}$. We claim that%
\begin{equation*}
\left( W_{\alpha },U,\varphi |_{W_{\alpha }}\right) :G_{\alpha
}/N\rightarrow H_{\alpha }/M
\end{equation*}%
is a pseudomorphism between TSI groups with a Polish cover. Indeed, $%
G_{\alpha }/N$ and $H_{\alpha }/M$ are TSI groups with a Polish cover by
Lemma \ref{Lemma:metric}; $W_{\alpha }=W\cap G_{\alpha }$ is an open
symmetric identity neighborhood of $G_{\alpha }$, and $U\subseteq W_{\infty
}\subseteq N\subseteq G_{\alpha }$ is an open identity neighborhood of $N$; $%
\varphi |_{W_{\alpha }}$ takes values in $H_{\alpha }$ and is $\left(
G_{\alpha },H_{\alpha }\right) $-uniformly continuous by the inductive
hypothesis; conditions (1)--(3) and (5) of Definition \ref%
{Definition:local-uniformly-continuous} are inherited from $\left(
W,U,\varphi \right) $, using that $W_{\alpha }\cap N=W\cap N=W_{\infty }$;
and the restriction of $\sigma $ to%
\begin{equation*}
\left\{ \left( z,x\right) \in N\times W_{\alpha }:zx\in W_{\alpha }\right\}
\end{equation*}%
witnesses Definition \ref{Definition:pseudomorphism}, since $\varphi \left(
zx\right) =\sigma \left( z,x\right) \varphi \left( x\right) $ and the
modulus of uniform continuity of $\sigma $ in the first variable is
independent of the second variable. Since%
\begin{equation*}
\mathrm{Ph}^{1}\left( \mathrm{Ph}^{\alpha }\left( G/N\right) \right) =%
\mathrm{Ph}^{\alpha +1}\left( G/N\right)
\end{equation*}%
and likewise for $H/M$, and since $d_{\alpha +1}$ is obtained from $%
d_{\alpha }$ and $d_{\infty }$ exactly as $d_{1}$ is obtained from $d_{0}$
and $d_{\infty }$, the case $\alpha =1$ applied to this pseudomorphism
yields (1) and (2) at $\alpha +1$.

Finally, let $\lambda <\omega _{1}$ be a limit ordinal, and suppose that (1)
and (2) hold at every $\beta <\lambda $. Since%
\begin{equation*}
G_{\lambda }=\bigcap_{\beta <\lambda }G_{\beta }\text{\quad and\quad }%
H_{\lambda }=\bigcap_{\beta <\lambda }H_{\beta }\text{,}
\end{equation*}%
assertion (1) at $\lambda $ is immediate from assertion (1) below $\lambda $%
. We prove (2). Let $\left( \lambda _{n}\right) $ and $\left( \mu
_{n}\right) $ be the increasing sequences of successor ordinals cofinal in $%
\lambda $ used to define $d_{\lambda }^{G}$ and $d_{\lambda }^{H}$,
respectively. We will use the following two facts. First, by definition,%
\begin{equation}
d_{\lambda }^{G}\geq 2^{-j-1}d_{\lambda _{j}}^{G}\text{\quad for every }j\in 
\mathbb{N}\text{.}  \label{eq:limit-metric-lower-bound}
\end{equation}%
Second, for $\beta \leq \gamma <\lambda $ the inclusion $G_{\gamma
}\rightarrow G_{\beta }$ is a continuous group homomorphism between Polish
groups, hence uniformly continuous with respect to the two-sided
uniformities; since $d_{\gamma }$ and $d_{\beta }$ are compatible \emph{%
two-sided invariant }metrics on $G_{\gamma }$ and $G_{\beta }$, these
uniformities are the ones induced by $d_{\gamma }$ and $d_{\beta }$. Hence,
for every $\eta >0$ there exists $\eta ^{\prime }>0$ such that, for $%
x,x^{\prime }\in G_{\gamma }$,%
\begin{equation}
d_{\gamma }\left( x,x^{\prime }\right) \leq \eta ^{\prime }\Rightarrow
d_{\beta }\left( x,x^{\prime }\right) \leq \eta \text{.}
\label{eq:nested-uniform-continuity}
\end{equation}

Now fix $\varepsilon >0$ and choose $N_{0}\in \mathbb{N}$ with $%
\sum_{n>N_{0}}2^{-n-1}<\varepsilon /2$; this is possible since $d_{\mu
_{n}}^{H}\leq 1$ for every $n$. For each $n\leq N_{0}$, the inductive
hypothesis at $\mu _{n}$ provides $\delta _{n}>0$ such that $x,x^{\prime
}\in W_{\mu _{n}}$ and $d_{\mu _{n}}^{G}\left( x,x^{\prime }\right) \leq
\delta _{n}$ imply%
\begin{equation*}
d_{\mu _{n}}^{H}\left( \varphi \left( x\right) ,\varphi \left( x^{\prime
}\right) \right) \leq \varepsilon /2\text{.}
\end{equation*}%
Choose $j(n)\in \mathbb{N}$ with $\lambda _{j(n) }\geq \mu _{n}$, and let $%
\eta _{n}>0$ be obtained from $\delta _{n}$ by %
\eqref{nested-uniform-continuity} applied to $\beta =\mu _{n}$ and $\gamma
=\lambda _{j(n)}$. Set%
\begin{equation*}
\delta :=\min_{n\leq N_{0}}2^{-j(n)-1}\eta _{n}>0\text{.}
\end{equation*}%
Suppose that $x,x^{\prime }\in W_{\lambda }$ satisfy $d_{\lambda }^{G}\left(
x,x^{\prime }\right) \leq \delta $. By \eqref{limit-metric-lower-bound}, $%
d_{\lambda _{j(n)}}^{G}\left( x,x^{\prime }\right) \leq 2^{j(n)+1}\delta
\leq \eta _{n}$ for every $n\leq N_{0}$, whence $d_{\mu _{n}}^{G}\left(
x,x^{\prime }\right) \leq \delta _{n}$ and therefore $d_{\mu _{n}}^{H}\left(
\varphi \left( x\right) ,\varphi \left( x^{\prime }\right) \right) \leq
\varepsilon /2$. Consequently,%
\begin{equation*}
d_{\lambda }^{H}\left( \varphi \left( x\right) ,\varphi \left( x^{\prime
}\right) \right) =\sum_{n\in \mathbb{N}}2^{-n-1}d_{\mu _{n}}^{H}\left(
\varphi \left( x\right) ,\varphi \left( x^{\prime }\right) \right) \leq 
\frac{\varepsilon }{2}\sum_{n\leq N_{0}}2^{-n-1}+\sum_{n>N_{0}}2^{-n-1}\leq
\varepsilon \text{.}
\end{equation*}%
This proves (2) at $\lambda $ and completes the induction.
\end{proof}

\section{Banach--Lie groups\label{Section:Banach-Lie}}

In this section, we recall some notions pertaining to Lie groups modeled
over possibly infinite-dimensional Banach spaces, known as Banach--Lie
groups or infinite-dimensional Lie groups; see \cite%
{gloeckner_infinite-dimensional_2026,glockner_banach-lie_2003,ando_large_2022,pestov_free_1993,omori_banach-lie_1978,neeb_infinite-dimensional_2004,ando_lie_2025}%
. The theory of Lie groups and Lie algebras is also treated in \cite%
{bourbaki_lie_1989,bourbaki_lie_2002,bourbaki_lie_2005,serre_lie_2006}.

\subsection{Banach--Lie algebras}

Recall that a Lie algebra \cite[Section 2]{upmeier_symmetric_1985} is a
vector space $X$ endowed with a binary operation $[\cdot ,\cdot ]$ (Lie
bracket) that is linear in each variable and satisfies, for $x,y,z\in X$,%
\begin{equation*}
\lbrack x,x]=0
\end{equation*}%
and%
\begin{equation*}
\lbrack \lbrack x,y],z]+[[y,z],x]+[[z,x],y]=0\text{.}
\end{equation*}%
A Banach--Lie algebra is a Banach space that is also a Lie algebra, such
that the Lie bracket is a bounded bilinear map. Following \cite[Definition
5.22]{bonfiglioli_topics_2012}, we say that the norm of a Banach--Lie
algebra $\mathfrak{g}$ is \emph{compatible }with the Lie bracket if%
\begin{equation*}
\left\Vert \left[ x,y\right] \right\Vert \leq \left\Vert x\right\Vert
\left\Vert y\right\Vert \text{\quad for }x,y\in \mathfrak{g}\text{.}
\end{equation*}%
Every Banach--Lie algebra admits an equivalent compatible norm: if $%
\left\Vert \left[ x,y\right] \right\Vert \leq M\left\Vert x\right\Vert
\left\Vert y\right\Vert $ with $M>0$, then $M\left\Vert \cdot \right\Vert $
is compatible \cite[Remark 5.27]{bonfiglioli_topics_2012}. For an
associative Banach algebra $A$ one has $\left\Vert \left[ x,y\right]
\right\Vert \leq 2\left\Vert x\right\Vert \left\Vert y\right\Vert $, so that
twice the norm of $A$ is compatible with the commutator bracket of $\mathrm{g%
}(A)$ \cite[Remark 5.26]{bonfiglioli_topics_2012}.

If $A$ is an (associative) Banach algebra, then it becomes a Banach--Lie
algebra $\mathrm{g}(A)$ with respect to the \emph{commutator product}%
\begin{equation*}
\lbrack x,y]:=xy-yx
\end{equation*}%
for $x,y\in A$ \cite[Example 2.14]{upmeier_symmetric_1985}. For a Banach
space $X$ one lets $L\left( X\right) $ be the space of bounded linear
operators on $X$ and $\mathrm{gl}\left( X\right) $ be $\mathrm{g}\left(
L\left( X\right) \right) $.

\subsection{Banach--Lie groups}

Recall that a \emph{Banach manifold} is a topological space $X$ together
with a Banach space $Z$, an open cover $\mathcal{U}$ of $X$ and, for $U\in 
\mathcal{U}$, a homeomorphism $\varphi _{U}$ (\emph{chart}) between $U$ and
an open subset $\varphi _{U}[U]$ of $Z$ such that, for $U,V\in \mathcal{U}$,
the \textquotedblleft change of coordinates\textquotedblright\ function%
\begin{equation*}
\varphi _{U}\circ \varphi _{V}^{-1}:\varphi _{V}(U\cap V)\rightarrow \varphi
_{U}(U\cap V)
\end{equation*}%
is smooth. A function between Banach manifolds is smooth if it is so
\textquotedblleft in the charts\textquotedblright. A topological group is a
Banach--Lie Polish group if it is also a Banach manifold, and the group
operations are smooth \cite[Definition 2.7]{ando_large_2022}.

As in the case of Lie groups, one can assign to a Banach--Lie group a \emph{%
Lie algebra }$\mathfrak{g}$ as its tangent space at the identity \cite[%
Section 2.2]{ando_large_2022}. It also admits an exponential function $%
\mathrm{\exp }:\mathfrak{g}\rightarrow G$ that establishes a diffeomorphism
between an open zero neighborhood of $\mathfrak{g}$ and an open identity
neighborhood of $G$ \cite[Proposition 7.1.2]%
{gloeckner_infinite-dimensional_2026}. This amounts to saying that $G$ is 
\emph{locally exponential }in the sense of \cite[Definition 7.1.1]%
{gloeckner_infinite-dimensional_2026}. In particular, $\mathrm{\exp }$ is
Lipschitz on a zero neighborhood, and its inverse is Lipschitz on an
identity neighborhood. A Banach--Lie group is furthermore \emph{exponential}
if the exponential map $\exp :\mathfrak{g}\rightarrow G$ is a global
diffeomorphism.

More generally, one can define Fr\'{e}chet--Lie groups \cite[Chapter 5]%
{gloeckner_infinite-dimensional_2026}, which are Fr\'{e}chet manifolds such
that the group operations are smooth. Generally, a Fr\'{e}chet--Lie group is
not locally exponential \cite[Example 5.6.11]%
{gloeckner_infinite-dimensional_2026}.

If $A$ is an associative unital Banach algebra, then one denotes by $\mathrm{%
G}(A)$ the (open) set of invertible elements of $A$. One has that $\mathrm{G}%
(A)$ is a Banach--Lie group \cite[Example 6.9]{upmeier_symmetric_1985},
whose Lie algebra can be identified with $\mathrm{g}(A)$. The exponential
map of $\mathrm{G}(A)$ is given by holomorphic functional calculus.

If $X$ is a Banach space, we set 
\begin{equation*}
\mathrm{GL}\left( X\right) :=\mathrm{G}\left( L\left( X\right) \right) \text{%
.}
\end{equation*}%
The Lie algebra of $\mathrm{GL}\left( X\right) $ is then $\mathrm{gl}\left(
X\right) $.

\subsection{Derivations}

If $A$ is a (not necessarily associative) algebra, then a bimodule over $A$
is a vector space $X$ endowed with bilinear operations $A\times X\rightarrow
X$ and $X\times A\rightarrow X$. A \emph{derivation} of $A$ into $X$ is then
a linear map%
\begin{equation*}
\delta :A\rightarrow X
\end{equation*}%
satisfying%
\begin{equation*}
\delta \left( xy\right) =\left( \delta x\right) y+x\left( \delta y\right)
\end{equation*}%
for $x,y\in A$ \cite[Section 2]{upmeier_symmetric_1985}.

If $A$ is a (not necessarily associative) Banach algebra, then the set $%
\mathrm{aut}(A)$ of continuous derivations of $A$ is a closed subalgebra of $%
\mathrm{gl}(A)$, whence a Banach--Lie algebra \cite[Proposition 2.15]%
{upmeier_symmetric_1985}. Likewise, the set $\mathrm{Aut}(A)$ of continuous
algebra automorphisms of $A$ is a closed subgroup of $\mathrm{GL}(A)$ \cite[%
Proposition 2.15]{upmeier_symmetric_1985}, whence a Banach--Lie group. The
Lie algebra of $\mathrm{Aut}(A)$ can be identified with $\mathrm{aut}(A)$,
where the exponential map $\mathrm{aut}(A)\rightarrow \mathrm{\mathrm{Aut}}%
(A)$ is the restriction of the exponential map $\mathrm{gl}(A)\rightarrow 
\mathrm{GL}(A)$ \cite[Proposition 2.15]{upmeier_symmetric_1985}.

If $A$ is an associative Banach algebra, then%
\begin{equation*}
\mathrm{ad}:\mathrm{g}(A)\rightarrow \mathrm{aut}(A)\text{, }x\mapsto 
\mathrm{ad}\left( x\right) =[x,-]
\end{equation*}%
is a continuous Lie algebra homomorphism \cite[Proposition 2.16]%
{upmeier_symmetric_1985}. More generally, if $\mathfrak{g}$ is a Banach--Lie
algebra, then%
\begin{equation*}
\mathrm{ad}:\mathfrak{g}\rightarrow \mathrm{aut}\left( \mathfrak{g}\right) 
\text{, }x\mapsto \mathrm{ad}\left( x\right) :=[x,-]
\end{equation*}%
is a continuous Lie algebra homomorphism. Its kernel is the \emph{center }of 
$\mathfrak{g}$. The image%
\begin{equation*}
\mathrm{int}\left( \mathfrak{g}\right) =\left\{ \mathrm{ad}\left( x\right)
:x\in \mathfrak{g}\right\}
\end{equation*}%
is an ideal of $\mathrm{aut}\left( \mathfrak{g}\right) $, whose elements are
called \emph{inner }derivations \cite[Proposition 2.17]%
{upmeier_symmetric_1985}.

If $A$ is an associative unital Banach algebra, one can also define a
continuous group homomorphism%
\begin{equation*}
\mathrm{Ad}:\mathrm{G}(A)\rightarrow \mathrm{\mathrm{Aut}}(A)\text{, }%
u\mapsto \mathrm{Ad}\left( u\right)
\end{equation*}%
where%
\begin{equation*}
\mathrm{Ad}\left( u\right) :x\mapsto uxu^{-1}\text{.}
\end{equation*}%
The diagram%
\begin{equation*}
\begin{array}{ccc}
\mathrm{G}(A) & \overset{\mathrm{Ad}}{\longrightarrow } & \mathrm{\mathrm{Aut%
}}(A) \\ 
\mathrm{\exp }\uparrow &  & \uparrow \mathrm{\exp } \\ 
\mathrm{g}(A) & \underset{\mathrm{ad}}{\longrightarrow } & \mathrm{aut}(A)%
\end{array}%
\end{equation*}%
commutes \cite[Proposition 2.18]{upmeier_symmetric_1985}.

\subsection{The Baker--Campbell--Hausdorff Formula\label{Subsection:BCH}}

Let $\mathfrak{g}$ be a Banach--Lie algebra whose norm is \emph{compatible }%
with the Lie bracket. We adopt the notation from \cite[Section 5.2.2]%
{bonfiglioli_topics_2012}. Define 
\begin{equation*}
\Delta :=\left\{ \left( x,y\right) \in \mathfrak{g}\times \mathfrak{g}:\max
\left\{ \left\Vert x\right\Vert ,\left\Vert y\right\Vert \right\}
<2^{-1}\log 2\right\} \text{.}
\end{equation*}%
The Hausdorff series converges absolutely on $\Delta $ and defines a
continuous function%
\begin{equation*}
\Delta \rightarrow \mathfrak{g}\text{, }\left( a,b\right) \mapsto a\ast b%
\text{;}
\end{equation*}%
see \cite[Theorem 5.30]{bonfiglioli_topics_2012}, and \cite%
{wojtynski_quasinilpotent_1998} for the Banach--Lie algebras on which the
Hausdorff series converges \emph{globally}. On every set 
\begin{equation*}
\left\{ \left( a,b\right) :\max \left\{ \left\Vert a\right\Vert , \left\Vert
b\right\Vert \right\} \leq r\right\}
\end{equation*}%
with $0<r<2^{-1}\log 2$, this function is uniformly continuous. Indeed, the
Fundamental Estimate \cite[Theorem 5.29]{bonfiglioli_topics_2012} and \cite[%
Theorem 5.30]{bonfiglioli_topics_2012} give normal, hence uniform,
convergence there. Each homogeneous term of the Hausdorff series is the
diagonal restriction of a continuous multilinear map, and is therefore
uniformly continuous on bounded sets. The partial sums are consequently
uniformly continuous, and so is their uniform limit.

If the norm of $\mathfrak{g}$ satisfies $\left\Vert \left[ x,y\right]
\right\Vert \leq M\left\Vert x\right\Vert \left\Vert y\right\Vert $ rather
than being compatible, the same conclusions hold with $2^{-1}\log 2$
replaced by $(2M)^{-1}\log 2$, as one sees by passing to the compatible norm 
$M\left\Vert \cdot \right\Vert $.

We record once and for all the normalization that will be used below. The
Baker--Campbell--Hausdorff formula will only be applied to the Banach--Lie
algebra $\mathrm{aut}(A)$ of *-derivations of a C*-algebra $A$, regarded as
a Banach--Lie subalgebra of $\mathrm{gl}(A)$. There the operator norm
satisfies $\left\Vert \left[ \delta ,\delta ^{\prime }\right] \right\Vert
\leq 2\left\Vert \delta \right\Vert \left\Vert \delta ^{\prime }\right\Vert $%
, so that \emph{twice }the operator norm is compatible with the Lie bracket.
Accordingly we set%
\begin{equation*}
\Delta _{A}:=\left\{ \left( x,y\right) \in \mathrm{aut}(A)\times \mathrm{aut}%
(A):\left\Vert x\right\Vert +\left\Vert y\right\Vert <2^{-1}\log 2\right\} 
\text{,}
\end{equation*}%
where $\left\Vert \cdot \right\Vert $ denotes the operator norm. By the
above, the Hausdorff series converges absolutely on $\Delta _{A}$, and
normally---hence uniformly---on%
\begin{equation*}
\left\{ \left( x,y\right) :\left\Vert x\right\Vert +\left\Vert y\right\Vert
\leq \eta \right\}
\end{equation*}%
for every $\eta <2^{-1}\log 2$. \emph{All the Hausdorff products occurring
below are formed within }$\Delta _{A}$. Only the existence of \emph{some }%
positive radius with these properties is used, so the precise value of the
constant is immaterial.

Suppose that $G$ is a Banach--Lie group with Banach--Lie algebra $\mathfrak{g%
}$. Then $G$ is locally exponential \cite[Proposition 7.1.2]%
{gloeckner_infinite-dimensional_2026}, and the local multiplication that $%
\exp $ transports from $G$ to a sufficiently small ball around $0$ in $%
\mathfrak{g}$ is given by the Hausdorff series \cite[Lemma 8.2.2 and Theorem
8.4.1]{gloeckner_infinite-dimensional_2026}. Thus, there exists $\rho >0$
such that%
\begin{equation*}
\exp \left( a\ast b\right) =\exp \left( a\right) \exp \left( b\right)
\end{equation*}%
whenever $\left\Vert a\right\Vert ,\left\Vert b\right\Vert \leq \rho $. This
is the only form of the Banach--Lie Baker--Campbell--Hausdorff theorem that
will be used below.

\subsection{Ideals\label{Subsection:ideals}}

Let $\mathfrak{g}$ be a Banach--Lie algebra whose norm is compatible with
its Lie bracket. An \emph{ideal} of $\mathfrak{g}$ is a subspace $\mathfrak{h%
}$ such that for $x\in \mathfrak{h}$ and $a\in \mathfrak{g}$, $\left[ a,x%
\right] \in \mathfrak{h}$. Thus, $\mathfrak{h}$ is invariant with respect to
all the inner derivations of $\mathfrak{g}$ \cite[Section 1.4]%
{bourbaki_lie_1989}. A \emph{Banach--Lie ideal }of $\mathfrak{g}$ is a (not
necessarily closed) ideal that is also a Banach--Lie algebra, with respect
to a norm $\left\Vert \cdot \right\Vert _{\mathfrak{h}}$ for which the
inclusion $\mathfrak{h}\rightarrow \mathfrak{g}$ is continuous and which
satisfies%
\begin{equation*}
\left\Vert \left[ a,x\right] \right\Vert _{\mathfrak{h}}\leq \left\Vert
a\right\Vert _{\mathfrak{g}}\left\Vert x\right\Vert _{\mathfrak{h}}
\end{equation*}%
for $x\in \mathfrak{h}$ and $a\in \mathfrak{g}$. After multiplying $%
\left\Vert \cdot \right\Vert _{\mathfrak{h}}$ by a constant, we also assume
that $\left\Vert x\right\Vert _{\mathfrak{g}}\leq \left\Vert x\right\Vert _{%
\mathfrak{h}}$ for $x\in \mathfrak{h}$. Thus, the restriction of inner
derivations from $\mathfrak{g}$ to $\mathfrak{h}$ defines a linear map%
\begin{equation*}
\mathfrak{g}\rightarrow \mathrm{aut}\left( \mathfrak{h}\right) \text{, }%
a\mapsto \mathrm{ad}\left( a\right) |_{\mathfrak{h}}
\end{equation*}%
of norm at most $1$, which factors through $\mathrm{int}\left( \mathfrak{g}%
\right) $.

The Fundamental Estimate and the convergence estimates for the Hausdorff
series \cite[Theorems 5.29 and 5.30]{bonfiglioli_topics_2012}, applied both
in $\mathfrak{g}$ and with the mixed estimate above, show that, if one sets%
\begin{equation*}
\Delta _{\mathfrak{h}}:=\left\{ \left( x,a\right) \in \mathfrak{h}\times 
\mathfrak{g}:\max \left\{ \left\Vert a\right\Vert _{\mathfrak{g}},\left\Vert
x\right\Vert _{\mathfrak{g}},\left\Vert x\right\Vert _{\mathfrak{h}}\right\}
<2^{-1}\log 2\right\} \text{,}
\end{equation*}%
then for $\left( x,a\right) \in \Delta _{\mathfrak{h}}$, the terms of $x\ast
a-a$ converge in $\mathfrak{h}$. Since the inclusion $\mathfrak{h}%
\rightarrow \mathfrak{g}$ is continuous, this sum agrees with the Hausdorff
sum computed in $\mathfrak{g}$. In particular, 
\begin{equation*}
x\ast a-a\in \mathfrak{h}\text{,}
\end{equation*}%
and the function%
\begin{eqnarray*}
\Delta _{\mathfrak{h}} &\rightarrow &\mathfrak{h} \\
\left( x,a\right) &\mapsto &x\ast a-a
\end{eqnarray*}%
is continuous. The same multilinear estimate used above shows that its
restriction to every subdomain obtained by replacing the strict bound in the
definition of $\Delta _{\mathfrak{h}}$ by $\leq r$, where $0<r<2^{-1}\log 2$%
, is uniformly continuous with respect to the product uniform structure.

\subsection{Unitary groups of C*-algebras}

Let $A$ be a separable \emph{unital} C*-algebra \cite%
{blackadar_operator_2006,pedersen_algebras_1979}. Define $U(A)$ to be its 
\emph{unitary group}, which is a Banach--Lie Polish group \cite[Section 3.2]%
{ando_large_2022}. Likewise, we can regard $A$ as a Banach--Lie algebra with
respect to the Lie bracket%
\begin{equation*}
\lbrack x,y]:=xy-yx\text{.}
\end{equation*}%
The Lie algebra of $U(A)$ can be identified with $iA_{\mathrm{sa}}$ with
exponential map given by holomorphic functional calculus; see \cite[Section
3.2]{ando_large_2022}.

The group of Banach algebra automorphisms of $A$ is a closed subgroup of $%
\mathrm{GL}(A)$ \cite[Proposition 2.15]{upmeier_symmetric_1985}, and the
group $\mathrm{\mathrm{Aut}}(A)$ of C*-algebra automorphisms of $A$ is in
turn closed in it, being the set of those Banach algebra automorphisms that
commute with the involution. Equipped with the operator norm topology, it is
a Banach--Lie group modeled on the Banach--Lie algebra $\mathrm{aut}(A)$ of
*-derivations of $A$. Indeed, if $\left\Vert \alpha-\mathrm{id}%
_{A}\right\Vert <2$, then the principal logarithm $\mathrm{Log}(\alpha)$ is
a *-derivation and $\exp(\mathrm{Log}(\alpha))=\alpha$; thus $\mathrm{Log}$
is a local chart at the identity inverse to the exponential map; see \cite[%
Theorem 8.7.7]{pedersen_algebras_1979}. By definition, the *-derivations are
the derivations mapping self-adjoint elements to self-adjoint elements \cite[%
Section 8.6]{pedersen_algebras_1979}.

\section{Derivation length\label{Section:derivation}}

In this section, we define the notion of \emph{derivation length }of a
unital C*-algebra, which measures the complexity of its space of
*-derivations modulo inner ones.

\subsection{Derivations}

Let $A$ be a separable unital C*-algebra. If $M$ is the double dual of $A$,
or any von Neumann algebra that contains $A$ as a $\sigma $-weakly dense
subalgebra, every *-derivation of $A$ extends to a $\sigma $-weakly
continuous derivation of $M$. Furthermore, such an extension is inner (as a
derivation of $M$). If $\delta $ is a *-derivation, then by \cite[Theorem 3.1%
]{kadison_derivations_1967} there exists a \emph{unique} $a\in M_{\mathrm{sa}%
}$ with $\delta =\mathrm{ad}\left( ia\right) $ and such that for every
central projection $p$ of $M$, 
\begin{equation*}
\left\Vert pa\right\Vert =\frac{1}{2}\left\Vert \delta |_{pM}\right\Vert 
\text{.}
\end{equation*}

We let $\mathrm{aut}(A)$ be the space of *-derivations of $A$. This is a
Banach space with respect to the operator norm, and a locally convex
topological vector space with respect to the strong operator topology. We
regard the subspace $\mathrm{inn}(A)$ of inner derivations as a separable
Banach space with respect to the norm induced by the surjective linear map $%
iA_{\mathrm{sa}}\rightarrow \mathrm{inn}\left( A\right) $, $a\mapsto \mathrm{%
ad}\left( a\right) $ with kernel $Z\left( iA_{\mathrm{sa}}\right) $. Thus,%
\begin{equation*}
\mathrm{out}(A):=\frac{\mathrm{aut}(A)}{\mathrm{inn}(A)}
\end{equation*}%
is a \emph{phantom topological vector space}. Furthermore, $\mathrm{\mathrm{%
Ball}}(\mathrm{inn}\left( A\right) )$ is dense in $\mathrm{\mathrm{Ball}}(%
\mathrm{aut}\left( A\right) )$ with respect to the strong operator topology;
see \cite[Theorem 8.6.5]{pedersen_algebras_1979} and \cite[Lemma 8.6.12]%
{pedersen_algebras_1979}.

\begin{remark}
\label{Remark:aut-not-Polish}The pair $\left( \mathrm{aut}(A),\mathrm{inn}%
(A)\right) $ is \emph{not }in general a group with a Polish cover in the
sense of Section \ref{Section:phantom}: neither the operator norm topology
nor the strong operator topology on $\mathrm{aut}(A)$ need be Polish. For
instance, let%
\begin{equation*}
A:=\left\{ x\in \prod_{n\in \mathbb{N}}M_{2}\left( \mathbb{C}\right) :\left(
x_{n}\right) \text{ converges to a scalar}\right\} \text{,}
\end{equation*}%
a separable unital C*-algebra. Every *-derivation of $A$ preserves each
summand $M_{2}\left( \mathbb{C}\right) $, whence $\mathrm{aut}(A)\cong \ell
_{\infty }(\mathbb{N},\mathfrak{su}(2))$ and $\mathrm{inn}(A)\cong c_{0}(%
\mathbb{N},\mathfrak{su}(2))$. In particular $\mathrm{aut}(A)$ is not
separable in the operator norm, and it is not Polish in the strong operator
topology either, being the increasing union of the countably many closed
balls $n\mathrm{\mathrm{Ball}}\left( \mathrm{aut}(A)\right) $, none of which
has nonempty interior.

Accordingly, we do \emph{not} apply the theory of Section \ref%
{Section:phantom} to $\mathrm{out}(A)$ itself, but only to its first phantom
subspace, where it does apply; see Lemma \ref%
{Lemma:phantom-derivation-order-1} and the convention introduced immediately
below. This is no loss, since all the invariants considered in this paper
are computed inside $\mathrm{Ph}^{1}\mathrm{out}\left( A\right) $.
\end{remark}

We set $\mathrm{Ph}^{0}\mathrm{aut}(A):=\mathrm{aut}(A)$, and we let $%
\mathrm{Ph}^{1}\mathrm{aut}(A)$ be the closure of $\mathrm{inn}(A)$ in $%
\mathrm{aut}(A)$ with respect to the operator norm, so that%
\begin{equation*}
\mathrm{Ph}^{1}\mathrm{out}(A)=\frac{\mathrm{Ph}^{1}\mathrm{aut}(A)}{\mathrm{%
inn}(A)}\text{.}
\end{equation*}%
By Lemma \ref{Lemma:phantom-derivation-order-1} below, $\mathrm{Ph}^{1}%
\mathrm{aut}(A)$ is a \emph{separable }Banach space in which $\mathrm{inn}%
(A) $ is a Polishable subspace; thus $\mathrm{Ph}^{1}\mathrm{out}(A)$ \emph{%
is }a phantom Banach space, and the theory of Section \ref{Section:phantom}
applies to it. For $1\leq \alpha <\omega _{1}$ we write $\alpha =1+\xi $ and
define $\mathrm{Ph}^{\alpha }\mathrm{aut}(A)$ to be the subspace of $\mathrm{%
Ph}^{1}\mathrm{aut}(A)$ determined by%
\begin{equation*}
\mathrm{Ph}^{\xi }\left( \mathrm{Ph}^{1}\mathrm{out}(A)\right) =\frac{%
\mathrm{Ph}^{1+\xi }\mathrm{aut}(A)}{\mathrm{inn}(A)}\text{,}
\end{equation*}%
which is consistent at $\alpha =1$. A *-derivation of $A$ is phantom of
order at least $\alpha $ if and only if it belongs to $\mathrm{Ph}^{\alpha }%
\mathrm{aut}(A)$. It follows from \cite[Theorem 3]{elliott_some_1977} that
whenever $\mathrm{out}\left( A\right) $ is nonzero, $\mathrm{Ph}^{1}\mathrm{%
out}\left( A\right) $ is a proper subspace of $\mathrm{out}\left( A\right) $.

\begin{lemma}
\label{Lemma:phantom-derivation-order-1}Let $A$ be a separable \emph{unital}%
\ C*-algebra. Then:

\begin{enumerate}
\item the norm-closure $\mathrm{Ph}^{1}\mathrm{aut}(A)$ of $\mathrm{inn}(A)$
in $\mathrm{aut}(A)$ coincides with the set of *-derivations $\delta $ of $A$
such that, for every $\varepsilon >0$, there exist $b\in A_{\mathrm{sa}}$
and a sequence $\left( z_{n}\right) $ in $\mathrm{\mathrm{Ball}}%
_{\varepsilon }\left( A_{\mathrm{sa}}\right) $ with $\mathrm{ad}\left(
iz_{n}\right) \rightarrow \delta +\mathrm{ad}\left( ib\right) $ in the
strong operator topology; in particular, $\mathrm{Ph}^{1}\mathrm{aut}(A)$ is
a separable Banach space in which $\mathrm{inn}(A)$ is a Polishable subspace;

\item the topology on $\mathrm{Ph}^{1}\mathrm{\mathrm{aut}}(A)$ induced by
the operator norm is the unique Polish group topology on $\mathrm{Ph}^{1}%
\mathrm{\mathrm{aut}}\left( A\right) $ making the inclusion $\mathrm{Ph}^{1}%
\mathrm{aut}\left( A\right) \rightarrow \mathrm{aut}\left( A\right) $
continuous.
\end{enumerate}
\end{lemma}

\begin{proof}
(1) For $k\in \omega $ let $V_{k}$ be the set of inner derivations of $A$ of
the form $\mathrm{ad}\left( ix\right) $ for some $x\in A_{\mathrm{sa}}$ with 
$\left\Vert x\right\Vert \leq 2^{-(k+1)}$. Then $\left( V_{k}\right) $ is a
basis of zero neighborhoods of $\mathrm{inn}(A)$. Furthermore, the closure $%
\overline{V}_{k}^{\mathrm{aut}(A)}$ of $V_{k}$ in $\mathrm{aut}(A)$ with
respect to the strong operator topology is contained in $\mathrm{\mathrm{Ball%
}}_{2^{-k}}\left( \mathrm{aut}\left( A\right) \right) $. Denote by $X$ the
closure of $\mathrm{inn}(A)$ with respect to the topology induced by the
operator norm, and by $Y$ the set of *-derivations described in the
statement. Suppose that $\delta \in Y$. Then for every $\varepsilon >0$
there exist a sequence $\left( z_{n}\right) $ in $\mathrm{\mathrm{Ball}}%
_{\varepsilon }\left( A_{\mathrm{sa}}\right) $ and $b\in A_{\mathrm{sa}}$
such that $\delta +\mathrm{ad}\left( ib\right) $ is the limit of the
sequence $\left( \mathrm{ad}\left( iz_{n}\right) \right) $ in the strong
operator topology. In particular, this implies that $\left\Vert \delta +%
\mathrm{ad}\left( ib\right) \right\Vert \leq 2\varepsilon $ and hence $%
\delta \in X$.

Conversely, suppose that $\delta \in X$. Fix $\varepsilon >0$. Then there
exists $b\in A_{\mathrm{sa}}$ such that $\left\Vert \delta +\mathrm{ad}%
\left( ib\right) \right\Vert <\varepsilon $. Thus, by \cite[Theorem 8.6.5]%
{pedersen_algebras_1979} and \cite[Lemma 8.6.12]{pedersen_algebras_1979},
there exists a sequence $\left( z_{n}\right) $ in $\mathrm{Ball}%
_{\varepsilon }\left( A_{\mathrm{sa}}\right) $ such that $\delta +\mathrm{ad}%
\left( ib\right) $ is the limit of $\left( \mathrm{ad}\left( iz_{n}\right)
\right) $ in the strong operator topology. Since $\varepsilon $ was
arbitrary, $\delta \in Y$. Thus $X=Y$. Finally, $X$ is separable, being the
norm-closure of the separable space $\mathrm{inn}(A)$, and $\mathrm{inn}(A)$
is a Polishable subspace of $X$, since it is a separable Banach space whose
inclusion into $X$ is continuous.

(2) Suppose that $\tau $ is a Polish group topology on $\mathrm{Ph}^{1}%
\mathrm{aut}\left( A\right) $ making the inclusion into $\mathrm{aut}\left(
A\right) $, endowed with the strong operator topology, continuous, and let $%
\tau _{0}$ be the operator norm topology. Both $\tau $ and $\tau _{0}$ are
Polish group topologies whose inclusions into the Hausdorff space $\mathrm{%
aut}\left( A\right) $ with the strong operator topology are continuous.
Hence the diagonal 
\begin{equation*}
\Gamma :=\left\{ \left( \delta ,\delta \right) :\delta \in \mathrm{Ph}^{1}%
\mathrm{aut}\left( A\right) \right\}
\end{equation*}%
is closed in the product of $\mathrm{Ph}^{1}\mathrm{aut}\left( A\right) $
endowed with $\tau $ and the same group endowed with $\tau _{0}$. Thus $%
\Gamma $ is a Polish group, and its two coordinate projections are
continuous bijective homomorphisms onto the two Polish groups. By the
Lusin--Souslin theorem \cite[Theorem 15.1]{kechris_classical_1995}, the
inverses of these projections are Borel. Therefore the identity map between
the two topologies is Borel in both directions. By \cite[Theorem 9.10]%
{kechris_classical_1995}, a Borel---hence Baire measurable---group
homomorphism between Polish groups is continuous, so $\tau =\tau _{0}$.
\end{proof}

It follows from Lemma \ref{Lemma:phantom-derivation-order-1} that $\mathrm{Ph%
}^{1}\mathrm{aut}\left( A\right) $ is a phantom Banach space. By this and 
\cite[Proposition 10.1]{lupini_complexity_2025}, $\mathrm{Ph}^{\alpha }%
\mathrm{aut}\left( A\right) $ is a phantom Banach space for every countable
successor ordinal $\alpha $.

\begin{definition}
\label{Definition:derivation-length} Let $A$ be a separable unital
C*-algebra. The \emph{derivation length }$\ell _{\mathrm{aut}}(A)\in \omega
_{1}^{\mathrm{Ph}}$ of $A$ is:

\begin{itemize}
\item $0$ if every *-derivation of $A$ is inner;

\item $1+\alpha $ where $\alpha $ is the phantom length of $\mathrm{Ph}^{1}%
\mathrm{out}\left( A\right) $, if $A$ has outer derivations.
\end{itemize}
\end{definition}

One can recursively characterize the phantom derivations of order $\alpha $
as follows:

\begin{definition}
Suppose that $A$ is a separable unital C*-algebra. We denote by $r$ an
element of $\left[ 0,\infty \right] $. A *-derivation $\delta $ of $A$ is:

\begin{enumerate}
\item $r$-phantom of order at least $1$ if, for every $\eta >0$, there
exists $a\in A_{\mathrm{sa}}$ such that $\left\Vert a\right\Vert \leq r$ and 
$\left\Vert \delta +\mathrm{ad}\left( ia\right) \right\Vert \leq \eta $;

\item $r$-phantom of order at least $\beta +1$ for some $\beta \geq 1$ if,
for every $\eta >0$, there exists $a\in A_{\mathrm{sa}}$ such that $%
\left\Vert a\right\Vert \leq r$ and $\delta +\mathrm{ad}\left( ia\right) $
is $\eta $-phantom of order at least $\beta $;

\item $r$-phantom of order at least $\lambda $ for some limit ordinal $%
\lambda $ if it is $r$-phantom of order at least $\beta $ for every $1\leq
\beta <\lambda $;

\item phantom of order at least $\alpha $ if and only if it is $\infty $%
-phantom of order at least $\alpha $.
\end{enumerate}
\end{definition}

\subsection{The von Neumann corona}

Suppose now that $A\subseteq B\left( H\right) $ is a concrete separable
unital C*-algebra with center $ZA$. Let $M$ be the $\sigma $-weak closure of 
$A$ in $B\left( H\right) $, and $ZM$ be its center. We define the \emph{%
reduced} versions%
\begin{equation*}
\rho A:=A_{\mathrm{sa}}/ZA_{\mathrm{sa}}
\end{equation*}%
and%
\begin{equation*}
\rho M:=M_{\mathrm{sa}}/ZM_{\mathrm{sa}}
\end{equation*}%
of $A$ and $M$, respectively. We can identify $\rho A$ with a subspace of $%
\rho M$, and let $\delta A$ be the norm-closure of $\rho A$ in $\rho M$. We
can thus consider the phantom Banach space 
\begin{equation*}
\partial A:=\delta A/\rho A\text{.}
\end{equation*}

\begin{lemma}
\label{Lemma:vN-corona}Let $A\subseteq B\left( H\right) $ be a concrete
separable unital C*-algebra. The continuous linear map%
\begin{equation*}
L:\delta A\rightarrow \mathrm{Ph}^{1}\mathrm{aut}(A)
\end{equation*}%
\begin{equation*}
x\mapsto \mathrm{ad}\left( ix\right)
\end{equation*}%
induces an iso-airomorphism%
\begin{equation*}
\partial A\rightarrow \mathrm{Ph}^{1}\mathrm{out}(A)
\end{equation*}%
and, in particular,%
\begin{equation*}
\ell _{\mathrm{Ph}}\left( \partial A\right) =\ell _{\mathrm{Ph}}\left( 
\mathrm{Ph}^{1}\mathrm{out}(A)\right) \text{.}
\end{equation*}
\end{lemma}

\begin{proof}
By \cite[Theorem 3.1]{kadison_derivations_1967}, every *-derivation of $A$
is implemented by an element of $iM_{\mathrm{sa}}$, uniquely modulo $iZM_{%
\mathrm{sa}}$. The norm formula for derivations on a von Neumann algebra 
\cite{zsido_norm_1973}, together with the Kaplansky Density Theorem \cite[%
I.9.1.3]{blackadar_operator_2006}, shows that%
\begin{equation*}
\left\Vert \mathrm{ad}\left( ix\right) \right\Vert =2\left\Vert x+ZM_{%
\mathrm{sa}}\right\Vert
\end{equation*}%
for $x\in M_{\mathrm{sa}}$.

Hence $x\mapsto \mathrm{ad}(ix)$ is a topological linear isomorphism from $%
\rho M$ onto $\mathrm{aut}(A)$ and maps $\rho A$ onto $\mathrm{inn}(A)$.

Taking closures gives the asserted topological linear isomorphism $L:\delta
A\rightarrow \mathrm{Ph}^{1}\mathrm{aut}(A)$. Thus $L$ and $L^{-1}$ are
continuous homomorphism lifts of the induced bijection and its inverse,
respectively. The induced map is therefore an iso-airomorphism, and the
equality of phantom lengths follows from Subsection \ref%
{Subsection:homomorphisms}.
\end{proof}

\subsection{Multipliers}

Suppose now that $A$ is a (not necessarily unital) separable C*-algebra. Let 
$M(A)$ be the \emph{multiplier algebra} of $A$ \cite[Section\ II.7.3]%
{blackadar_operator_2006}. Then $\mathrm{\mathrm{Ball}}\left( M(A)\right) $
is a Polish space with respect to the \emph{strict topology}, generated by
the seminorms%
\begin{equation*}
\nu _{a}\left( x\right) :=\left\Vert xa\right\Vert +\left\Vert ax\right\Vert
\end{equation*}%
for $a\in A$.

Derivations and *-derivations on $A$ are defined as in the unital case \cite[%
Section 8.6]{pedersen_algebras_1979}. In this case, one defines the space $%
\mathrm{inn}(M(A))$ of *-derivations of $A$ of the form $\mathrm{ad}\left(
ix\right) $ for $x\in M(A)_{\mathrm{sa}}$. Two *-derivations $\delta ,\delta
^{\prime }$ of $A$ are \emph{multiplier equivalent }if $\delta ^{\prime
}-\delta \in \mathrm{inn}(M(A))$.

The operator norm $\left\Vert \delta \right\Vert $ of a *-derivation of $A$
is equal to its spectral radius \cite[Theorem 8.6.5]{pedersen_algebras_1979}%
. The corresponding implementation statement also holds when $A$ is not
unital. Indeed, let $\widetilde{A}$ be the unitization of $A$ (and take $%
\widetilde{A}=A$ when $A$ is unital), and extend $\delta $ to $\widetilde{A}$
by setting $\delta(1)=0$. If $M$ contains $A$ as a $w^{\ast }$-dense
subalgebra, the inclusion of $A$ in $M$ extends to a unital representation
of $\widetilde{A}$ whose weak closure is still $M$. Applying \cite[Theorem
3.1]{kadison_derivations_1967} to this unitization shows that $\delta $ is $%
\sigma $-weakly continuous and extends to a *-derivation of $M$.
Furthermore, there exists a \emph{unique} $a\in M_{\mathrm{sa}}$ with $%
\delta =\mathrm{ad}\left( ia\right) $ and such that for every central
projection $p$ of $M$, 
\begin{equation*}
\left\Vert pa\right\Vert =\frac{1}{2}\left\Vert \delta |_{pM}\right\Vert 
\text{.}
\end{equation*}

\begin{proposition}
\label{Proposition:multiplier-equivalence-derivations-analytic}Let $A$ be a
separable C*-algebra. For every $n\in \mathbb{N}$, the relation of
multiplier equivalence on the closed ball $n\mathrm{\mathrm{Ball}}(\mathrm{%
aut}(A))$ is analytic.
\end{proposition}

\begin{proof}
Fix $n,m\in \mathbb{N}$. The space $n\mathrm{\mathrm{Ball}}(\mathrm{aut}(A))$%
, with the strong operator topology, is Polish: it is a closed subspace of
the norm-$n$ ball of bounded linear maps from $A$ to $A$, endowed with the
topology of pointwise convergence. The space $m\mathrm{\mathrm{Ball}}(M(A)_{%
\mathrm{sa}})$ is Polish in the strict topology. Define%
\begin{align*}
R_{n,m}:=\{&(\delta ,\delta ^{\prime },x)\in n\mathrm{\mathrm{Ball}}(\mathrm{%
aut}(A))^{2}\times m\mathrm{\mathrm{Ball}}(M(A)_{\mathrm{sa}}): \\
&\delta ^{\prime }-\delta =\mathrm{ad}(ix)\}.
\end{align*}%
Fix a dense sequence $(a_{k})$ in $A$. Since all the maps under
consideration are bounded, membership in $R_{n,m}$ is equivalent to%
\begin{equation*}
\delta ^{\prime }(a_{k})-\delta (a_{k})=i[x,a_{k}] \quad\text{for every }%
k\in \mathbb{N}.
\end{equation*}%
The evaluation maps on $n\mathrm{\mathrm{Ball}}(\mathrm{aut}(A))$ are
continuous, while on the bounded multiplier ball the maps $x\mapsto xa_{k}$
and $x\mapsto a_{k}x$ are strict-to-norm continuous. Hence $R_{n,m}$ is
closed. Finally, multiplier equivalence on $n\mathrm{\mathrm{Ball}}(\mathrm{%
aut}(A))$ is the analytic set%
\begin{equation*}
\bigcup_{m\in \mathbb{N}}\mathrm{proj}_{1,2}(R_{n,m}),
\end{equation*}%
which concludes the proof.
\end{proof}

\section{Automorphism length\label{Section:automorphism}}

We now recall the definition of the \emph{automorphism length }of a
separable unital C*-algebra $A$ used in work in preparation by Bergfalk,
Lupini, and Panagiotopoulos \cite{bergfalk_definable_2026}.\ This measures
the complexity of describing the group of automorphisms of $A$ modulo inner
ones.

\subsection{Automorphisms}

Let $A$ be a separable C*-algebra, and let $M(A)$ be its multiplier algebra.
We write $U(A):=U(M(A))$ for the unitary group of $M(A)$, consistently with
the notation of Section \ref{Section:Banach-Lie} when $A$ is unital. On the
bounded set $\mathrm{Ball}(M(A))$, multiplication is jointly continuous in
the strict topology, and the involution is strictly continuous. Consequently,%
\begin{equation*}
U(A)=\left\{u\in \mathrm{Ball}(M(A)):u^{\ast }u=uu^{\ast }=1\right\}
\end{equation*}%
is strictly closed in $\mathrm{Ball}(M(A))$. It is therefore a Polish group
in the strict topology.

$\mathrm{Aut}(A)$ is a Polish group endowed with the topology of pointwise
convergence. The normal subgroup 
\begin{equation*}
\mathrm{Inn}(A)=\left\{ \mathrm{Ad}\left( u\right) :u\in U(A)\right\}
\end{equation*}%
of \emph{inner }automorphisms is a Polish subgroup in the sense of Section %
\ref{Section:phantom}. Indeed, the kernel of the continuous homomorphism%
\begin{equation*}
u\mapsto \mathrm{Ad}\left( u\right) \text{.}
\end{equation*}%
is closed in $U(A)$, so $U(A)/\ker(\mathrm{Ad})$ is a Polish group. The
induced homomorphism from this quotient into $\mathrm{Aut}(A)$ is continuous
and injective, with image $\mathrm{Inn}(A)$. In particular, $\mathrm{Inn}(A)$
is Borel in $\mathrm{Aut}(A)$ by the Lusin--Souslin theorem \cite[Theorem
15.1]{kechris_classical_1995}. We define the group with a Polish cover 
\begin{equation*}
\mathrm{Out}(A):=\frac{\mathrm{\mathrm{Aut}}(A)}{\mathrm{Inn}(A)}\text{.}
\end{equation*}%
We have the corresponding phantom subgroups%
\begin{equation*}
\mathrm{Ph}^{\alpha }\mathrm{Out}(A)=\frac{\mathrm{Ph}^{\alpha }\mathrm{%
\mathrm{Aut}}(A)}{\mathrm{Inn}(A)}
\end{equation*}%
for $\alpha <\omega _{1}$.\ The elements of $\mathrm{Ph}^{\alpha }\mathrm{Aut%
}(A)$ are called phantom automorphisms of order at least $\alpha $. In
particular, $\mathrm{Ph}^{0}\mathrm{\mathrm{Aut}}(A)$ is the closure of $%
\mathrm{Inn}(A)$ within $\mathrm{\mathrm{Aut}}(A)$ with respect to the
topology of pointwise convergence. The elements of $\mathrm{Ph}^{0}\mathrm{%
\mathrm{Aut}}(A)$ are also called \emph{approximately inner }automorphisms.

\begin{definition}
Let $A$ be a separable C*-algebra. The \emph{automorphism length }$\ell _{%
\mathrm{Aut}}(A)\in \omega _{1}^{\mathrm{Ph}}$ of $A$ is the phantom length
of $\mathrm{Out}(A)$.
\end{definition}

A \emph{unital }separable C*-algebra has automorphism length zero if and
only if every central sequence is trivial or, equivalently, if and only if
it has continuous trace; see \cite[Theorem 2.4]{akemann_central_1979} and
Theorem \ref{Theorem:outer-unital} below. The following is the analogue for
automorphisms of Lemma \ref{Lemma:phantom-derivation-order-1}.

\begin{lemma}
\label{Lemma:phantom-automorphism-1}Let $A$ be a separable \emph{unital}\
C*-algebra. Then:

\begin{enumerate}
\item $\mathrm{Ph}^{1}\mathrm{\mathrm{Aut}}(A)$ is the norm-closure of $%
\mathrm{Inn}(A)$;

\item the Polish group topology on $\mathrm{Ph}^{1}\mathrm{\mathrm{Aut}}(A)$
is induced by the norm.
\end{enumerate}
\end{lemma}

\begin{proof}
For $k\in \omega $ let $V_{k}$ be the set of inner automorphisms $\beta $ of 
$A$ such that $\beta =\mathrm{Ad}\left( u\right) $ for some $u\in U(A)$ with 
$\left\Vert u-1\right\Vert \leq 2^{-(k+1)}$. Then $\left( V_{k}\right) $ is
a basis of identity neighborhoods of $\mathrm{Inn}(A)$. Furthermore, the
closure $\overline{V}_{k}^{\mathrm{\mathrm{Aut}}(A)}$ of $V_{k}$ in $\mathrm{%
\mathrm{Aut}}(A)$ is contained in the set of automorphisms $\beta $ of $A$
such that $\left\Vert \beta -\mathrm{id}_{A} \right\Vert \leq 2^{-k}$.
Denote by $H$ the closure of $\mathrm{Inn}(A)$ with respect to the uniform
distance. Let $\left( u_{n}\right) _{n\in \omega }$ be a dense sequence in $%
U(A)$. Then we have that $\alpha \in H$ if and only if for every $k\in
\omega $ there exists $n\in \omega $ such that $\left\Vert \alpha -\mathrm{Ad%
}\left( u_{n}\right) \right\Vert \leq 2^{-k}$. This shows that $H$ is $%
\boldsymbol{\Pi }_{3}^{0}$ in $\mathrm{\mathrm{Aut}}(A)$. This implies that $%
\mathrm{Ph}^{1}\mathrm{\mathrm{Aut}}(A)\subseteq H$. Conversely, suppose
that $\alpha \in H$. Fix $\varepsilon >0$. Let $\xi >0$ be such that for
every $\lambda \in \mathbb{C}$ with $\left\vert \lambda \right\vert \leq \xi 
$ one has that $\left\vert \exp \left( i\lambda \right) -1\right\vert \leq
\varepsilon $. Let also $\eta \in (0,1)$ be such that%
\begin{equation*}
-\log(1-\eta)=\log\frac{1}{1-\eta}<\xi \text{.}
\end{equation*}%
Then there exists $u\in U(A)$ such that $\left\Vert \mathrm{id}_{A}-\mathrm{%
Ad}\left( u\right) \circ \alpha \right\Vert <\eta $. Thus, by \cite[Theorem
8.7.7]{pedersen_algebras_1979}, 
\begin{equation*}
\delta :=\mathrm{Log}\left( \mathrm{Ad}\left( u\right) \circ \alpha \right)
\end{equation*}%
is a *-derivation of $A$, which satisfies%
\begin{equation*}
\mathrm{Ad}\left( u\right) \circ \alpha =\exp \left( \delta \right) \text{.}
\end{equation*}%
By the choice of $\eta $, $\left\Vert \delta \right\Vert <\xi $. Indeed,%
\begin{equation*}
\left\Vert \delta \right\Vert =\left\Vert \mathrm{Log}\left( \mathrm{Ad}%
\left( u\right) \circ \alpha \right) \right\Vert \leq \sum_{n\geq 1}\frac{%
\left\Vert \mathrm{Ad}\left( u\right) \circ \alpha -\mathrm{id}%
_{A}\right\Vert ^{n}}{n}\leq \sum_{n\geq 1}\frac{\eta ^{n}}{n}=-\log \left(
1-\eta \right) <\xi \text{.}
\end{equation*}%
Setting%
\begin{equation*}
\beta _{t}:=\exp \left( t\delta \right)
\end{equation*}%
for $t\in \mathbb{R}$ produces a uniformly continuous one-parameter group of
automorphisms $\left( \beta _{t}\right) _{t\in \mathbb{R}}$ of $A$ \cite[%
Chapter 8]{pedersen_algebras_1979}. By \cite[Proposition 8.6.14]%
{pedersen_algebras_1979}, there exists a sequence $\left( v^{(n) }\right)
_{n\in \mathbb{N}}$ of uniformly continuous representations $%
(v_{t}^{(n)})_{t\in \mathbb{R}}$ with values in $U(A)$ such that, 
\begin{equation*}
(\mathrm{Ad}(v_{t}^{(n)}))_{n\in \mathbb{N}}
\end{equation*}%
converges to $\beta _{t}$ in $\mathrm{\mathrm{Aut}}(A)$ for every $t\in 
\mathbb{R}$, uniformly on compact subsets of $\mathbb{R}$. In particular,%
\begin{equation*}
(\mathrm{Ad}(v_{1}^{(n)}))_{n\in \mathbb{N}}
\end{equation*}%
converges to $\mathrm{Ad}\left( u\right) \circ \alpha $ in $\mathrm{\mathrm{%
Aut}}(A)$. The proofs of \cite[Proposition 8.6.14]{pedersen_algebras_1979}, 
\cite[Lemma 8.6.12]{pedersen_algebras_1979}, and \cite[Theorem 8.6.5]%
{pedersen_algebras_1979}, show that in fact%
\begin{equation*}
v_{1}^{(n)}=\mathrm{\exp }(ih^{(n)})
\end{equation*}%
for some $h^{(n)}\in A_{\mathrm{sa}}$ with 
\begin{equation*}
\Vert h^{(n)}\Vert \leq \left\Vert \delta \right\Vert <\xi \text{.}
\end{equation*}%
Thus, by the choice of $\xi $ we have that the norm-distance between $%
v_{1}^{(n)}$ and $1$ is at most $\varepsilon $. Since $\varepsilon $ was
arbitrary, this concludes the proof.
\end{proof}

\subsection{The unitary corona}

Let $A$ be a separable unital C*-algebra, realized as a C*-subalgebra of $%
B\left( H\right) $. As before, define $M$ to be the $\sigma $-weak closure
of $A$ in $B\left( H\right) $. Define $ZU(A)$ to be the center of $U(A)$,
which is equal to $U(A)\cap ZM=U(ZA)$. Define%
\begin{equation*}
\rho U(A):=U(A)/ZU(A)
\end{equation*}%
and%
\begin{equation*}
\rho U(M):=U(M)/ZU(M)\text{,}
\end{equation*}%
both endowed with the quotient topology induced by the norm topology. We can
identify $\rho U(A)$ with a subgroup of $\rho U(M)$. Let 
\begin{equation*}
\delta U(A)\subseteq \rho U(M)
\end{equation*}%
be the norm-closure of $\rho U(A)$, which is a Polish group, and consider
the group with a Polish cover 
\begin{equation*}
\partial U(A):=\frac{\delta U(A)}{\rho U(A)}\text{.}
\end{equation*}%
For $u\in U(M)$ and $X\subseteq M$ we define%
\begin{equation*}
d\left( u,X\right) =\inf_{x\in X}\left\Vert u-x\right\Vert \text{.}
\end{equation*}%
Note that, if $B$ is a unital C*-algebra and $u\in U(B)$, then%
\begin{equation*}
d\left( u,ZB\right) \leq d\left( u,U(ZB)\right) \leq 2d\left( u,ZB\right) 
\text{;}
\end{equation*}%
see \cite[Lemma 5.1]{kadison_derivations_1967}. We will use this for $B=A$
and for $B=M$.

\begin{proposition}
\label{Proposition:unitary-corona}Let $A\subseteq B\left( H\right) $ be a
concrete unital C*-algebra, and denote by $M$ its $w^{\ast }$-closure. The
continuous homomorphism 
\begin{equation*}
\Theta :\delta U(A)\rightarrow \mathrm{Ph}^{1}\mathrm{\mathrm{Aut}}(A)
\end{equation*}%
\begin{equation*}
u\mapsto \mathrm{Ad}\left( u\right)
\end{equation*}%
induces an iso-airomorphism%
\begin{equation*}
\partial U(A)\rightarrow \mathrm{Ph}^{1}\mathrm{Out}(A)
\end{equation*}%
and, in particular,%
\begin{equation*}
\ell _{\mathrm{Ph}}\left( \partial U(A)\right) =\ell _{\mathrm{Ph}}\left( 
\mathrm{Ph}^{1}\mathrm{Out}(A)\right) \text{.}
\end{equation*}
\end{proposition}

\begin{proof}
If $u\in \delta U(A)$ then for every $\varepsilon >0$ there exist $w\in ZU(M)
$ and $v\in U(A)$ such that $\left\Vert u-vw\right\Vert <\varepsilon $.
Thus, we have that $\left\Vert \mathrm{Ad}\left( u\right) -\mathrm{Ad}\left(
v\right) \right\Vert \leq 2\varepsilon $ and $\mathrm{Ad}\left( u\right) \in 
\mathrm{Ph}^{1}\mathrm{\mathrm{Aut}}(A)$. The kernel of $\mathrm{Ad}$ on $%
U(M)$ is $ZU(M)$, so $\Theta $ is injective, and clearly $\Theta [\rho U(A)]=%
\mathrm{Inn}(A)$.

It remains to prove that $\Theta $ is surjective. If $\alpha \in \mathrm{Ph}%
^{1}\mathrm{Aut}(A)$, choose $u\in U(A)$ such that $\left\Vert \mathrm{Ad}%
\left( u\right) \circ \alpha -\mathrm{id}_{A} \right\Vert <2$. By \cite[%
Theorem 8.7.7]{pedersen_algebras_1979}, $\mathrm{Ad}\left( u\right) \circ
\alpha =\exp \left( \delta \right) $ for some *-derivation $\delta $ of $A$.
As recalled at the beginning of Section \ref{Section:derivation}, $\delta $
extends to a $\sigma $-weakly continuous derivation of $M$, and \cite[%
Theorem 3.1]{kadison_derivations_1967} provides $a\in M_{\mathrm{sa}}$ with $%
\delta =\mathrm{ad}\left( ia\right) $. Since $\exp \left( \mathrm{ad}\left(
ia\right) \right) =\mathrm{Ad}\left( \exp \left( ia\right) \right) $, the
unitary $w:=\exp \left( ia\right) \in U(M)$ satisfies $\mathrm{Ad}\left(
u\right) \circ \alpha =\mathrm{Ad}\left( w\right) $.

We claim that $wZU(M)\in \delta U(A)$. Fix $\varepsilon >0$. Since $\mathrm{%
Ad}\left( w\right) \in \mathrm{Ph}^{1}\mathrm{\mathrm{Aut}}(A)$, there
exists $v\in U(A)$ such that $\left\Vert \mathrm{Ad}\left( w\right) -\mathrm{%
Ad}\left( v\right) \right\Vert <\varepsilon /2$. Set $x:=v^{\ast }w\in U(M) $%
. Since%
\begin{equation*}
\mathrm{Ad}\left( x\right) -\mathrm{id}_{A} =\mathrm{Ad}\left( v^{\ast
}\right) \circ \left( \mathrm{Ad}\left( w\right) -\mathrm{Ad}\left( v\right)
\right)
\end{equation*}%
and $\mathrm{Ad}\left( v^{\ast }\right) $ is isometric, we obtain that%
\begin{equation*}
\sup_{a\in \mathrm{\mathrm{Ball}}\left( A\right) }\left\Vert
xa-ax\right\Vert =\left\Vert \mathrm{Ad}\left( x\right) -\mathrm{id}_{A}
\right\Vert =\left\Vert \mathrm{Ad}\left( w\right) -\mathrm{Ad}\left(
v\right) \right\Vert <\varepsilon /2\text{,}
\end{equation*}%
where the first equality holds since $\left\Vert xax^{\ast }-a\right\Vert
=\left\Vert xa-ax\right\Vert $ for every $a\in A$. By the Kaplansky Density
Theorem \cite[I.9.1.3]{blackadar_operator_2006}, the left-hand side is equal
to the norm of $\mathrm{ad}\left( x\right) $ regarded as a derivation of $M$%
. Hence, by the formula for the norm of a derivation on a von Neumann
algebra \cite{zsido_norm_1973},%
\begin{equation*}
2d\left( x,ZM\right) =\left\Vert \mathrm{ad}\left( x\right) \right\Vert <
\varepsilon /2\text{.}
\end{equation*}%
By the inequality above applied to $B=M$, this gives $d\left( x,ZU(M)\right)
\leq 2d\left( x,ZM\right) <\varepsilon /2$, and hence there exists $z\in
ZU(M)$ such that%
\begin{equation*}
\left\Vert w-vz\right\Vert =\left\Vert x-z\right\Vert <\varepsilon \text{.}
\end{equation*}%
As $\varepsilon >0$ was arbitrary, this shows that $wZU(M)$ belongs to the
norm-closure $\delta U(A)$ of $\rho U(A)$. Hence $u^{\ast }wZU(M)\in \delta
U(A)$ and $\Theta (u^{\ast }wZU(M))=\alpha $, proving surjectivity. Since $%
\Theta $ is a continuous bijective homomorphism between Polish groups, its
inverse is continuous. Thus $\Theta $ and $\Theta ^{-1}$ are continuous
homomorphism lifts of the induced bijection and its inverse. The induced map
is therefore an iso-airomorphism, and the equality of phantom lengths
follows from Subsection \ref{Subsection:homomorphisms}.
\end{proof}

\subsection{Phantom length one}

Let $A$ be a C*-algebra. If a *-derivation $\delta $ is of the form $\mathrm{%
ad}\left( a\right) $ for some $a\in iM\left( A\right) _{\mathrm{sa}}$, then 
\begin{equation*}
\exp \left( \delta \right) =\mathrm{Ad}\left( \exp \left( a\right) \right)
\end{equation*}
is an inner *-automorphism of $A$. The converse holds in certain cases.

\begin{proposition}
\label{Proposition:derivation-automorphism}Let $A\subseteq B\left( H\right) $
be a concrete C*-algebra acting nondegenerately, and let $M$ be its $w^{\ast
}$-closure. Suppose that the center $ZM$ of $M$ is contained in $M(A)$. Then
a *-derivation $\delta $ of $A$ satisfying $\left\Vert \delta \right\Vert
<2\pi $ is of the form $\mathrm{ad}\left( a\right) $ for some $a\in iM\left(
A\right) _{\mathrm{sa}}$ if and only if $\mathrm{\exp }\left( \delta \right) 
$ is an inner automorphism.
\end{proposition}

\begin{proof}
The proof is similar to that of \cite[Theorem 4.6]{kadison_derivations_1967}%
. Let $U(A)$ be the unitary group of $M(A)$. Let $\delta $ be a *-derivation
of $A$ such that $\mathrm{\exp }\left( \delta \right) $ is an inner
automorphism. Since $\left\Vert \delta \right\Vert <2\pi $, there exists a
self-adjoint element $b$ of $M$ such that $\delta =\mathrm{ad}\left(
ib\right) $ and $\left\Vert b\right\Vert <\pi $, whence $\exp \left( \delta
\right) =\mathrm{Ad}\left( \exp \left( ib\right) \right) $. If $\exp \left(
\delta \right) $ is inner, there exists a unitary multiplier $u$ of $A$ such
that $\exp \left( \delta \right) =\mathrm{Ad}\left( u\right) $. This shows
that%
\begin{equation*}
z:=\mathrm{\exp }\left( ib\right) u^{\ast }
\end{equation*}%
belongs to $U\left( ZM\right) \subseteq U(A)$. Thus, 
\begin{equation*}
\exp \left( ib\right) \in U(A)\text{.}
\end{equation*}%
Since $\left\Vert b\right\Vert <\pi $ we have that $-1$ does not belong to
the spectrum of $\exp \left( ib\right) $. Hence, by continuous functional
calculus we have $b\in M(A)$, concluding the proof. The other implication is
obvious.
\end{proof}

\begin{corollary}
Let $A$ be a separable \emph{primitive} C*-algebra. Then a *-derivation $%
\delta $ of $A$ satisfying $\left\Vert \delta \right\Vert <2\pi $ is of the
form $\mathrm{ad}\left( a\right) $ for some $a\in iM\left( A\right) _{%
\mathrm{sa}}$ if and only if $\mathrm{\exp }\left( \delta \right) $ is inner.
\end{corollary}

From \cite[Theorem 5.3 and Corollary 5.5]{kadison_derivations_1967} one has
the following:

\begin{proposition}
\label{Proposition:KLR}Let $A$ be a separable \emph{unital} C*-algebra.

\begin{enumerate}
\item $\ell _{\mathrm{Aut}}\left( A\right) \leq 1\Leftrightarrow \ell _{%
\mathrm{aut}}\left( A\right) \leq 1$;

\item if $A\subseteq B\left( H\right) $ acts nondegenerately, has $w^{\ast }$%
-closure $M$, and the center $ZM$ of $M$ is contained in $A$, then $\ell _{%
\mathrm{aut}}\left( A\right) \leq 1$.
\end{enumerate}
\end{proposition}

\begin{proof}
(1) is the content of the equivalence of (i) and (iv) in \cite[Theorem 5.3]%
{kadison_derivations_1967}.

(2) follows from \cite[Corollary 5.5]{kadison_derivations_1967}.
\end{proof}

\subsection{From derivations to automorphisms\label%
{Subsection:from-derivations-to-automorphisms}}

Let $A$ be a separable unital C*-algebra. We now construct a pseudomorphism $%
\mathrm{Ph}^{1}\mathrm{out}(A)\rightarrow \mathrm{Ph}^{1}\mathrm{Out}(A)$ in
the sense of Definition \ref{Definition:pseudomorphism}. Recall that $%
\mathrm{Ph}^{1}\mathrm{Aut}(A)$ is the norm-closure of $\mathrm{Inn}(A)$ in $%
\mathrm{\mathrm{Aut}}(A)$ by Lemma \ref{Lemma:phantom-automorphism-1}, and
that $\mathrm{Ph}^{1}\mathrm{aut}(A)$ is the norm-closure of $\mathrm{inn}%
(A) $ in $\mathrm{aut}(A)$ by Lemma \ref{Lemma:phantom-derivation-order-1};
in particular, both these groups are endowed with the topology induced by
the operator norm.

Set%
\begin{equation*}
r_{0}:=2^{-1}\log 2
\end{equation*}%
and define%
\begin{equation*}
W:=\mathrm{\mathrm{Ball}}_{r_{0}}\left( \mathrm{Ph}^{1}\mathrm{aut}\left(
A\right) \right) \text{, }U:=\mathrm{\mathrm{Ball}}_{r_{0}/2}\left( \mathrm{%
inn}\left( A\right) \right) \text{, }\varphi :=\exp |_{W}\text{.}
\end{equation*}%
Thus, $W$ is an open symmetric zero neighborhood in $\mathrm{Ph}^{1}\mathrm{%
aut}\left( A\right) $, while $U$ is an open zero neighborhood in $\mathrm{inn%
}\left( A\right) $. Since the inclusion $\mathrm{inn}\left( A\right)
\rightarrow \mathrm{aut}\left( A\right) $ has norm at most $2$, we have that%
\begin{equation*}
U\subseteq W\cap \mathrm{inn}\left( A\right) \text{.}
\end{equation*}%
As $\exp \left( \mathrm{ad}\left( ih\right) \right) =\mathrm{Ad}\left( \exp
\left( ih\right) \right) $ for $h\in A_{\mathrm{sa}}$, and $\exp $ is
norm-continuous, $\varphi $ maps $W\cap \mathrm{inn}\left( A\right) $ to $%
\mathrm{Inn}\left( A\right) $ and $W$ to $\mathrm{Ph}^{1}\mathrm{\mathrm{Aut}%
}\left( A\right) $. Furthermore, $\varphi $ is Lipschitz, as%
\begin{equation*}
\left\Vert \exp \left( \delta \right) -\exp \left( \delta ^{\prime }\right)
\right\Vert \leq \left\Vert \delta -\delta ^{\prime }\right\Vert \exp \left(
\max \left\{ \left\Vert \delta \right\Vert ,\left\Vert \delta ^{\prime
}\right\Vert \right\} \right) \leq \sqrt{2}\left\Vert \delta -\delta
^{\prime }\right\Vert
\end{equation*}%
for $\delta ,\delta ^{\prime }\in W$. In particular, $\varphi $ is $\left( 
\mathrm{Ph}^{1}\mathrm{aut}\left( A\right) ,\mathrm{Ph}^{1}\mathrm{\mathrm{%
Aut}}\left( A\right) \right) $-uniformly continuous, as required in
Definition \ref{Definition:local-uniformly-continuous}. Likewise, $\varphi
|_{U}$ is $\left( \mathrm{inn}\left( A\right) ,\mathrm{Inn}\left( A\right)
\right) $-uniformly continuous. As $r_{0}<\log 2$, the function $\varphi $
is moreover injective; see Lemma \ref{Lemma:exp-homeomorphism} below, where
this is proved in the form in which it will be used.

The following argument was suggested by Ilja Gogi\'{c}. Define 
\begin{equation*}
\sigma :\mathrm{inn}(A)\times \mathrm{aut}(A)\rightarrow \mathrm{Inn}(A)
\end{equation*}%
\begin{equation*}
\left( \mathrm{ad}\left( ih\right) ,\delta \right) \mapsto \mathrm{Ad}\left(
u_{1}\right)
\end{equation*}%
where:

\begin{enumerate}
\item $h\in A_{\mathrm{sa}}$;

\item $\alpha _{t}=\exp \left( t\delta \right) $ for $t\in \mathbb{R}$;

\item $u_{t}$ for $t\in \mathbb{R}$ is the solution of%
\begin{equation*}
\frac{d}{dt}u_{t}=u_{t}\alpha _{t}\left( ih\right)
\end{equation*}%
subject to $u_{0}=1$.
\end{enumerate}

\begin{proposition}
\label{Proposition:from-derivation-to-automorphism}Let $A\subseteq B\left(
H\right) $ be a separable unital C*-algebra with $w^{\ast }$-closure $M$.
Then:

\begin{enumerate}
\item $\mathrm{Ad}\left( u_{1}\right) $ is well-defined and independent of
the choice of the self-adjoint element $h$;

\item for every $\delta \in \mathrm{aut}(A)$, the function 
\begin{equation*}
\mathrm{inn}\left( A\right) \rightarrow \mathrm{Inn}\left( A\right) \text{, }%
\mathrm{ad}\left( ih\right) \mapsto \sigma \left( \mathrm{ad}\left(
ih\right) ,\delta \right)
\end{equation*}%
is uniformly continuous on $U$ with modulus independent of $\delta $;

\item for $\delta ,\delta ^{\prime }\in \mathrm{aut}(A)$, $h\in A_{\mathrm{sa%
}}$, and $u_{1}\in U(A)$ such that%
\begin{equation*}
\sigma \left( \mathrm{ad}\left( ih\right) ,\delta \right) =\mathrm{Ad}\left(
u_{1}\right) \text{\quad and\quad }\delta ^{\prime }=\mathrm{ad}\left(
ih\right) +\delta \text{,}
\end{equation*}%
one has%
\begin{equation*}
\mathrm{\exp }\left( \delta ^{\prime }\right) =\mathrm{Ad}\left(
u_{1}\right) \circ \exp \left( \delta \right) \text{;}
\end{equation*}

\item $\left( W,U,\varphi \right) :\mathrm{Ph}^{1}\mathrm{out}\left(
A\right) \rightarrow \mathrm{Ph}^{1}\mathrm{Out}\left( A\right) $ is a
pseudomorphism.
\end{enumerate}
\end{proposition}

\begin{proof}
The differential equation 
\begin{equation*}
\frac{d}{dt}u_{t}=u_{t}\alpha _{t}\left( ih\right)
\end{equation*}%
subject to $u_{0}=1$ has a unique solution for all $t\in \mathbb{R}$.
Indeed, \cite[Section IV.1, Proposition 1.1]{lang_fundamentals_1999} gives
local existence and uniqueness, and the identity $\left\Vert \alpha
_{t}\left( ih\right) \right\Vert =\left\Vert h\right\Vert $ allows the local
solution to be continued stepwise across every bounded time interval. The
resulting global solution produces a uniformly continuous unitary cocycle $%
\left( u_{t}\right) $ in $A$. Furthermore, if $h,h^{\prime }$ are two
self-adjoint elements of $A$ with $\mathrm{ad}\left( ih\right) =\mathrm{ad}%
\left( ih^{\prime }\right) $, then $h-h^{\prime }\in ZA$. Thus, if $\left(
u_{t}\right) ,\left( u_{t}^{\prime }\right) $ are the corresponding
uniformly continuous unitary cocycles in $A$, then $\left( u_{t}^{\ast
}u_{t}^{\prime }\right) $ is a uniformly continuous unitary cocycle in $ZA$.
Thus, $\mathrm{Ad}\left( u_{t}\right) =\mathrm{Ad}\left( u_{t}^{\prime
}\right) $ for every $t\in \mathbb{R}$ and in particular $\mathrm{Ad}\left(
u_{1}\right) =\mathrm{Ad}\left( u_{1}^{\prime }\right) $. This shows that $%
\sigma $ is well-defined.\ If $\varepsilon >0$ and $h,h^{\prime }$ are
self-adjoint elements of $A$ such that%
\begin{equation*}
\left\Vert \mathrm{ad}\left( i\left( h-h^{\prime }\right) \right)
\right\Vert _{\mathrm{inn}(A)}=\left\Vert \mathrm{ad}\left( ih\right) -%
\mathrm{ad}\left( ih^{\prime }\right) \right\Vert _{\mathrm{inn}(A)}\leq
\varepsilon \text{,}
\end{equation*}%
then, for every $\eta >0$, by the definition of the quotient norm on $%
\mathrm{inn}(A)$ one can choose representatives satisfying $\left\Vert
h-h^{\prime }\right\Vert <\varepsilon +\eta $. When $\mathrm{ad}(ih)$ and $%
\mathrm{ad}(ih^{\prime })$ belong to $U$, these representatives can be
chosen uniformly bounded. Since every $\alpha _{t}$ is isometric, the
Banach-space Gr\"{o}nwall estimate for two differential equations \cite[%
Theorem 2.1]{howard_gronwall_2025} gives a modulus independent of $\delta
\in \mathrm{aut}(A)$; see also \cite[Section IV.1, Corollary 1.2]%
{lang_fundamentals_1999} and \cite[Section 1]{deimling_ordinary_1977}.

We claim that for every $t\in \mathbb{R}$,%
\begin{equation}
\alpha _{t}^{\prime }=\mathrm{\exp }\left( t\delta ^{\prime }\right) =%
\mathrm{Ad}\left( u_{t}\right) \circ \alpha _{t}\text{.}
\label{eq:identity-differential}
\end{equation}%
For $t=0$ the identity holds. Notice that%
\begin{equation*}
\frac{d}{dt}\alpha _{t}=\alpha _{t}\delta =\delta \alpha _{t}\text{.}
\end{equation*}%
The derivative with respect to $t$ of the left-hand side of %
\eqref{identity-differential} yields%
\begin{equation*}
\exp \left( t\delta ^{\prime }\right) \delta ^{\prime }=\alpha _{t}^{\prime
}\delta ^{\prime }\text{.}
\end{equation*}%
For $a\in A$, we let $L\left( a\right) $ and $R\left( a\right) $ be the left
and right multiplication operators. The derivative of the right-hand side
gives%
\begin{eqnarray*}
&&\frac{d\,\mathrm{Ad}\left( u_{t}\right)}{dt}\,\alpha _{t}+\mathrm{Ad}%
\left( u_{t}\right) \frac{d\alpha _{t}}{dt} \\
&=&L\left( u_{t}\alpha _{t}\left( ih\right) \right) R\left( u_{t}^{\ast
}\right) \alpha _{t}-L\left( u_{t}\right) R\left( \alpha _{t}\left(
ih\right) u_{t}^{\ast }\right) \alpha _{t}+\mathrm{Ad}\left( u_{t}\right)
(\alpha _{t}\delta ) \\
&=&\mathrm{Ad}\left( u_{t}\right) \alpha _{t}\left( L\left( ih\right)
-R\left( ih\right) +\delta \right) \\
&=&\mathrm{Ad}\left( u_{t}\right) \alpha _{t}\left( \mathrm{ad}\left(
ih\right) +\delta \right) =\left( \mathrm{Ad}\left( u_{t}\right) \circ
\alpha _{t}\right) \delta ^{\prime }\text{.}
\end{eqnarray*}%
Thus, both the left-hand side and the right-hand side of %
\eqref{identity-differential} are solutions to the same differential
equation in $L(A)$ with the same initial value. Local uniqueness \cite[%
Section IV.1, Proposition 1.1]{lang_fundamentals_1999} shows that the set of
times at which the two solutions agree is open; it is closed by continuity
and contains $0$. Hence it is all of $\mathbb{R}$, and we obtain the desired
equality. In particular, for $t=1$ we obtain $\alpha _{1}^{\prime }=\mathrm{%
Ad}\left( u_{1}\right) \circ \alpha _{1}$.

Finally, (4) follows from (1), (2), and (3), together with the observations
preceding the statement. Indeed, $W$ is an open symmetric zero neighborhood
in $\mathrm{Ph}^{1}\mathrm{aut}\left( A\right) $, $U$ is an open zero
neighborhood in $\mathrm{inn}\left( A\right) $ with $U\subseteq W\cap 
\mathrm{inn}\left( A\right) $, the function $\varphi $ is $\left( \mathrm{Ph}%
^{1}\mathrm{aut}\left( A\right) ,\mathrm{Ph}^{1}\mathrm{\mathrm{Aut}}\left(
A\right) \right) $-uniformly continuous and maps $W\cap \mathrm{inn}\left(
A\right) $ to $\mathrm{Inn}\left( A\right) $, and%
\begin{equation*}
\varphi \left( 0\right) =\mathrm{id}_{A} \text{, }\varphi \left( -\delta
\right) =\varphi \left( \delta \right) ^{-1}
\end{equation*}%
for $\delta \in W$. Thus $\left( W,U,\varphi \right) $ is a local uniformly
continuous map in the sense of Definition \ref%
{Definition:local-uniformly-continuous}, and $\sigma $ witnesses that it is
a pseudomorphism in the sense of Definition \ref{Definition:pseudomorphism}
by (2) and (3).
\end{proof}

\section{A dichotomy for derivations and automorphisms\label%
{Section:dichotomy}}

In this section we obtain a dichotomy for the complexity of the relation of
unitary equivalence of automorphisms and multiplier equivalence of
*-derivations of a given C*-algebra, building on previous work of Elliott 
\cite{elliott_some_1977}, Akemann--Pedersen \cite{akemann_central_1979}, and
the author \cite{lupini_unitary_2014}. In the proof of such a dichotomy, we
will apply recent results of Shani \cite[Theorem 1.1]{shani_generic_2024}
from descriptive set theory.

\subsection{Only inner derivations}

A major result in the study of derivations is the result from \cite%
{elliott_some_1977,akemann_central_1979} characterizing C*-algebras with
only inner derivations.

\begin{theorem}[Akemann--Elliott--Pedersen]
\label{Theorem:AEP}Let $A$ be a separable C*-algebra. The following
assertions are equivalent:

\begin{enumerate}
\item every derivation of $A$ is inner (in the multiplier algebra);

\item every summable central sequence of $A_{+}$ is trivial;

\item $A=B\oplus C$ where $B$ has continuous trace and $C$ is a direct sum
of simple C*-algebras.
\end{enumerate}
\end{theorem}

In item (3), ``direct sum'' means the restricted (that is, $c_{0}$-) direct
sum, as in \cite[Theorem 3.9]{akemann_central_1979}. Item (1) is the
condition, denoted (B) in \cite[Theorem 3.9]{akemann_central_1979}, that
every derivation of $M(A)$ be inner; the two formulations agree because
every derivation of $A$ extends uniquely to a derivation of $M(A)$, and
conversely every derivation of $M(A)$ restricts to one of $A$, compatibly
with innerness in $M(A)$.

We will give a \emph{complexity-theoretic }refinement of such a
characterization, in the context of Borel complexity of sets and equivalence
relations \cite{gao_invariant_2009}. Recall that $E_{0}$ denotes the
relation of \emph{tail equivalence }on the space $2^{\mathbb{N}}$ of binary
sequences. The relation $E_{0}^{\mathbb{N}}$ on $\left( 2^{\mathbb{N}%
}\right) ^{\mathbb{N}}$ is the countable product of copies of $E_{0}$. It
can thus be seen as the relation of \emph{tail equivalence of countably many
binary sequences}. In terms of such a relation, we obtain the following
Dichotomy Theorem for the relations of multiplier equivalence of derivations
and unitary equivalence of \emph{derivable} automorphisms. Recall that an
automorphism of $A$ is derivable if it is of the form $\exp \left( \delta
\right) $ for some *-derivation $\delta $ \cite[Section 8.7]%
{pedersen_algebras_1979}.

The set of derivable automorphisms need not be Borel in $\mathrm{Aut}(A)$.
Accordingly, when we say that an equivalence relation is Borel reducible to
unitary equivalence of derivable automorphisms, we mean that there is a
Borel map into the Polish space $\mathrm{Aut}(A)$, endowed with the topology
of pointwise convergence, whose range consists of derivable automorphisms
and which is a reduction to unitary equivalence. Thus, this terminology does
not presuppose a standard Borel structure on the set of all derivable
automorphisms.

\begin{theorem}
\label{Theorem:main}Let $A$ be a separable C*-algebra. Let $E$ be either the
relation of multiplier equivalence of (contractive) *-derivations of $A$, or
the relation of unitary equivalence of derivable automorphisms of $A$. Then
either $E$ is trivial, or the relation $E_{0}^{\mathbb{N}}$ of tail
equivalence of countably many binary sequences is Borel reducible to $E$.
\end{theorem}

In the proof of Theorem \ref{Theorem:main} we will provide several
equivalent characterizations of separable C*-algebras with outer
derivations; see Theorem \ref{Theorem:main-long}.\ One of these is Property
APE, a variation of Property AEP introduced in \cite{lupini_unitary_2014}. A
major ingredient in the proof of Theorem \ref{Theorem:main} is Shani's
Dichotomy Theorem for the relation $E_{0}^{\mathbb{N}}$ \cite[Theorem 1.1]%
{shani_generic_2024}.

\subsection{Derivation and automorphism lengths}

In the \emph{unital} case, one can rephrase Theorem \ref{Theorem:main} in
terms of derivation and automorphism length. The following result
characterizing the separable unital C*-algebras $A$ of automorphism length $%
0 $, i.e., such that $\mathrm{Inn}(A)$ is closed in $\mathrm{Aut}(A)$, is a
consequence of \cite[Theorem 0.8]{raeburn_crossed_1988}, \cite[Theorem 1.2]%
{lupini_unitary_2014}, and \cite[Theorem 6.3.1]{gao_invariant_2009}; see
also \cite[Theorem 3.1]{phillips_outer_1987}.

\begin{theorem}
\label{Theorem:outer-unital}Let $A$ be a separable unital C*-algebra. The
following assertions are equivalent:

\begin{enumerate}
\item $A$ has continuous trace;

\item $A$ has automorphism length $0$, i.e., $\mathrm{Inn}(A)$ is closed in $%
\mathrm{\mathrm{Aut}}(A)$;

\item the relation of unitary equivalence of automorphisms of $A$ is
classifiable by countable structures;

\item $E_{0}$ is not Borel reducible to the relation of unitary equivalence
of automorphisms of $A$.
\end{enumerate}
\end{theorem}

We extend this by characterizing the separable unital C*-algebras $A$ whose
automorphism length is at most $1/2$ or at most $1$, and ruling out $%
1/2+\varepsilon $ as a possible value of the automorphism length.

\begin{theorem}
\label{Theorem:main-unital}Let $A$ be a separable unital C*-algebra. The
following assertions are equivalent:

\begin{enumerate}
\item $\mathrm{out}(A)=0$, i.e., every derivation of $A$ is inner;

\item $A$ has automorphism length at most $1/2$, i.e., $\mathrm{Inn}(A)$ is $%
\boldsymbol{\Sigma }_{2}^{0}$ in $\mathrm{\mathrm{Aut}}(A)$;

\item $\mathrm{Inn}(A)$ is $\boldsymbol{\Sigma }_{3}^{0}$ in $\mathrm{%
\mathrm{Aut}}(A)$;

\item $E_{0}^{\mathbb{N}}$ is not Borel reducible to the coset relation of $%
\mathrm{Inn}(A)$ in $\mathrm{\mathrm{Aut}}(A)$.
\end{enumerate}

Furthermore, the following assertions are equivalent:

\renewcommand{\labelenumi}{(\alph{enumi})}

\begin{enumerate}
\item $A$ has derivation length at most $1$, i.e., $\mathrm{inn}(A)$ is
norm-closed in $\mathrm{aut}(A)$;

\item $\mathrm{Inn}(A)$ is norm-closed in $\mathrm{\mathrm{Aut}}(A)$;

\item $A$ has automorphism length at most $1$, i.e., $\mathrm{Inn}(A)$ is $%
\boldsymbol{\Pi }_{3}^{0}$ in $\mathrm{Aut}(A)$.
\end{enumerate}
\end{theorem}

\begin{proof}
(1)$\Rightarrow $(2). If $\alpha $ is an automorphism of $A$ such that $%
\left\Vert \alpha -\mathrm{id}_{A}\right\Vert <2$, then $\alpha $ is \emph{%
derivable }\cite[Theorem 8.7.7]{pedersen_algebras_1979}. Thus, there exists
a *-derivation $\delta $ of $A$ such that $\alpha =\exp \left( \delta
\right) $. Since by assumption $\delta $ is inner, also $\alpha $ must be
inner. Let $\left( u_{n}\right) $ be a dense sequence in $U(A)$.\ Then we
have that an automorphism $\alpha $ of $A$ is inner if and only if there
exists $n\in \mathbb{N}$ such that $\left\Vert \mathrm{Ad}\left(
u_{n}\right) -\alpha \right\Vert \leq 1$.\ This shows that $\mathrm{Inn}(A)$
is $\boldsymbol{\Sigma }_{2}^{0}$ in $\mathrm{Aut}(A)$.

(2)$\Rightarrow $(3). Obvious.

(3)$\Rightarrow $(4). Suppose, toward a contradiction, that $E_{0}^{\mathbb{N%
}}$ is Borel reducible to the coset relation of $\mathrm{Inn}\left( A\right) 
$ in $\mathrm{\mathrm{Aut}}\left( A\right) $. By \cite[Proposition 4.14]%
{lupini_looking_2024}, $E_{0}^{\mathbb{N}}$ is the coset relation of a
Polish subgroup of exact complexity class $\boldsymbol{\Pi }_{3}^{0}$. Thus,
by \cite[Theorem 3.3(3)]{lupini_complexity_2025}, it is not potentially $%
\boldsymbol{\Sigma }_{3}^{0}$. If $\mathrm{Inn}(A)$ were $\boldsymbol{\Sigma 
}_{3}^{0}$ in $\mathrm{\mathrm{Aut}}(A)$, then its coset relation would be
potentially $\boldsymbol{\Sigma }_{3}^{0}$, again by \cite[Theorem 3.3(3)]%
{lupini_complexity_2025}. This contradicts the fact that potential $%
\boldsymbol{\Sigma }_{3}^{0}$ complexity is preserved downward under Borel
reducibility.

(4)$\Rightarrow $(1). If $A$ has an outer derivation, then by Theorem \ref%
{Theorem:main}, $E_{0}^{\mathbb{N}}$ is Borel reducible to the relation of
unitary equivalence of automorphisms of $A$.

The equivalence of (a) and (b) is the content of Proposition \ref%
{Proposition:KLR}, i.e., the equivalence of (i) and (iv) in \cite[Theorem 5.3%
]{kadison_derivations_1967}. The equivalence of (b) and (c) follows from
Lemma \ref{Lemma:phantom-automorphism-1}.
\end{proof}

\begin{corollary}
\label{Corollary:automorphism-length}The automorphism length of a separable
unital C*-algebra $A$ cannot be equal to $1/2+\varepsilon $. Furthermore%
\begin{equation*}
\ell _{\mathrm{\mathrm{Aut}}}\left( A\right) \leq 1\Rightarrow \ell _{%
\mathrm{aut}}\left( A\right) \leq \ell _{\mathrm{\mathrm{Aut}}}\left(
A\right) \text{.}
\end{equation*}
\end{corollary}

\begin{proof}
By \cite[Theorem 1.1]{lupini_complexity_2025}, if $\mathrm{Inn}(A)$ is $%
\boldsymbol{\Pi }_{3}^{0}$ and not closed, then its complexity class is
either $\boldsymbol{\Sigma }_{2}^{0}$, or $D(\boldsymbol{\Pi }_{2}^{0})$, or 
$\boldsymbol{\Pi }_{3}^{0}$, corresponding to the values $1/2$, $%
1/2+\varepsilon $, and $1$ for the automorphism length. By Theorem \ref%
{Theorem:main-unital}, such a complexity class cannot be $D(\boldsymbol{\Pi }%
_{2}^{0})$.

If $\ell _{\mathrm{\mathrm{Aut}}}\left( A\right) \leq 1/2+\varepsilon $ then 
$\ell _{\mathrm{\mathrm{Aut}}}\left( A\right) \leq 1/2$ as just established.
This implies $\mathrm{out}\left( A\right) =0$ and $\ell _{\mathrm{aut}%
}\left( A\right) =0$ by Theorem \ref{Theorem:main-unital}. If $\ell _{%
\mathrm{\mathrm{Aut}}}\left( A\right) \leq 1$ then $\ell _{\mathrm{aut}%
}\left( A\right) \leq 1$ by Proposition \ref{Proposition:KLR}.
\end{proof}

The conclusion of Corollary \ref{Corollary:automorphism-length} does not
hold if $A$ is not unital. A construction of a separable continuous trace
C*-algebra such that $\mathrm{Inn}(A)$ has complexity class $D(\boldsymbol{%
\Pi }_{2}^{0})$ in $\mathrm{\mathrm{Aut}}(A)$ is described in work in
preparation \cite{bergfalk_definable_2026}.

\subsection{Property APE}

\emph{Property AEP} for a separable C*-algebra $A$ was introduced in \cite[%
Definition 4.4]{lupini_unitary_2014}. This notion is inspired by the work of
Elliott and Akemann--Pedersen \cite{elliott_some_1977,akemann_central_1979}
characterizing C*-algebras with only inner derivations. It is asserted in 
\cite[Theorem 4.5]{lupini_unitary_2014} that Property AEP is equivalent to
the property of having an outer derivation. A correction needed in the proof
there is explained in Remark \ref{Remark:AEP+} and is incorporated in the
proof below.

Property AEP was defined in terms of the relation of $\ell _{1}$-equivalence
of sequences. It was motivated by the goal of showing that, for C*-algebras
with outer derivations, the relation of unitary equivalence of automorphisms
of C*-algebras is not classifiable by countable structures. We now introduce
a variant, which we term Property Akemann--Pedersen--Elliott (APE).\ 

Suppose that $A$ is a separable C*-algebra and $\boldsymbol{a}=\left(
a_{n}\right) $ is a dense sequence in $\mathrm{\mathrm{Ball}}(A)$. Suppose
that $\boldsymbol{x}:=\left( x_{n}\right) $ is a sequence of pairwise
orthogonal positive contractions of $A$ such that, for every $n\in \mathbb{N}
$ and $i\leq n$,%
\begin{equation*}
\left\Vert \left[ x_{n},a_{i}\right] \right\Vert \leq 2^{-n}\text{.}
\end{equation*}%
As remarked in \cite[Theorem 4.5]{lupini_unitary_2014}, if $\boldsymbol{t}%
\in \left( 0,1\right) ^{\mathbb{N}}$, then the series%
\begin{equation*}
\sum_{n\in \mathbb{N}}t_{n}x_{n}
\end{equation*}%
converges in the strong operator topology to a self-adjoint element $x_{%
\boldsymbol{t}}$ of $A^{\ast \ast }$.

For $n\in \mathbb{N}$, let $\alpha _{\boldsymbol{t},n}$ be the inner
automorphism associated with the unitary%
\begin{equation*}
\exp (i\sum_{k\leq n}t_{k}x_{k})\text{.}
\end{equation*}%
Then the sequence $\left( \alpha _{\boldsymbol{t},n}\right) _{n\in \mathbb{N}%
}$ converges to an automorphism $\alpha _{\boldsymbol{t}}$ of $A$. This is a 
\emph{derivable }automorphism, associated with the *-derivation $\delta _{%
\boldsymbol{t}}:=\mathrm{ad}(ix_{\boldsymbol{t}})$ induced by the positive
contraction $x_{\boldsymbol{t}}\in A^{\ast \ast }$.

Since the $x_{n}$'s are pairwise orthogonal, if $\boldsymbol{t},\boldsymbol{s%
}\in \left( 0,1\right) ^{\mathbb{N}}$ are such that $\boldsymbol{t}-%
\boldsymbol{s}\in c_{0}$, then the series%
\begin{equation*}
\sum_{n\in \mathbb{N}}\left( t_{n}-s_{n}\right) x_{n}
\end{equation*}%
converges to a self-adjoint contraction $x_{\boldsymbol{t}-\boldsymbol{s}%
}\in A$. Then we have that%
\begin{equation*}
\delta _{\boldsymbol{t}}=\mathrm{ad}\left( ix_{\boldsymbol{t}-\boldsymbol{s}%
}\right) +\delta _{\boldsymbol{s}}
\end{equation*}%
and%
\begin{equation*}
\alpha _{\boldsymbol{t}}=\mathrm{Ad}\left( \mathrm{\exp }\left( ix_{%
\boldsymbol{t}-\boldsymbol{s}}\right) \right) \circ \alpha _{\boldsymbol{s}}%
\text{.}
\end{equation*}

\begin{definition}
\label{Definition:APE-sequence}Let $A$ be a separable C*-algebra, and $%
\boldsymbol{a}=\left( a_{n}\right) $ be a dense sequence in $\mathrm{\mathrm{%
Ball}}(A)$. A sequence $\boldsymbol{x}=\left( x_{n}\right) $ as above is an
APE sequence for $\boldsymbol{a}$ if for every uncountable well-separated $%
Y\subseteq \left( 0,1\right) ^{\mathbb{N}}$ there exist $\boldsymbol{s},%
\boldsymbol{t}\in Y$ such that $\alpha _{\boldsymbol{t}}$ and $\alpha _{%
\boldsymbol{s}}$ are not unitarily equivalent.
\end{definition}

\begin{definition}
\label{Definition:APE-property}Let $A$ be a separable C*-algebra.\ Then $A$
satisfies Property APE if for every dense sequence $\boldsymbol{a}$ in $%
\mathrm{\mathrm{Ball}}(A)$ there exists an APE sequence for $\boldsymbol{a}$.
\end{definition}

We now consider the implications that Property APE has for the relations of
multiplier equivalence of derivations and unitary equivalence of
automorphisms.

The set $\exp \left( \mathrm{aut}(A)\right) $ of derivable automorphisms of $%
A$ need not be Borel in $\mathrm{\mathrm{Aut}}(A)$. We therefore work
instead with the following two relations, both of which live on \emph{Polish 
}spaces: the relation $E_{\mathrm{aut}}$ of multiplier equivalence on $%
\mathrm{\mathrm{Ball}}\left( \mathrm{aut}(A)\right) $, and the relation $E_{%
\mathrm{Aut}}$ on $2\mathrm{\mathrm{Ball}}\left( \mathrm{aut}(A)\right) $
defined by%
\begin{equation*}
\delta \mathrel{E_{\mathrm{Aut}}}\delta ^{\prime }\Leftrightarrow \exp
\left( \delta \right) \circ \exp \left( \delta ^{\prime }\right) ^{-1}\in 
\mathrm{Inn}(A)\text{.}
\end{equation*}%
The relation $E_{\mathrm{aut}}$ is analytic by Proposition \ref%
{Proposition:multiplier-equivalence-derivations-analytic}. The relation $E_{%
\mathrm{Aut}}$ is Borel. Indeed, on a uniformly norm-bounded set of
derivations, the exponential map is continuous from the strong operator
topology to the topology of pointwise convergence: this follows directly
from its power series, which converges uniformly on norm-bounded sets.
Therefore the map $\left( \delta ,\delta ^{\prime }\right)\mapsto \exp
\left( \delta \right) \circ \exp \left( \delta ^{\prime }\right) ^{-1}$ is
continuous on $2\mathrm{Ball}(\mathrm{aut}(A))^{2}$. The subgroup $\mathrm{%
Inn}(A)$ is Borel in $\mathrm{Aut}(A)$, as observed above, so its preimage $%
E_{\mathrm{Aut}}$ is Borel.

\begin{lemma}
\label{Lemma:APE-orthogonality}Let $A$ be a separable C*-algebra. Consider
the following conditions:

\begin{enumerate}
\item $A$ satisfies Property APE;

\item neither $E_{\mathrm{aut}}$ nor $E_{\mathrm{Aut}}$ is $c_{0}$%
-orthogonal;

\item the relation $E_{0}^{\mathbb{N}}$ of tail equivalence of countably
many binary sequences is Borel reducible to $E_{\mathrm{aut}}$ and to $E_{%
\mathrm{Aut}}$, and hence to the relation of multiplier equivalence of
(contractive) *-derivations of $A$ and to the relation of unitary
equivalence of derivable automorphisms of $A$ in the sense specified before
Theorem \ref{Theorem:main}.
\end{enumerate}

Then (1)$\Rightarrow $(2)$\Rightarrow $(3).
\end{lemma}

\begin{proof}
(1)$\Rightarrow $(2). Let $\boldsymbol{a}$ be a dense sequence in $\mathrm{%
\mathrm{Ball}}(A)$ and let $\boldsymbol{x}$ be an APE sequence for $%
\boldsymbol{a}$. Since $x_{\boldsymbol{t}}$ is a positive contraction, $%
\left\Vert \delta _{\boldsymbol{t}}\right\Vert \leq 2$, so that $\boldsymbol{%
t}\mapsto \delta _{\boldsymbol{t}}$ maps $\left( 0,1\right) ^{\mathbb{N}}$
into $2\mathrm{\mathrm{Ball}}\left( \mathrm{aut}(A)\right) $, and $%
\boldsymbol{t}\mapsto \delta _{\boldsymbol{t}}/2$ maps it into $\mathrm{%
\mathrm{Ball}}\left( \mathrm{aut}(A)\right) $; both maps are Borel, and both
induce functions on the quotients by $c_{0}$, since $\boldsymbol{t}-%
\boldsymbol{s}\in c_{0}$ implies $\delta _{\boldsymbol{t}}=\mathrm{ad}\left(
ix_{\boldsymbol{t}-\boldsymbol{s}}\right) +\delta _{\boldsymbol{s}}$ with $%
x_{\boldsymbol{t}-\boldsymbol{s}}\in A$ and $\alpha _{\boldsymbol{t}}=%
\mathrm{Ad}\left( \exp \left( ix_{\boldsymbol{t}-\boldsymbol{s}}\right)
\right) \circ \alpha _{\boldsymbol{s}}$, as observed before Definition \ref%
{Definition:APE-sequence}. By that definition, for every uncountable
well-separated $Y\subseteq \left( 0,1\right) ^{\mathbb{N}}$ there are $%
\boldsymbol{s},\boldsymbol{t}\in Y$ with $\alpha _{\boldsymbol{t}}$ and $%
\alpha _{\boldsymbol{s}}$ not unitarily equivalent, that is, with $\neg
\left( \delta _{\boldsymbol{t}}\mathrel{E_{\mathrm{Aut}}}\delta _{%
\boldsymbol{s}}\right) $. For such $\boldsymbol{s},\boldsymbol{t}$ we also
have $\neg \left( \delta _{\boldsymbol{t}}/2\mathrel{E_{\mathrm{aut}}}\delta
_{\boldsymbol{s}}/2\right) $: otherwise $\delta _{\boldsymbol{t}}-\delta _{%
\boldsymbol{s}}=\mathrm{ad}\left( ia\right) $ for some $a\in M(A)_{\mathrm{sa%
}}$, since $\mathrm{inn}(M(A))$ is a linear subspace, and then, the $x_{n}$
being pairwise orthogonal, $\delta _{\boldsymbol{t}}$ and $\delta _{%
\boldsymbol{s}}$ commute, whence%
\begin{equation*}
\alpha _{\boldsymbol{t}}\circ \alpha _{\boldsymbol{s}}^{-1}=\exp \left(
\delta _{\boldsymbol{t}}-\delta _{\boldsymbol{s}}\right) =\exp \left( 
\mathrm{ad}\left( ia\right) \right) =\mathrm{Ad}\left( \exp \left( ia\right)
\right) \text{,}
\end{equation*}%
which is inner. Thus neither $E_{\mathrm{aut}}$ nor $E_{\mathrm{Aut}}$ is $%
c_{0}$-orthogonal.

(2)$\Rightarrow $(3). The relation $E_{\mathrm{aut}}$ is analytic and $E_{%
\mathrm{Aut}}$ is Borel, hence both are analytic equivalence relations on
Polish spaces. Thus Corollary \ref{Corollary:Shani} applies to each of them
and yields Borel reductions of $E_{0}^{\mathbb{N}}$ to $E_{\mathrm{aut}}$
and to $E_{\mathrm{Aut}}$. The first is a Borel reduction of $E_{0}^{\mathbb{%
N}}$ to multiplier equivalence of contractive *-derivations. Composing the
second with the exponential map produces a Borel function into $\mathrm{Aut}%
(A)$ whose range is contained in $\exp \left( 2\mathrm{\mathrm{Ball}}\left( 
\mathrm{aut}(A)\right) \right) $. It is a Borel reduction of $E_{0}^{\mathbb{%
N}}$ to unitary equivalence of derivable automorphisms in the sense
specified before Theorem \ref{Theorem:main}.
\end{proof}

\subsection{The primitive case}

Recall that an ideal of a C*-algebra is \emph{primitive} if it is the kernel
of an irreducible representation. A C*-algebra is primitive if the trivial
ideal is primitive, i.e., if it has a faithful irreducible representation.

\begin{lemma}
\label{Lemma:primitive}Let $A$ be a separable primitive nonsimple
infinite-dimensional C*-algebra. Then $A$ satisfies Property APE.
\end{lemma}

\begin{proof}
We follow the ideas of \cite[Lemma 4.6]{lupini_unitary_2014}. Since $A$ is
primitive, we can fix a faithful irreducible representation $\pi $ of $A$ on
a separable Hilbert space. Then $\pi $ has a unique extension to a $\sigma $%
-weakly continuous representation of $A^{\ast \ast }$ on $H$, which we still
denote by $\pi $.

Let $\boldsymbol{a}=\left( a_{n}\right) $ be a dense sequence in $\mathrm{%
\mathrm{Ball}}(A)$, and fix a strictly positive contraction $b_{0}$ of $A$;
such an element exists because every separable C*-algebra is $\sigma $%
-unital \cite[Proposition II.4.2.4]{blackadar_operator_2006}. Following the
argument of \cite[Lemma 3.2]{akemann_central_1979}, one can define a
sequence $\boldsymbol{x}=\left( x_{n}\right) $ of pairwise orthogonal
positive contractions and an $\varepsilon >0$ such that for every $k,n\in 
\mathbb{N}$ such that $k\leq n$, one has that 
\begin{equation*}
\left\Vert x_{n}b_{0}\right\Vert >\varepsilon
\end{equation*}%
and 
\begin{equation*}
\left\Vert \left[ x_{n},a_{k}\right] \right\Vert <2^{-n}\text{.}
\end{equation*}%
We claim that $\boldsymbol{x}$ is an APE sequence for $\boldsymbol{a}$.\ Fix 
$\theta >0$ and a $\theta $-separated uncountable subset $Y\subseteq \left(
0,1\right) ^{\mathbb{N}}$. Suppose, toward a contradiction, that for every $%
\boldsymbol{s},\boldsymbol{t}\in Y$, the automorphisms $\alpha _{\boldsymbol{%
t}}$ and $\alpha _{\boldsymbol{s}}$ are unitarily equivalent. This means
that $\alpha _{\boldsymbol{t}}\circ \alpha _{\boldsymbol{s}}^{-1}$ is inner.
The latter automorphism is equal to $\exp \left( \delta _{\boldsymbol{t}%
}-\delta _{\boldsymbol{s}}\right) $. Since the elements $x_{n}$ are pairwise
orthogonal positive contractions,%
\begin{equation*}
\left\Vert \delta _{\boldsymbol{t}}-\delta _{\boldsymbol{s}}\right\Vert
=\left\Vert \mathrm{ad}\left( i\left( x_{\boldsymbol{t}}-x_{\boldsymbol{s}%
}\right) \right) \right\Vert \leq 2\left\Vert x_{\boldsymbol{t}}-x_{%
\boldsymbol{s}}\right\Vert \leq 2<2\pi \text{.}
\end{equation*}%
Furthermore, as $\pi $ is irreducible, the $w^{\ast }$-closure of $\pi
\lbrack A]$ is $B\left( H\right) $, whose center consists of the scalar
multiples of the identity and is therefore contained in the multiplier
algebra of $\pi \lbrack A]$. Thus, Proposition \ref%
{Proposition:derivation-automorphism} applies and gives that 
\begin{equation*}
\delta _{\boldsymbol{t}}-\delta _{\boldsymbol{s}}=\mathrm{ad}\left( i\left(
x_{\boldsymbol{t}}-x_{\boldsymbol{s}}\right) \right)
\end{equation*}%
is inner in the multiplier algebra. Thus, there exists $z_{\boldsymbol{t},%
\boldsymbol{s}}$ in the center of $A^{\ast \ast }$ such that 
\begin{equation*}
x_{\boldsymbol{t}}-x_{\boldsymbol{s}}+z_{\boldsymbol{t},\boldsymbol{s}}\in
M(A)\text{.}
\end{equation*}%
Since $\pi $ is irreducible, it is in particular a \emph{factor
representation}. Thus, the relative commutant of $\pi \lbrack A]$ in $%
B\left( H\right) $ consists only of scalar multiples of the identity. In
particular, we have that $\pi \left( z_{\boldsymbol{t},\boldsymbol{s}%
}\right) \in \mathbb{C}1$ for every $\boldsymbol{s},\boldsymbol{t}\in Y$.
This implies that%
\begin{equation*}
\pi \left( x_{\boldsymbol{t}}-x_{\boldsymbol{s}}\right) \in \pi \lbrack M(A)]
\end{equation*}%
and hence%
\begin{equation*}
\pi \left( b_{0}\left( x_{\boldsymbol{t}}-x_{\boldsymbol{s}}\right) \right)
\in \pi \lbrack A]\text{.}
\end{equation*}%
Fix $\boldsymbol{s}\in Y$. For distinct $\boldsymbol{t},\boldsymbol{t}%
^{\prime }\in Y$ there exists $m\in \mathbb{N}$ such that%
\begin{equation*}
\left\vert t_{m}-t_{m}^{\prime }\right\vert \geq \theta \text{.}
\end{equation*}%
Since the elements $x_{n}$ are pairwise orthogonal, we have that%
\begin{equation*}
\left( x_{\boldsymbol{t}}-x_{\boldsymbol{t}^{\prime }}\right) x_{m}=\left(
t_{m}-t_{m}^{\prime }\right) x_{m}^{2}\text{.}
\end{equation*}%
Moreover, since $b_{0}$ and $x_{m}$ are self-adjoint contractions,%
\begin{equation*}
\left\Vert b_{0}x_{m}^{2}\right\Vert =\left\Vert x_{m}^{2}b_{0}\right\Vert
\geq \left\Vert b_{0}x_{m}^{2}b_{0}\right\Vert =\left\Vert
x_{m}b_{0}\right\Vert ^{2}>\varepsilon ^{2}\text{.}
\end{equation*}%
As $b_{0}x_{m}^{2}$ belongs to $A$, on which $\pi $ is isometric, and $\pi $
is contractive on $A^{\ast \ast }$, we conclude that%
\begin{eqnarray*}
\left\Vert \pi \left( b_{0}\left( x_{\boldsymbol{t}}-x_{\boldsymbol{s}%
}\right) \right) -\pi \left( b_{0}\left( x_{\boldsymbol{t}^{\prime }}-x_{%
\boldsymbol{s}}\right) \right) \right\Vert &=&\left\Vert \pi \left(
b_{0}\left( x_{\boldsymbol{t}}-x_{\boldsymbol{t}^{\prime }}\right) \right)
\right\Vert \\
&\geq &\left\Vert \pi \left( b_{0}\left( x_{\boldsymbol{t}}-x_{\boldsymbol{t}%
^{\prime }}\right) \right) \pi \left( x_{m}\right) \right\Vert \\
&=&\left\vert t_{m}-t_{m}^{\prime }\right\vert \left\Vert \pi \left(
b_{0}x_{m}^{2}\right) \right\Vert \\
&=&\left\vert t_{m}-t_{m}^{\prime }\right\vert \left\Vert
b_{0}x_{m}^{2}\right\Vert \geq \varepsilon ^{2}\theta \text{.}
\end{eqnarray*}%
As $Y$ is uncountable, this contradicts the separability of $\pi \lbrack A]$.
\end{proof}

\subsection{Lifting}

A crucial step in the proof of the main result of \cite{lupini_unitary_2014}
consists of showing that Property AEP is \emph{liftable}. Here we prove that
the same holds for Property APE.

\begin{lemma}
\label{Lemma:lifting}Let $A$ and $B$ be separable C*-algebras, and $\pi
:A\rightarrow B$ be a surjective *-homomorphism. If $B$ satisfies Property
APE, then so does $A$.
\end{lemma}

\begin{proof}
Denote by $J$ the kernel of $\pi $. Let $\boldsymbol{a}=\left( a_{n}\right) $
be a dense sequence in $\mathrm{\mathrm{Ball}}(A)$. Observe that $\pi \left( 
\boldsymbol{a}\right) =\left( \pi \left( a_{n}\right) \right) $ is a dense
sequence in $\mathrm{\mathrm{Ball}}\left( B\right) $. Since by hypothesis $B$
satisfies Property APE, we can consider an APE sequence $\boldsymbol{y}%
=\left( y_{n}\right) $ for $\pi \left( \boldsymbol{a}\right) $. Replacing $%
y_{n}$ by $c_{n}y_{n}$, where $c_{n}=1-2^{-n-2}$, we can assume that 
\begin{equation*}
\left\Vert \left[ y_{n},\pi \left( a_{i}\right) \right] \right\Vert <2^{-n}
\end{equation*}%
for every $i\leq n$. This replacement preserves the property of being an APE
sequence, since the coordinatewise map 
\begin{equation*}
\left( t_{n}\right) \mapsto \left( c_{n}t_{n}\right)
\end{equation*}%
maps well-separated sets to well-separated sets, and the automorphisms
defined by the rescaled sequence with parameters $\left( t_{n}\right) $ are
the automorphisms defined by the original sequence with parameters $\left(
c_{n}t_{n}\right) $. Then by \cite[Lemma 10.1.12]{loring_lifting_1997} we
can find a sequence $\boldsymbol{z}=\left( z_{n}\right) $ of pairwise
orthogonal positive contractions in $A$ such that $\pi \left( z_{n}\right)
=y_{n}$ for every $n\in \omega $.

Let $\left( e_{k}\right) $ be an \emph{increasing approximate unit }for $J$
that is \emph{quasicentral} in $A$. This means that $\left( e_{k}\right) $
is an increasing sequence of positive elements of $J$ such that, for every $%
b\in J$ and $a\in A$, the sequence $\left( [e_{k},a]\right) _{k\in \mathbb{N}%
}$ converges to $0$ and the sequences $\left( e_{k}b\right) $ and $\left(
be_{k}\right) $ converge to $b$.

By \cite[II.5.1.1]{blackadar_operator_2006}, for every $x\in A$ one has that%
\begin{equation*}
\left\Vert \pi \left( x\right) \right\Vert =\lim_{k}\left\Vert x\left(
1-e_{k}\right) \right\Vert =\inf_{k}\left\Vert x\left( 1-e_{k}\right)
\right\Vert \text{.}
\end{equation*}%
Fix $n\in \mathbb{N}$ and apply this to $x:=z_{n}a_{i}-a_{i}z_{n}$ for $%
i\leq n$. Considering that $\left( e_{k}\right) $ is quasicentral in $A$,
for every $n\in \mathbb{N}$ one can find $k_{n}\in \mathbb{N}$ such that,
setting%
\begin{equation*}
x_{n}:=z_{n}^{\frac{1}{2}}\left( 1-e_{k_{n}}\right) z_{n}^{\frac{1}{2}}
\end{equation*}%
one has that%
\begin{equation*}
\left\Vert \lbrack x_{n},a_{i}]\right\Vert <2^{-n}
\end{equation*}%
for $i\leq n$. Suppose that $\boldsymbol{s},\boldsymbol{t}\in \left(
0,1\right) ^{\mathbb{N}}$ are such that the corresponding automorphisms $%
\alpha _{\boldsymbol{s}}$ and $\alpha _{\boldsymbol{t}}$ defined in terms of
the sequence $\boldsymbol{x}$ are unitarily equivalent. Then the
automorphisms $\beta _{\boldsymbol{s}}$ and $\beta _{\boldsymbol{t}}$
obtained from $\boldsymbol{s}$ and $\boldsymbol{t}$ in terms of the sequence 
$\boldsymbol{y}$ are also unitarily equivalent. Since $\boldsymbol{y}$ is an
APE\ sequence for $\pi (\boldsymbol{a})$, this implies that $\boldsymbol{x}$
is an APE sequence for $\boldsymbol{a}$.
\end{proof}

\begin{corollary}
\label{Corollary:notT1}If $A$ is a separable C*-algebra whose primitive
spectrum $\mathrm{Prim}(A)$ is not $T_{1}$, then $A$ has Property APE.
\end{corollary}

\begin{proof}
Since the primitive spectrum of $A$ is not $T_{1}$, $A$ has a nonsimple
primitive quotient. The conclusion follows from Lemma \ref{Lemma:primitive}
and liftability of Property APE (Lemma \ref{Lemma:lifting}).
\end{proof}

\subsection{A stronger property}

As in the proof of the main result of \cite{lupini_unitary_2014}, it is
convenient to consider a stronger property, called Property AEP$^{+}$. The
following is a corrected version of \cite[Definition 4.9]%
{lupini_unitary_2014}.

\begin{definition}
\label{Definition:AEP+}A C*-algebra $A$ has Property AEP$^{+}$ if there
exist a sequence $\left( \pi _{n}\right) $ of irreducible representations of 
$A$, a positive contraction $b_{0}$ of $A$, and a central sequence $\left(
x_{n}\right) $ of pairwise orthogonal positive contractions of $A$, such
that:

\begin{enumerate}
\item for every $\lambda \in \mathbb{C}$, 
\begin{equation*}
\liminf_{n\rightarrow \infty}\left\Vert \pi _{n}\left( \left( x_{n}-\lambda
\right) b_{0}\right) \right\Vert >0;
\end{equation*}

\item for every pair of \emph{distinct} natural numbers $n$ and $m$, $%
x_{n}\in \ker\left( \pi _{m}\right) $.
\end{enumerate}
\end{definition}

\begin{remark}
\label{Remark:AEP+}There is a typo in \cite[Definition 4.9]%
{lupini_unitary_2014}, where condition (1) is stated as the requirement that
the corresponding sequence \textquotedblleft does not converge to
zero\textquotedblright. The stronger condition in Definition \ref%
{Definition:AEP+} is the one needed in the proof of \cite[Proposition 4.11]%
{lupini_unitary_2014}. There is also a missing hypothesis in \cite[Lemma 4.13%
]{lupini_unitary_2014}: the spectrum condition alone does not imply Property
AEP$^{+}$ (for example, consider $C(\omega +1)$). Corrected versions of the
source results, with the uniform estimate required by Definition \ref%
{Definition:AEP+}, are given in Lemmas \ref{Lemma:quotient} and \ref%
{Lemma:omega} below.
\end{remark}

With the caveat from Remark \ref{Remark:AEP+}, the following perturbation
lemma is obtained in \cite[Lemma 4.10]{lupini_unitary_2014}.

\begin{lemma}
\label{Lemma:perturbation}For every $\varepsilon >0$ there exists $\rho >0$
such that for every C*-algebra $A$ and every pair of positive contractions $%
x,b$ of $A$ and $\mu \in \mathbb{C}$ such that 
\begin{equation*}
\left\Vert b\right\Vert \geq \varepsilon
\end{equation*}%
and%
\begin{equation*}
\left\Vert \left( \mathrm{\exp }\left( ix\right) -\mu \right) b\right\Vert
\leq \rho
\end{equation*}%
one must have that%
\begin{equation*}
\left\Vert \left( x-\lambda \right) b\right\Vert \leq \varepsilon
\end{equation*}%
for some $\lambda \in \mathbb{C}$.
\end{lemma}

With this lemma, we can show that Property AEP$^{+}$ is indeed stronger than
Property APE.

\begin{lemma}
\label{Lemma:APE+}Let $A$ be a separable C*-algebra. If $A$ satisfies
Property AEP$^{+}$, then it satisfies Property APE.
\end{lemma}

\begin{proof}
We follow the outline of the proof of \cite[Proposition 4.11]%
{lupini_unitary_2014}. Let $\left( \pi _{n}\right) $ be a sequence of
irreducible representations of $A$, $b_{0}$ be a positive contraction of $A$%
, and $\left( x_{n}\right) $ be a sequence of pairwise orthogonal positive
contractions of $A$ that witness Property AEP$^{+}$ for $A$. Rescaling $%
b_{0} $, we may assume that it has norm $1$. For $n\in \mathbb{N}$ and $%
\mu\in\mathbb{C}$, set 
\begin{equation*}
f_{n}(\mu):=\left\Vert \pi _{n}\left( \left( x_{n}-\mu \right) b_{0}\right)
\right\Vert\text{.}
\end{equation*}
We claim that 
\begin{equation*}
\liminf_{n\rightarrow\infty}\inf_{\mu\in\mathbb{C}}f_{n}(\mu)>0.
\end{equation*}
Indeed, condition (1) with $\mu=0$ gives, after discarding finitely many
terms, a constant $c>0$ such that 
\begin{equation*}
\left\Vert \pi _{n}(b_{0})\right\Vert \geq f_{n}(0)\geq c.
\end{equation*}
If the claim failed, there would be a strictly increasing sequence $\left(
n_{k}\right) $ and scalars $\mu _{k}$ such that $f_{n_{k}}(\mu _{k})$
converges to zero. Since 
\begin{equation*}
\left\vert \mu _{k}\right\vert \left\Vert \pi _{n_{k}}(b_{0})\right\Vert
\leq f_{n_{k}}(\mu _{k})+\left\Vert \pi _{n_{k}}(x_{n_{k}}b_{0}) \right\Vert
,
\end{equation*}
the sequence $\left( \mu _{k}\right) $ is bounded. After passing to a
subsequence, suppose that $\mu _{k}$ converges to $\mu$. The estimate 
\begin{equation*}
\left\vert f_{n}(\mu)-f_{n}(\nu)\right\vert \leq \left\vert
\mu-\nu\right\vert \left\Vert \pi _{n}(b_{0})\right\Vert \leq \left\vert
\mu-\nu\right\vert
\end{equation*}
then gives $f_{n_{k}}(\mu)\rightarrow0$, contradicting Definition \ref%
{Definition:AEP+}(1).

Let $\boldsymbol{a}=\left( a_{n}\right) $ be a dense sequence in $\mathrm{%
Ball}(A)$. After discarding finitely many terms and then passing to a common
subsequence of the pairs $\left( \pi _{n},x_{n}\right) $, we may therefore
assume that there exists an $\varepsilon >0$ such that for every $\mu \in 
\mathbb{C}$ and $i,n\in \mathbb{N}$ with $i\leq n$,%
\begin{equation*}
\left\Vert \pi _{n}\left( \left( x_{n}-\mu \right) b_{0}\right) \right\Vert
\geq \varepsilon
\end{equation*}%
and%
\begin{equation*}
\left\Vert \lbrack x_{n},a_{i}]\right\Vert <2^{-n}\text{.}
\end{equation*}%
Thus, for $\delta >0$, $\mu \in \mathbb{C}$, $n\in \mathbb{N}$, and $t\in
\lbrack \delta ,1)$ we have%
\begin{equation*}
\left\Vert \pi _{n}\left( \left( tx_{n}-\mu \right) b_{0}\right) \right\Vert
\geq \delta \varepsilon \text{.}
\end{equation*}%
We claim that $\boldsymbol{x}$ is an APE sequence for $\boldsymbol{a}$.
Indeed, let $Y$ be an uncountable $\delta $-separated subset of $\left(
0,1\right) ^{\mathbb{N}}$ for some $\delta >0$. Suppose, toward a
contradiction, that $\alpha _{\boldsymbol{s}}$ and $\alpha _{\boldsymbol{t}}$
are unitarily equivalent for all $\boldsymbol{s}$ and $\boldsymbol{t}$ in $Y$%
. Taking a skew-adjoint logarithm of the resulting central unitary, we
obtain an element $z_{\boldsymbol{t},\boldsymbol{s}}$ in the center of $%
A^{\ast \ast }$ such that 
\begin{equation*}
\exp \left( i\left( x_{\boldsymbol{t}}-x_{\boldsymbol{s}}\right) +z_{%
\boldsymbol{t},\boldsymbol{s}}\right) \in M(A)
\end{equation*}%
and hence%
\begin{equation*}
y_{\boldsymbol{t},\boldsymbol{s}}:=\exp \left( i\left( x_{\boldsymbol{t}}-x_{%
\boldsymbol{s}}\right) +z_{\boldsymbol{t},\boldsymbol{s}}\right) b_{0}\in A%
\text{.}
\end{equation*}%
Fix $\boldsymbol{s}\in Y$ and let $\boldsymbol{t}$ and $\boldsymbol{t}%
^{\prime }$ be distinct elements of $Y$.\ Then there exists $m\in \mathbb{N}$
such that $\left\vert t_{m}-t_{m}^{\prime }\right\vert \geq\delta$.
Interchanging $\boldsymbol{t}$ and $\boldsymbol{t}^{\prime }$ if necessary,
we may assume that%
\begin{equation*}
t_{m}-t_{m}^{\prime }\geq \delta \text{.}
\end{equation*}%
Since $\pi _{m}$ is an irreducible representation of $A$, there exists $\mu
_{\boldsymbol{t},\boldsymbol{t}^{\prime }}\in \mathbb{C}$ such that%
\begin{equation*}
\pi _{m}\left( \exp \left( z_{\boldsymbol{t}^{\prime },\boldsymbol{s}}-z_{%
\boldsymbol{t},\boldsymbol{s}}\right) \right) =\mu _{\boldsymbol{t},%
\boldsymbol{t}^{\prime }}1
\end{equation*}%
is a scalar multiple of the identity.

Let $\rho >0$ be obtained from $\delta \varepsilon /2$ as in Lemma \ref%
{Lemma:perturbation}. Since $\left\Vert \pi _{n}(b_{0})\right\Vert \geq
f_{n}(0)\geq\varepsilon$, the perturbation lemma, applied in $\pi _{n}[A]$,
shows that for every $t\in \lbrack \delta ,1)$, $n\in \mathbb{N}$, and $\mu
\in \mathbb{C}$,%
\begin{equation*}
\left\Vert \pi _{n}\left( \left( \exp \left( itx_{n}\right) -\mu \right)
b_{0}\right) \right\Vert \geq \rho \text{.}
\end{equation*}%
Thus, we have%
\begin{eqnarray*}
\left\Vert y_{\boldsymbol{t},\boldsymbol{s}}-y_{\boldsymbol{t}^{\prime },%
\boldsymbol{s}}\right\Vert &=&\left\Vert \left( \exp \left( i\left( x_{%
\boldsymbol{t}}-x_{\boldsymbol{t}^{\prime }}\right) \right) -\exp \left( z_{%
\boldsymbol{t}^{\prime },\boldsymbol{s}}-z_{\boldsymbol{t},\boldsymbol{s}%
}\right) \right) b_{0}\right\Vert \\
&\geq &\left\Vert \pi _{m}\left( \left( \exp \left( i\left( x_{\boldsymbol{t}%
}-x_{\boldsymbol{t}^{\prime }}\right) \right) -\mu _{\boldsymbol{t},%
\boldsymbol{t}^{\prime }}\right) b_{0}\right) \right\Vert \\
&\geq &\left\Vert \pi _{m}\left( \left( \exp \left( i\left(
t_{m}-t_{m}^{\prime }\right) x_{m}\right) -\mu _{\boldsymbol{t},\boldsymbol{t%
}^{\prime }}\right) b_{0}\right) \right\Vert \geq \rho \text{.}
\end{eqnarray*}%
As this holds for every distinct $\boldsymbol{t}$ and $\boldsymbol{t}%
^{\prime }$ in $Y$, it contradicts the separability of $A$, concluding the
proof.
\end{proof}

The following two lemmas are corrected forms of \cite[Lemmas 4.12 and 4.13]%
{lupini_unitary_2014}. The first follows from the constructions in \cite[%
Lemmas 3.6 and 3.7]{akemann_central_1979}; the second uses the construction
in the proof of \cite[Theorem 1, pp. 139--141]{elliott_some_1977}.

\begin{lemma}
\label{Lemma:quotient}Let $A$ be a separable C*-algebra with $T_{1}$
primitive spectrum $\mathrm{Prim}(A)$. Let $\left( \xi _{n}\right) $ be a
convergent sequence of pairwise distinct separated points in $\mathrm{Prim}%
(A)$, where separated is meant in the sense of \cite[Section 1]%
{dixmier_points_1961}, and let 
\begin{equation*}
F:=\left\{ \rho\in\mathrm{Prim}(A):\left( \xi _{n}\right) \text{ converges
to }\rho\right\}
\end{equation*}
be its nonempty closed set of limits. Let $I$ be the closed ideal of $A$
corresponding to $F$. If either the quotient $A/I$ does not have continuous
trace, or the multiplier algebra of $A/I$ has nontrivial center, then $A$
has Property \emph{AEP}$^{+}$.
\end{lemma}

\begin{proof}
Put $B=A/I$, and let $q:A\rightarrow B$ be the quotient map. We give the
quantitative part of the constructions in \cite[Lemmas 3.6 and 3.7]%
{akemann_central_1979}, since this is where the uniformity in the scalar is
obtained. Fix a dense sequence $\left(d_r\right)$ in $A$.

Suppose first that $B$ does not have continuous trace. By \cite[Theorem 2.4]%
{akemann_central_1979}, $B_{+}$ has a nontrivial central sequence. Lemma 3.4
of the same paper lifts it to a nontrivial central sequence $\left(
c_{j}\right)$ of positive contractions in $A$. Fix a strictly positive
contraction $b_{0}$ in $A$. The strict-topology criterion and \cite[Lemma 2.1%
]{akemann_central_1979} show that, after passing to a subsequence, there is $%
\varepsilon>0$ such that 
\begin{equation*}
\inf_{w\in Z(M(B))}\|(q(c_j)-w)q(b_0)\|>\varepsilon \qquad(j\in\mathbb{N}).
\end{equation*}
Here is the boundedness point in this argument. Since $0\leq q(c_j)\leq1$,
the infimum may be restricted to positive central contractions: applying the
pointwise metric projection $\mathbb{C}\rightarrow[0,1]$ to a central
multiplier $w$ can only decrease $\|(q(c_j)-w)q(b_0)\|$. If the displayed
separation failed, positive central contractions $w_j$ could be chosen for
which this norm tends to zero. The differences $q(c_j)-w_j$ are bounded and
self-adjoint, so the corresponding left products also tend to zero; since $%
q(b_0)$ is strictly positive, the strict-topology criterion would make the
sequence $q(c_j)$ trivial. This proves the display. Since $\widetilde
q(Z(M(A)))\subseteq Z(M(B))$, it follows that 
\begin{equation}  \label{eq:quotient-central-separation}
\left\Vert q\left(\left(c_{j}-z\right)b_{0}\right)\right\Vert >\varepsilon
\end{equation}
for every $j$ and every $z\in Z(M(A))$, where $\widetilde q:M(A)\rightarrow
M(B)$ denotes the multiplier extension. By passing to a further subsequence,
we also arrange that 
\begin{equation*}
\|[c_j,d_r]\|<2^{-j}\qquad(r\leq j).
\end{equation*}

For each $P\in\mathrm{Prim}(A)$ choose an irreducible representation $\sigma
_P$ with kernel $P$, and abbreviate $\sigma _{\xi_k}$ to $\sigma_k$. We use
the same symbols for their strict extensions to $M(A)$. We claim that for
every $j$ and $k_{0}$ there is $k>k_{0}$ such that 
\begin{equation}  \label{eq:quotient-fibre-separation}
\inf_{\lambda\in\mathbb{C}} \left\Vert\sigma
_{k}\left(\left(c_{j}-\lambda\right)b_{0}\right) \right\Vert>\frac{%
\varepsilon}{2}.
\end{equation}
Otherwise choose, for all $k>k_{0}$, a scalar $\lambda _{k}$ for which the
displayed norm is at most $\varepsilon/2$. The sequence $\left(\lambda
_{k}\right)$ is bounded: if $\rho\in F$, strict positivity gives $\|\sigma
_{\rho}(b_{0})\|>0$, and lower semicontinuity gives a uniform positive lower
bound for $\|\sigma _{k}(b_{0})\|$ eventually, while 
\begin{equation*}
|\lambda _{k}|\,\|\sigma _{k}(b_{0})\| \leq \|\sigma _{k}((c_{j}-\lambda
_{k})b_{0})\|+ \|\sigma _{k}(c_{j}b_{0})\|.
\end{equation*}
Pass to a subsequence for which $\lambda _{k}\rightarrow\lambda$. Since 
\begin{equation*}
\|\sigma_k((c_j-\lambda)b_0)\| \leq\|\sigma_k((c_j-\lambda_k)b_0)\|
+|\lambda_k-\lambda|\,\|\sigma_k(b_0)\|,
\end{equation*}
the right-hand side has limit superior at most $\varepsilon/2$. Since this
subsequence converges to every $\rho\in F$, lower semicontinuity of $%
P\mapsto\|a+P\|$ gives 
\begin{equation*}
\left\Vert\sigma _{\rho}\left(\left(c_{j}-\lambda\right)b_{0}\right)
\right\Vert\leq\frac{\varepsilon}{2}\qquad(\rho\in F).
\end{equation*}
Taking the supremum over $F=\mathrm{Prim}(B)$ contradicts %
\eqref{quotient-central-separation} with $z=\lambda1$.

Using \eqref{quotient-fibre-separation}, choose recursively a subsequence,
still denoted by $\left(\xi _{j}\right)$, for which the corresponding $%
\sigma _{j}$ satisfies that estimate with $c_{j}$. As in \cite[Lemma 3.6]%
{akemann_central_1979}, separatedness permits a further subsequence and
pairwise disjoint open sets $G_{j}\subseteq\mathrm{Prim}(A)$ with $\xi
_{j}\in G_{j}$. Let $J_{j}$ be the ideal corresponding to $G_{j}$. Because $%
\mathrm{Prim}(A)$ is $T_1$, $\xi_j$ is maximal. Hence $\sigma_j(A)$ is
simple, and $\sigma_j(J_j)=\sigma_j(A)$ because $J_j\not\subseteq\xi_j$. If $%
\left(u_{j,\alpha}\right)_\alpha$ is a quasicentral approximate unit of $J_j$%
, it follows that $\sigma_j(u_{j,\alpha})$ is an approximate unit of $%
\sigma_j(A)$. Choosing $\alpha$ sufficiently large and setting 
\begin{equation*}
x_j=u_{j,\alpha}c_ju_{j,\alpha},
\end{equation*}
we obtain positive contractions $x_j\in J_j$ such that 
\begin{equation*}
\|[x_{j},d_{r}]\|<2^{-j}\quad(r\leq j), \qquad \left\Vert\sigma
_{j}\left(\left(c_{j}-x_{j}\right)b_{0}\right) \right\Vert<\frac{\varepsilon%
}{4}.
\end{equation*}
This is the localization step in the second paragraph of the proof of \cite[%
Lemma 3.6]{akemann_central_1979}. Consequently $\left(x_{j}\right)$ is
central, the $x_{j}$'s are pairwise orthogonal, $x_{j}\in\ker(\sigma _{l})$
for $j\ne l$, and 
\begin{equation}  \label{eq:quotient-final-separation}
\inf_{\lambda\in\mathbb{C}} \left\Vert\sigma
_{j}\left(\left(x_{j}-\lambda\right)b_{0}\right) \right\Vert>\frac{%
\varepsilon}{4}\qquad(j\in\mathbb{N}).
\end{equation}

It remains to treat the case $Z(M(B))\ne\mathbb{C}1$. The multiplier-lifting
result used in \cite[Lemma 3.7]{akemann_central_1979} gives a positive
contraction $h\in M(A)$ such that $\widetilde q(h)$ is central and, for two
points $\rho _{0},\rho _{1}\in F$, its scalar values are respectively $0$
and $1$. For a strictly positive contraction $b_{0}$ put 
\begin{equation*}
c=\min\{\|\sigma _{\rho _{0}}(b_{0})\|, \|\sigma _{\rho _{1}}(b_{0})\|\}>0.
\end{equation*}
For all sufficiently large $k$ one has 
\begin{equation}  \label{eq:center-fibre-separation}
\inf_{\lambda\in\mathbb{C}} \|\sigma _{k}((h-\lambda)b_{0})\|\geq\frac c4.
\end{equation}
For if this failed, boundedness of the witnessing scalars and passage to a
convergent subsequence, followed by lower semicontinuity at both $\rho _0$
and $\rho _1$, would give 
\begin{equation*}
|\lambda|\,\|\sigma _{\rho _0}(b_{0})\|\leq c/4, \qquad
|1-\lambda|\,\|\sigma _{\rho _1}(b_{0})\|\leq c/4,
\end{equation*}
which is impossible because $\max\{|\lambda|,|1-\lambda|\}\geq1/2$.

We now arrange centrality before making the cut-downs. For $m\in\mathbb{N}$
set 
\begin{equation*}
K_m:=\bigcup_{r\leq m}\left\{P\in\mathrm{Prim}(A):
\|\sigma_P([h,d_r])\|\geq2^{-m}\right\}.
\end{equation*}
These sets are compact and disjoint from $F$, because $\widetilde q(h)$ is
central. Moreover 
\begin{equation*}
F_n:=F\cup\{\xi_k:k\geq n\}
\end{equation*}
is closed for every $n$; this is the separated-point argument in the first
paragraph of the proof of \cite[Lemma 3.7]{akemann_central_1979}. Hence, by
compactness, $K_m\cap F_n=\varnothing$ for some $n$. A diagonal passage to a
subsequence therefore gives 
\begin{equation*}
\|\sigma_k([h,d_r])\|<2^{-k}\qquad(r\leq k).
\end{equation*}
After passing to one more subsequence, choose pairwise disjoint open sets $%
G_k$ with $\xi_k\in G_k$, and let $J_k$ be their corresponding ideals.

Finally perform the two quasicentral cut-downs from the second paragraph of
the proof of \cite[Lemma 3.7]{akemann_central_1979}: first by an approximate
unit of $\ker(\sigma _k)$ and then by one of $J_k$. They produce positive
contractions $x_{k}\in J_{k}$ with 
\begin{equation*}
\|[x_{k},d_{r}]\|<2^{-k}\quad(r\leq k), \qquad \|\sigma _k((h-x_k)b_0)\|<c/8.
\end{equation*}
Here the preliminary cut-down preserves the $\sigma _k$-image, and the final
one puts $x_k$ in $J_k$. Thus the $x_k$'s are central and pairwise
orthogonal, have the required cross-kernel property, and %
\eqref{center-fibre-separation} gives 
\begin{equation*}
\inf_{k}\inf_{\lambda\in\mathbb{C}} \|\sigma _k((x_k-\lambda)b_0)\|\geq c/8.
\end{equation*}
Together with \eqref{quotient-final-separation}, this proves Property AEP$%
^{+}$ in both cases.
\end{proof}

\begin{lemma}
\label{Lemma:omega}Let $A$ be a separable C*-algebra with $\mathrm{Prim}(A)$
homeomorphic to $\omega +1$, i.e., the one-point compactification of a
countable discrete space. If $A$ has an outer derivation, then $A$ has
Property \emph{AEP}$^{+}$.
\end{lemma}

\begin{proof}
Write 
\begin{equation*}
\mathrm{Prim}(A)=\left\{ t_{n}:n\in\mathbb{N}\right\}\mathbin{\dot\cup}
\left\{ t_{\infty}\right\},
\end{equation*}
where the points $t_{n}$ are isolated. For every $s\in\mathrm{Prim}(A)$
choose an irreducible representation $\pi_s$ with kernel $s$, and abbreviate 
$\pi_{t_n}$ to $\pi_n$. By Theorem \ref{Theorem:AEP}, there is a nontrivial
summable central sequence $\left(a_{j}\right)$ of positive contractions. Fix
a strictly positive contraction $b_{0}\in A$. By the strict-topology
criterion and \cite[Lemma 2.1]{akemann_central_1979}, after passing to a
subsequence there is $\varepsilon>0$ such that 
\begin{equation}  \label{eq:omega-central-distance}
\inf\left\{\|(a_j-z)b_0\|:z\in Z(M(A))_{\mathrm{sa}},\ \|z\|\leq1\right\}
>4\varepsilon
\end{equation}
for every $j$. Otherwise central self-adjoint contractions could be chosen
so that $a_j-z_j$ converges strictly to zero, contradicting nontriviality.

We recall explicitly the fiber-localization step in the proof of \cite[%
Theorem 1, pp. 138 and 141]{elliott_some_1977}. For a positive contraction $%
a $ it is enough, in minimizing 
\begin{equation*}
\|\pi_t((a-\lambda)b_0)\|,
\end{equation*}
to take $\lambda\in[0,1]$. Suppose that every fiber admits a scalar $%
\lambda_t \in[0,1]$ satisfying $\|\pi_t((a-\lambda_t)b_0)\|\leq2\varepsilon$%
. The local upper-semicontinuity step used on page 138 gives a neighborhood $%
U_t$ on which the same scalar satisfies 
\begin{equation*}
\|\pi_s((a-\lambda_t)b_0)\|<3\varepsilon\qquad(s\in U_t).
\end{equation*}
Choose a locally finite partition of unity $\left(f_i\right)$ subordinate to
a refinement of this cover and scalars $\lambda_i\in[0,1]$ associated with
the corresponding members. By the Dauns--Hofmann theorem, 
\begin{equation*}
z:=\sum_i f_i\lambda_i
\end{equation*}
is a self-adjoint central contraction, and the fiberwise triangle inequality
gives $\|(a-z)b_0\|\leq3\varepsilon$. Applied to $a=a_j$, this contradicts %
\eqref{omega-central-distance}. Hence for every $j$ some $s_j\in \mathrm{Prim%
}(A)$ satisfies 
\begin{equation*}
\inf_{\lambda\in\mathbb{C}}\|\pi_{s_j}((a_j-\lambda)b_0)\|>2\varepsilon.
\end{equation*}
The function 
\begin{equation*}
t\longmapsto\inf_{\lambda\in[0,1]} \|\pi_t((a_j-\lambda)b_0)\|
\end{equation*}
is lower semicontinuous: the norm is jointly lower semicontinuous in $%
\left(t,\lambda\right)$ and the parameter set is compact. Since the isolated
points $t_n$ are dense, after decreasing $\varepsilon$ we may therefore
choose indices $k(j)$ such that 
\begin{equation*}
\inf_{\lambda\in\mathbb{C}}\left\Vert \pi _{k(j)}\left( \left(
a_{j}-\lambda\right)b_{0}\right) \right\Vert \geq\varepsilon
\end{equation*}
for every $j$.

Each fixed value occurs only finitely often among the indices $k(j)$.
Indeed, otherwise the image of a subsequence of $\left( a_{j}\right) $ in
the corresponding simple quotient would be a summable central sequence. It
is trivial by \cite[Lemma 3.3]{akemann_central_1979}, contradicting the
displayed estimate. Passing to a further subsequence, we may therefore
assume that the indices $k(j)$ are pairwise distinct.

Since $t_{k(j)}$ is isolated and closed, the simple quotient with primitive
spectrum $\left\{t_{k(j)}\right\}$ is a direct summand of $A$. Let $p_{j}$
be the corresponding central multiplier projection and set 
\begin{equation*}
x_{j}:=p_{j}a_{j}.
\end{equation*}
The projections $p_{j}$ are pairwise orthogonal and $\pi_{k(j)}(p_j)=1$.
Hence $\left( x_{j}\right) $ is a central sequence of pairwise orthogonal
positive contractions, $x_{j}\in\ker(\pi _{k(l)})$ when $j\neq l$, and 
\begin{equation*}
\inf_{\lambda\in\mathbb{C}}\left\Vert \pi _{k(j)}\left( \left(
x_{j}-\lambda\right)b_{0}\right) \right\Vert \geq\varepsilon.
\end{equation*}
Thus $\left( \pi _{k(j)}\right) $, $b_{0}$, and $\left( x_{j}\right) $
witness Property AEP$^{+}$.
\end{proof}

\subsection{Proof of Theorem \protect\ref{Theorem:main}}

The following result provides in particular a proof of Theorem \ref%
{Theorem:main}. It is a refinement of the Akemann--Elliott--Pedersen Theorem
characterizing C*-algebras with outer derivations; see Theorem \ref%
{Theorem:AEP}.

\begin{theorem}
\label{Theorem:main-long}Let $A$ be a separable C*-algebra. The following
assertions are equivalent:

\begin{enumerate}
\item $A$ has an outer derivation;

\item $A$ has an outer derivable automorphism;

\item $A$ has Property \emph{APE};

\item the relation of multiplier equivalence of (contractive) *-derivations
of $A$ is not $c_{0}$-orthogonal;

\item the coded relation $E_{\mathrm{Aut}}$ on $2\mathrm{\mathrm{Ball}}(%
\mathrm{aut}(A))$ is not $c_{0}$-orthogonal;

\item the relation $E_{0}^{\mathbb{N}}$ of tail equivalence of countably
many binary sequences is Borel reducible to the relation of multiplier
equivalence of (contractive) *-derivations of $A$;

\item the relation $E_{0}^{\mathbb{N}}$ is Borel reducible to the coded
relation $E_{\mathrm{Aut}}$ on $2\mathrm{\mathrm{Ball}}(\mathrm{aut}(A))$,
and hence to unitary equivalence of derivable automorphisms in the sense
specified before Theorem \ref{Theorem:main}.
\end{enumerate}
\end{theorem}

\begin{proof}
(1)$\Rightarrow $(3). Suppose that $A$ has an outer derivation. If $\mathrm{%
Prim}(A)$ is not $T_{1}$, the conclusion follows from Corollary \ref%
{Corollary:notT1}. We may therefore assume that $\mathrm{Prim}(A)$ is $T_{1}$%
.

By the contrapositive of \cite[Theorem 16]{elliott_some_1977}, either $%
\mathrm{Prim}(A)$ is not Hausdorff, or there is an ideal $J$ of $A$ such
that the quotient $B=A/J$ has $\mathrm{Prim}(B)$ homeomorphic to $\omega+1$
and $B$ has an outer derivation. If such a quotient exists, then $B$ has
Property AEP$^{+}$ by Lemma \ref{Lemma:omega}, hence Property APE by Lemma %
\ref{Lemma:APE+}; liftability (Lemma \ref{Lemma:lifting}) shows that $A$ has
Property APE%
%TCIMACRO{\TeXButton{@}{\@}}%
%BeginExpansion
\@%
%EndExpansion
. We may therefore assume that no such quotient exists, in which case $%
\mathrm{Prim}(A)$ is not Hausdorff.

It remains to consider the case in which $\mathrm{Prim}(A)$ is $T_{1}$ but
not Hausdorff. Choose distinct points $\rho _{0},\rho _{1}$ which do not
have disjoint neighborhoods. Using the density of the separated points \cite[%
Proposition 1]{dixmier_points_1961}, and passing to a subsequence as in \cite%
[Lemma 12]{elliott_some_1977}, one obtains a sequence $\left( \xi
_{n}\right) $ of pairwise distinct separated points which converges to every
point of the closed set 
\begin{equation*}
F:=\left\{ \rho\in\mathrm{Prim}(A):\left( \xi _{n}\right) \text{ converges
to }\rho\right\},
\end{equation*}
where $\rho _{0},\rho _{1}\in F$. Let $I$ be the closed ideal corresponding
to $F$. Since $F$ contains at least two points, $A/I$ is nonsimple. If $A/I$
has continuous trace, then \cite[Lemma 3.1 and Theorem 3.9]%
{akemann_central_1979} imply that the center of $M(A/I)$ is nontrivial. If $%
A/I$ does not have continuous trace, the other alternative in Lemma \ref%
{Lemma:quotient} applies. Thus, in both cases, $A$ has Property AEP$^{+}$,
and hence Property APE by Lemma \ref{Lemma:APE+}. This concludes the proof
of the implication (1)$\Rightarrow $(3).

All the other implications are either trivial or follow from Lemma \ref%
{Lemma:APE-orthogonality}.
\end{proof}

\section{Comparison\label{Section:comparison}}

We conclude with the proof of the comparison theorem between derivation
length and automorphism length of a separable unital C*-algebra; see Theorem %
\ref{Theorem:equality}.

\subsection{From derivations to automorphisms}

Let $A$ be a separable unital C*-algebra. We now establish that the
derivation length of $A$ is less than or equal to the automorphism length of 
$A$. In the proof of the following lemmas, we consider the pseudomorphism 
\begin{equation*}
\left( W,U,\varphi \right) :\mathrm{Ph}^{1}\mathrm{out}\left( A\right)
\rightarrow \mathrm{Ph}^{1}\mathrm{Out}\left( A\right)
\end{equation*}
defined in Proposition \ref{Proposition:from-derivation-to-automorphism},
where%
\begin{equation*}
W=\mathrm{\mathrm{Ball}}_{r_{0}}\left( \mathrm{Ph}^{1}\mathrm{aut}\left(
A\right) \right) \text{, }U=\mathrm{\mathrm{Ball}}_{r_{0}/2}\left( \mathrm{%
inn}\left( A\right) \right) \text{, }\varphi =\exp |_{W}\text{, }%
r_{0}=2^{-1}\log 2\text{.}
\end{equation*}%
We use the intrinsic quotient norm 
\begin{equation*}
\left\Vert \mathrm{ad}\left( ih\right) \right\Vert _{\mathrm{inn}\left(
A\right)} :=\inf_{z\in ZA_{\mathrm{sa}}}\left\Vert h-z\right\Vert
\end{equation*}%
on $\mathrm{inn}\left( A\right)$, and the complete two-sided invariant
quotient metric 
\begin{equation*}
q\left( \mathrm{Ad}\left( u\right) ,\mathrm{Ad}\left( v\right) \right)
:=\inf_{z\in U\left( ZA\right) }\left\Vert u-vz\right\Vert
\end{equation*}%
on $\mathrm{Inn}\left( A\right)$. In particular, 
\begin{equation}
\left\Vert d\right\Vert \leq 2\left\Vert d\right\Vert _{\mathrm{inn}\left(
A\right)}\quad \text{for }d\in \mathrm{inn}\left( A\right)
\label{eq:ambient-intrinsic-derivation}
\end{equation}%
and 
\begin{equation*}
\left\Vert \mathrm{Ad}\left( u\right) -\mathrm{id}_{A}\right\Vert \leq
2q\left( \mathrm{Ad}\left( u\right) ,1\right) \text{.}
\end{equation*}%
After truncation by $1$, the operator norm divided by $2$ on $\mathrm{Ph}^{1}%
\mathrm{aut}\left( A\right)$, the intrinsic norm on $\mathrm{inn}\left(
A\right)$, the uniform norm divided by $2$ on $\mathrm{Ph}^{1}\mathrm{Aut}%
\left( A\right)$, and $q$ on $\mathrm{Inn}\left( A\right)$ witness that $%
\mathrm{Ph}^{1}\mathrm{out}\left( A\right)$ and $\mathrm{Ph}^{1}\mathrm{Out}%
\left( A\right)$ are TSI groups with a Polish cover. Indeed, the conjugation
action on $\mathrm{Inn}\left( A\right)$ is isometric for $q$, since every
automorphism of $A$ acts isometrically on $A$ and maps $U\left( ZA\right) $
onto itself. This rescaling by $1/2$ serves only to normalize the metrics as
required by Definition \ref{Definition:TSI}; all norms and radii in the
estimates below refer to the operator norm on $\mathrm{Ph}^{1}\mathrm{aut}%
\left( A\right) $ and to $\left\Vert \cdot \right\Vert _{\mathrm{inn}\left(
A\right) }$ on $\mathrm{inn}\left( A\right) $. We denote the corresponding
canonical metrics from Lemma \ref{Lemma:metric} by $d_{\xi }^{\mathrm{Ph}^{1}%
\mathrm{out}\left( A\right) }$ and $d_{\xi }^{\mathrm{Ph}^{1}\mathrm{Out}%
\left( A\right) }$, respectively.

Adopt the notation from Proposition \ref{Proposition:local-comparison}. The
following lemma is the content of \cite[Theorem 4.3]%
{kadison_derivations_1967}.

\begin{lemma}
\label{Lemma:inner-inner}Let $A$ be a separable unital C*-algebra. If $%
\delta \in \mathrm{aut}\left( A\right) $ satisfies $\exp \left( t\delta
\right) \in \mathrm{Inn}\left( A\right) $ for uncountably many $t\in \mathbb{%
R}$, then $\delta \in \mathrm{inn}\left( A\right) $.
\end{lemma}

\begin{lemma}
\label{Lemma:comparison-ordinal}Let $A$ be a separable unital C*-algebra.
For $\alpha <\omega _{1}$, if $\ell _{\mathrm{Aut}}\left( A\right) \leq
1+\alpha $ then $\ell _{\mathrm{aut}}\left( A\right) \leq 1+\alpha $.
\end{lemma}

\begin{proof}
Suppose that $\delta \in W_{\alpha }\subseteq \mathrm{Ph}^{1+\alpha }\mathrm{%
aut}\left( A\right) $. Then for every $t\in \lbrack 0,1]$, we have $t\delta
\in W_{\alpha }\subseteq \mathrm{Ph}^{1+\alpha }\mathrm{aut}\left( A\right) $%
. Thus, for every $t\in \left[ 0,1\right] $, we have%
\begin{equation*}
\exp \left( t\delta \right) \in \mathrm{Ph}^{1+\alpha }\mathrm{Aut}\left(
A\right) \text{.}
\end{equation*}%
Thus, for every $t\in \left[ 0,1\right] $, $\exp \left( t\delta \right) $ is
an inner automorphism of $A$. By Lemma \ref{Lemma:inner-inner}, $\delta \in 
\mathrm{inn}\left( A\right) $. As this holds for an arbitrary $\delta \in
W_{\alpha }$ and $W_{\alpha }$ is an open subset of $\mathrm{Ph}^{1+\alpha }%
\mathrm{aut}\left( A\right) $, while $\mathrm{inn}\left( A\right) $ is dense
in every canonical phantom numerator, it follows that $\mathrm{Ph}^{1+\alpha
}\mathrm{aut}\left( A\right) =\mathrm{inn}\left( A\right) $.
\end{proof}

\begin{lemma}
\label{Lemma:comparison-sigma}Let $A$ be a separable unital C*-algebra and $%
\beta <\omega _{1}$. If 
\begin{equation*}
\ell _{\mathrm{Ph}}\left( \mathrm{Ph}^{1}\mathrm{Out}\left( A\right) \right)
\leq \beta +1/2,
\end{equation*}%
then 
\begin{equation*}
\ell _{\mathrm{Ph}}\left( \mathrm{Ph}^{1}\mathrm{out}\left( A\right) \right)
\leq \beta +1/2.
\end{equation*}
\end{lemma}

\begin{proof}
Since $\mathrm{Inn}\left( A\right) $ is $\boldsymbol{\Sigma }_{2}^{0}$ in $%
\mathrm{Ph}^{1+\beta }\mathrm{Aut}\left( A\right) $, there is an identity
neighborhood $V$ in $\mathrm{Inn}\left( A\right) $ that is \emph{closed} in $%
\mathrm{Ph}^{1+\beta }\mathrm{Aut}\left( A\right) $. By continuity of 
\begin{equation*}
\exp :\mathrm{inn}\left( A\right) \rightarrow \mathrm{Inn}\left( A\right)
\end{equation*}%
there exists $r\in \left( 0,r_{0}/2\right) $ such that 
\begin{equation*}
\mathrm{Ball}_{r}\left( \mathrm{inn}\left( A\right) \right) \subseteq U
\end{equation*}%
and%
\begin{equation*}
\exp \left( \mathrm{Ball}_{r}\left( \mathrm{inn}\left( A\right) \right)
\right) \subseteq V\text{.}
\end{equation*}

We claim that 
\begin{equation*}
\overline{\mathrm{Ball}_{r}\left( \mathrm{inn}\left( A\right) \right) }^{%
\mathrm{Ph}^{1+\beta }\mathrm{aut}\left( A\right) }\subseteq \mathrm{inn}%
\left( A\right) \text{.}
\end{equation*}%
Indeed, suppose that%
\begin{equation*}
\delta \in \overline{\mathrm{Ball}_{r}\left( \mathrm{inn}\left( A\right)
\right) }^{\mathrm{Ph}^{1+\beta }\mathrm{aut}\left( A\right) }\text{.}
\end{equation*}%
Then there exists a sequence $\left( \delta _{n}\right) $ in $\mathrm{Ball}%
_{r}\left( \mathrm{inn}\left( A\right) \right) $ such that%
\begin{equation*}
\delta _{n}\longrightarrow \delta \quad \text{in }\mathrm{Ph}^{1+\beta }%
\mathrm{aut}\left( A\right) \text{.}
\end{equation*}%
Since the inclusion $\mathrm{inn}(A)\rightarrow \mathrm{aut}(A)$ has norm at
most $2$, the choice $2r<r_{0}$ also ensures that $t\delta _{n}$ and $%
t\delta $ belong to $W_{\beta }$ for every $t\in \lbrack 0,1]$. For every
such $t$ we obtain%
\begin{equation*}
t\delta _{n}\longrightarrow t\delta \quad \text{in }\mathrm{Ph}^{1+\beta }%
\mathrm{aut}\left( A\right)
\end{equation*}%
and, for every $n\in \mathbb{N}$, 
\begin{equation*}
t\delta _{n}\in \mathrm{Ball}_{r}\left( \mathrm{inn}\left( A\right) \right) 
\text{.}
\end{equation*}%
Since $\exp |_{W_{\beta }}:W_{\beta }\rightarrow \mathrm{Ph}^{1+\beta }%
\mathrm{Aut}\left( A\right) $ is $(\mathrm{Ph}^{1+\beta }\mathrm{aut}\left(
A\right) ,\mathrm{Ph}^{1+\beta }\mathrm{Aut}\left( A\right) )$-uniformly
continuous by Proposition \ref{Proposition:local-comparison}, we obtain 
\begin{equation*}
\exp \left( t\delta _{n}\right) \longrightarrow \exp \left( t\delta \right)
\quad \text{in }\mathrm{Ph}^{1+\beta }\mathrm{Aut}\left( A\right)
\end{equation*}%
and, for every $k\in \mathbb{N}$,%
\begin{equation*}
\exp \left( t\delta _{k}\right) \in V\text{.}
\end{equation*}%
Since $V$ is closed in $\mathrm{Ph}^{1+\beta }\mathrm{Aut}\left( A\right) $,
this implies $\exp \left( t\delta \right) \in V\subseteq \mathrm{Inn}\left(
A\right) $. Since this holds for every $t\in \left[ 0,1\right] $, Lemma \ref%
{Lemma:inner-inner} implies that $\delta \in \mathrm{inn}(A)$. This
concludes the proof that%
\begin{equation*}
\overline{\mathrm{Ball}_{r}\left( \mathrm{inn}\left( A\right) \right) }^{%
\mathrm{Ph}^{1+\beta }\mathrm{aut}\left( A\right) }\subseteq \mathrm{inn}%
\left( A\right) \text{.}
\end{equation*}%
This shows that $\mathrm{inn}\left( A\right) $ has a zero neighborhood that
is closed in $\mathrm{Ph}^{1+\beta }\mathrm{aut}\left( A\right) $. This
implies that $\mathrm{inn}\left( A\right) \in \boldsymbol{\Sigma }_{2}^{0}(%
\mathrm{Ph}^{1+\beta }\mathrm{aut}\left( A\right) )$, concluding the proof.
\end{proof}

\begin{lemma}
\label{Lemma:comparison-D}Let $A$ be a separable unital C*-algebra and $%
\beta <\omega _{1}$. Suppose that $\ell _{\mathrm{Ph}}\left( \mathrm{Ph}^{1}%
\mathrm{Out}\left( A\right) \right) \leq \beta +1/2+\varepsilon $. Then $%
\ell _{\mathrm{Ph}}\left( \mathrm{Ph}^{1}\mathrm{out}\left( A\right) \right)
\leq \beta +1/2+\varepsilon $.
\end{lemma}

\begin{proof}
Consider%
\begin{equation*}
S:=\left\{ \delta \in W_{\beta }:\exp \left( \delta \right) \in \mathrm{Inn}%
\left( A\right) \right\} \text{.}
\end{equation*}%
Since 
\begin{equation*}
\exp |_{W_{\beta }}:W_{\beta }\rightarrow \mathrm{Ph}^{1+\beta }\mathrm{Aut}%
\left( A\right)
\end{equation*}%
is $(\mathrm{Ph}^{1+\beta }\mathrm{aut}\left( A\right) ,\mathrm{Ph}^{1+\beta
}\mathrm{Aut}\left( A\right) )$-uniformly continuous by Proposition \ref%
{Proposition:local-comparison}, we obtain that%
\begin{equation*}
S\in D(\boldsymbol{\Pi }_{2}^{0})(\mathrm{Ph}^{1+\beta }\mathrm{aut}\left(
A\right) )\subseteq \boldsymbol{\Sigma }_{3}^{0}(\mathrm{Ph}^{1+\beta }%
\mathrm{aut}\left( A\right) )\text{.}
\end{equation*}%
Consider now the multiplicative Polish group $\mathbb{R}_{+}$ of strictly
positive real numbers. We have a continuous action%
\begin{equation*}
\mathbb{R}_{+}\curvearrowright \mathrm{Ph}^{1+\beta }\mathrm{aut}\left(
A\right)
\end{equation*}%
defined by%
\begin{equation*}
\left( t,\delta \right) \mapsto t\delta \text{.}
\end{equation*}%
We consider the corresponding \emph{Vaught transform }defined as in \cite[%
Section 3.2]{gao_invariant_2009}. Then by \cite[Theorem 3.2.7]%
{gao_invariant_2009}, if $V=\left( 0,1\right) \subseteq \mathbb{R}_{+}$,%
\begin{equation*}
S^{\Delta V}\cap S\in \boldsymbol{\Sigma }_{3}^{0}(\mathrm{Ph}^{1+\beta }%
\mathrm{aut}\left( A\right) )\text{.}
\end{equation*}%
Then we have that%
\begin{equation*}
S^{\Delta V}\cap S=W_{\beta }\cap \mathrm{inn}\left( A\right) \text{.}
\end{equation*}%
Indeed, if $\delta \in W_{\beta }\cap \mathrm{inn}\left( A\right) $ then,
for every $t\in V$, $t\delta \in \mathrm{inn}\left( A\right) $ and hence $%
\exp \left( t\delta \right) \in \mathrm{Inn}\left( A\right) $ and $t\delta
\in S$. Since this holds for every $t\in V$, we have $\delta \in S^{\Delta
V} $. Conversely, if $\delta \in S^{\Delta V}\cap S$ then%
\begin{equation*}
\left\{ t\in V:t\delta \in S\right\} =\left\{ t\in V:\exp \left( t\delta
\right) \in \mathrm{Inn}\left( A\right) \right\}
\end{equation*}%
is nonmeager in $V$ and, in particular, uncountable. This shows that $\delta
\in \mathrm{inn}\left( A\right) $ by Lemma \ref{Lemma:inner-inner}. Since $%
W_{\beta }$ is open in $\mathrm{Ph}^{1+\beta }\mathrm{aut}\left( A\right) $, 
$W_{\beta }\cap \mathrm{inn}\left( A\right) $ is open in $\mathrm{inn}\left(
A\right) $. Since $\mathrm{inn}(A)$ is a separable Banach space, it is the
union of countably many translates of $W_{\beta }\cap \mathrm{inn}\left(
A\right) $. Thus, this implies that%
\begin{equation*}
\mathrm{inn}\left( A\right) \in \boldsymbol{\Sigma }_{3}^{0}(\mathrm{Ph}%
^{1+\beta }\mathrm{aut}\left( A\right) )\text{.}
\end{equation*}%
By \cite[Theorem 6.1]{lupini_complexity_2025} this actually implies%
\begin{equation*}
\mathrm{inn}\left( A\right) \in D(\boldsymbol{\Pi }_{2}^{0})(\mathrm{Ph}%
^{1+\beta }\mathrm{aut}\left( A\right) )\text{,}
\end{equation*}%
concluding the proof.
\end{proof}

\begin{proposition}
\label{Proposition:from-derivations-to-automorphisms}Let $A$ be a separable
unital C*-algebra. Then%
\begin{equation*}
\ell _{\mathrm{aut}}\left( A\right) \leq \ell _{\mathrm{Aut}}\left( A\right) 
\text{.}
\end{equation*}
\end{proposition}

\begin{proof}
If $\ell _{\mathrm{Aut}}(A)\leq 1$, the conclusion follows from Corollary %
\ref{Corollary:automorphism-length}. If $\ell _{\mathrm{Aut}}(A)>1$, the
conclusion follows by combining Lemma \ref{Lemma:comparison-ordinal}, Lemma %
\ref{Lemma:comparison-sigma}, and Lemma \ref{Lemma:comparison-D}, together
with Theorem \ref{Theorem:lupini-complexity}.
\end{proof}

\subsection{From automorphisms to derivations}

We continue to denote by $A$ a separable unital C*-algebra. We also consider
the pseudomorphism 
\begin{equation*}
\left( W,U,\varphi \right) :\mathrm{Ph}^{1}\mathrm{out}\left( A\right)
\rightarrow \mathrm{Ph}^{1}\mathrm{Out}\left( A\right)
\end{equation*}%
defined in Proposition \ref{Proposition:from-derivation-to-automorphism},
and adopt the notation from Proposition \ref{Proposition:local-comparison}.
Recall from Lemma \ref{Lemma:phantom-derivation-order-1} and Lemma \ref%
{Lemma:phantom-automorphism-1} that $\mathrm{Ph}^{1}\mathrm{aut}\left(
A\right) $ is the norm-closure of $\mathrm{inn}\left( A\right) $ in $\mathrm{%
aut}\left( A\right) $ and $\mathrm{Ph}^{1}\mathrm{\mathrm{Aut}}\left(
A\right) $ is the norm-closure of $\mathrm{Inn}\left( A\right) $ in $\mathrm{%
\mathrm{Aut}}\left( A\right) $, and that both these groups are endowed with
the topology induced by the operator norm. The following lemma is what
allows one to transfer complexity classes from $\mathrm{Ph}^{1}\mathrm{aut}%
\left( A\right) $ to $\mathrm{Ph}^{1}\mathrm{\mathrm{Aut}}\left( A\right) $
along $\varphi $. Its point is that $\varphi $ is not merely injective and
continuous on the relevant sets, but a homeomorphism onto a \emph{closed }%
subset of $\mathrm{Ph}^{1}\mathrm{\mathrm{Aut}}\left( A\right) $.

\begin{lemma}
\label{Lemma:exp-homeomorphism}Let $A$ be a separable unital C*-algebra, and
let $0<s<\log 2$. Define%
\begin{equation*}
C_{s}:=\left\{ \delta \in \mathrm{Ph}^{1}\mathrm{aut}\left( A\right)
:\left\Vert \delta \right\Vert \leq s\right\}
\end{equation*}%
and%
\begin{equation*}
O_{s}:=\left\{ \delta \in \mathrm{inn}\left( A\right) :\left\Vert \delta
\right\Vert _{\mathrm{inn}\left( A\right) }<s/2\right\} \text{.}
\end{equation*}%
Then:

\begin{enumerate}
\item $\exp $ maps $C_{s}$ to $\mathrm{Ph}^{1}\mathrm{\mathrm{Aut}}\left(
A\right) $, and $O_{s}$ to $\mathrm{Inn}\left( A\right) $;

\item $\exp $ is injective on $C_{s}$, and it induces a homeomorphism from $%
C_{s}$ onto $\exp \left( C_{s}\right) $ with respect to the topologies of $%
\mathrm{Ph}^{1}\mathrm{aut}\left( A\right) $ and $\mathrm{Ph}^{1}\mathrm{%
\mathrm{Aut}}\left( A\right) $;

\item $\exp \left( C_{s}\right) $ is closed in $\mathrm{Ph}^{1}\mathrm{%
\mathrm{Aut}}\left( A\right) $;

\item $O_{s}\subseteq C_{s}$, and $\exp \left( O_{s}\right) $ is an identity
neighborhood in $\mathrm{Inn}\left( A\right) $.
\end{enumerate}
\end{lemma}

\begin{proof}
Throughout the proof we use that $\mathrm{Ph}^{1}\mathrm{aut}\left( A\right) 
$ is norm-closed in $\mathrm{aut}\left( A\right) $, and hence that $C_{s}$
is norm-closed in $\mathrm{aut}\left( A\right) $.

(1) If $h\in A_{\mathrm{sa}}$ then $\exp \left( \mathrm{ad}\left( ih\right)
\right) =\mathrm{Ad}\left( \exp \left( ih\right) \right) \in \mathrm{Inn}%
\left( A\right) $. If $\delta \in \mathrm{Ph}^{1}\mathrm{aut}\left( A\right) 
$ then $\delta $ is the limit in norm of a sequence $\left( \mathrm{ad}%
\left( ih_{n}\right) \right) $ with $h_{n}\in A_{\mathrm{sa}}$, and hence $%
\exp \left( \delta \right) $ is the limit in norm of $\left( \mathrm{Ad}%
\left( \exp \left( ih_{n}\right) \right) \right) $, as $\exp $ is
norm-continuous. Thus $\exp \left( \delta \right) \in \mathrm{Ph}^{1}\mathrm{%
\mathrm{Aut}}\left( A\right) $.

(2) For $\delta \in C_{s}$ we have that%
\begin{equation*}
\left\Vert \exp \left( \delta \right) -\mathrm{id}_{A} \right\Vert \leq
\sum_{n\geq 1}\frac{\left\Vert \delta \right\Vert ^{n}}{n!}\leq e^{s}-1<1%
\text{.}
\end{equation*}%
The series defining $\mathrm{Log}$ converges uniformly on $\left\{ \beta
:\left\Vert \beta -\mathrm{id}_{A} \right\Vert \leq c\right\} $ for every $%
c<1$, and hence defines a norm-continuous function on $\left\{ \beta
:\left\Vert \beta -\mathrm{id}_{A} \right\Vert <1\right\} $. Furthermore,
since $\left\Vert \delta \right\Vert <\log 2$, the series obtained by
substituting the exponential series into the logarithmic series converges
absolutely, and therefore the formal identity between the exponential and
the logarithmic series gives%
\begin{equation*}
\mathrm{Log}\left( \exp \left( \delta \right) \right) =\delta \text{.}
\end{equation*}%
In particular $\exp $ is injective on $C_{s}$, and the inverse of $\exp
|_{C_{s}}$ is the restriction of $\mathrm{Log}$ to $\exp \left( C_{s}\right) 
$. As both $\exp $ and $\mathrm{Log}$ are norm-continuous, and the
topologies of $\mathrm{Ph}^{1}\mathrm{aut}\left( A\right) $ and $\mathrm{Ph}%
^{1}\mathrm{\mathrm{Aut}}\left( A\right) $ are induced by the norm, this
proves (2).

(3) Suppose that $\left( \delta _{n}\right) $ is a sequence in $C_{s}$ such
that $\left( \exp \left( \delta _{n}\right) \right) $ converges to some $%
\beta $ in $\mathrm{Ph}^{1}\mathrm{\mathrm{Aut}}\left( A\right) $, and hence
in norm. Then $\left\Vert \beta -\mathrm{id}_{A} \right\Vert \leq e^{s}-1<1$%
, and therefore%
\begin{equation*}
\delta _{n}=\mathrm{Log}\left( \exp \left( \delta _{n}\right) \right)
\longrightarrow \mathrm{Log}\left( \beta \right) =:\delta
\end{equation*}%
in norm. Since $C_{s}$ is norm-closed, $\delta \in C_{s}$, and $\beta =\exp
\left( \delta \right) \in \exp \left( C_{s}\right) $ by (2).

(4) If $\delta \in O_{s}$, then the inclusion $\mathrm{inn}\left( A\right)
\rightarrow \mathrm{aut}\left( A\right) $ has norm at most $2$, so%
\begin{equation*}
\left\Vert \delta \right\Vert \leq 2\left\Vert \delta \right\Vert _{\mathrm{%
inn}\left( A\right) }<s\text{.}
\end{equation*}%
Hence $\delta \in C_{s}$. Suppose now that $v\in U(A)$ satisfies $\left\Vert
v-1\right\Vert <2\sin \left( s/4\right) $. Since $\left\vert \exp \left(
i\vartheta \right) -1\right\vert =2\left\vert \sin \left( \vartheta
/2\right) \right\vert $ for $\vartheta \in (-\pi ,\pi ]$, and $s/4<\pi /2$,
the spectrum of $v$ is contained in%
\begin{equation*}
\left\{ \exp \left( i\vartheta \right) :\left\vert \vartheta \right\vert
<s/2\right\} \text{.}
\end{equation*}%
Thus, by continuous functional calculus, $v=\exp \left( ih\right) $ for some 
$h\in A_{\mathrm{sa}}$ with $\left\Vert h\right\Vert <s/2$, whence%
\begin{equation*}
\mathrm{Ad}\left( v\right) =\exp \left( \mathrm{ad}\left( ih\right) \right)
\in \exp \left( O_{s}\right) \text{,}
\end{equation*}%
as $\left\Vert \mathrm{ad}\left( ih\right) \right\Vert _{\mathrm{inn}\left(
A\right) }\leq \left\Vert h\right\Vert <s/2$. Since $\mathrm{Inn}\left(
A\right) $ is endowed with the quotient topology induced by the surjective
homomorphism $U(A)\rightarrow \mathrm{Inn}(A)$, $v\mapsto \mathrm{Ad}\left(
v\right) $, which is therefore an open map, this shows that $\exp \left(
O_{s}\right) $ is an identity neighborhood in $\mathrm{Inn}\left( A\right) $.
\end{proof}

For $\xi <\omega _{1}$, the recursive identity for phantom subgroups
recalled in Lemma \ref{Lemma:metric} gives%
\begin{equation*}
\mathrm{Ph}^{\xi }\left( \mathrm{Ph}^{1}\mathrm{out}\left( A\right) \right) =%
\frac{\mathrm{Ph}^{1+\xi }\mathrm{aut}\left( A\right) }{\mathrm{inn}\left(
A\right) }
\end{equation*}%
and%
\begin{equation*}
\mathrm{Ph}^{\xi }\left( \mathrm{Ph}^{1}\mathrm{Out}\left( A\right) \right) =%
\frac{\mathrm{Ph}^{1+\xi }\mathrm{Aut}\left( A\right) }{\mathrm{Inn}\left(
A\right) }\text{.}
\end{equation*}

We will also use the shift identity 
\begin{equation}
\ell _{\mathrm{Ph}}\left( G/N\right) =1+\ell _{\mathrm{Ph}}\left( \mathrm{Ph}%
^{1}\left( G/N\right) \right)  \label{eq:phantom-length-shift}
\end{equation}%
whenever the phantom length of $G/N$ is at least $1$. This follows from 
\begin{equation*}
\mathrm{Ph}^{\xi }\left( \mathrm{Ph}^{1}\left( G/N\right) \right) =\mathrm{Ph%
}^{1+\xi }\left( G/N\right)
\end{equation*}%
and the definition of the terminal Borel class; see \cite[Section 2]%
{solecki_polish_1999}.

\begin{lemma}
\label{Lemma:intrinsic-local-estimates}The following assertions hold.

\begin{enumerate}
\item There exist constants $C>0$ and $t_{0}>0$ such that, whenever $w\in 
\mathrm{Inn}\left( A\right)$ and $q\left( w,1\right) <t<t_{0}$, there exists 
$e\in \mathrm{inn}\left( A\right)$ in the principal exponential chart such
that 
\begin{equation*}
w=\exp \left( e\right) \quad \text{and}\quad \left\Vert e\right\Vert _{%
\mathrm{inn}\left( A\right)}\leq Ct\text{.}
\end{equation*}

\item The Banach--Lie algebra $\mathrm{inn}\left( A\right) $ is an ideal of $%
\mathrm{Ph}^{1}\mathrm{aut}\left( A\right) $, and 
\begin{equation}
\left\Vert \left[ \delta ,e\right] \right\Vert _{\mathrm{inn}\left( A\right)
}\leq \left\Vert \delta \right\Vert \left\Vert e\right\Vert _{\mathrm{inn}%
\left( A\right) }  \label{eq:intrinsic-bracket-estimate}
\end{equation}%
for $\delta \in \mathrm{Ph}^{1}\mathrm{aut}\left( A\right) $ and $e\in 
\mathrm{inn}\left( A\right) $. Consequently, for every sufficiently small $%
R>0$, there is an increasing modulus $\omega _{R}$, with $\omega _{R}\left(
t\right) \rightarrow 0$ as $t\rightarrow 0$, such that 
\begin{equation}
e\ast d-d\in \mathrm{inn}\left( A\right) \quad \text{and}\quad \left\Vert
e\ast d-d\right\Vert _{\mathrm{inn}\left( A\right) }\leq \omega
_{R}(\left\Vert e\right\Vert _{\mathrm{inn}\left( A\right) })
\label{eq:intrinsic-BCH-modulus}
\end{equation}%
whenever $d\in \mathrm{inn}\left( A\right) $, $\left\Vert d\right\Vert _{%
\mathrm{inn}\left( A\right) }\leq R$, and $e\in \mathrm{inn}\left( A\right) $
is sufficiently small.
\end{enumerate}
\end{lemma}

\begin{proof}
(1) Write $w=\mathrm{Ad}\left( u\right) $. By the definition of $q$, after
replacing a central unitary by its adjoint if necessary, there exists $z\in
U\left( ZA\right) $ such that 
\begin{equation*}
\left\Vert uz-1\right\Vert <t\text{.}
\end{equation*}%
For $t_{0}<1$, the principal logarithm of $uz$ is defined. Hence $uz=\exp
\left( ih\right) $ for some $h\in A_{\mathrm{sa}}$ satisfying 
\begin{equation*}
\left\Vert h\right\Vert \leq 2\arcsin \left( t/2\right) \leq Ct\text{.}
\end{equation*}%
Set $e=\mathrm{ad}\left( ih\right) $. Since $z$ is central, 
\begin{equation*}
\exp \left( e\right) =\mathrm{Ad}\left( \exp \left( ih\right) \right) =%
\mathrm{Ad}\left( uz\right) =w
\end{equation*}%
and $\left\Vert e\right\Vert _{\mathrm{inn}\left( A\right)}\leq \left\Vert
h\right\Vert \leq Ct$. Shrinking $t_{0}$ places $e$ in the principal chart.

(2) If $e=\mathrm{ad}\left( ih\right) $, then 
\begin{equation*}
\left[ \delta ,e\right] =\mathrm{ad}\left( i\delta \left( h\right) \right) 
\text{.}
\end{equation*}%
Every derivation annihilates $ZA$. Therefore, for $z\in ZA_{\mathrm{sa}}$, 
\begin{equation*}
\mathrm{dist}\left( \delta \left( h\right) ,ZA_{\mathrm{sa}}\right) \leq
\left\Vert \delta \left( h-z\right) \right\Vert \leq \left\Vert \delta
\right\Vert \left\Vert h-z\right\Vert \text{.}
\end{equation*}%
Taking the infimum over $z$ proves \eqref{intrinsic-bracket-estimate}.

We now deduce the conclusion concerning the Hausdorff series from the
Banach--Lie ideal estimates in Subsection \ref{Subsection:ideals}. Recall
from \eqref{ambient-intrinsic-derivation} that $\left\Vert d\right\Vert \leq
2\left\Vert d\right\Vert _{\mathrm{inn}\left( A\right) }$, and fix $0<R<%
\frac{1}{4}\log 2$, so that the ambient norms of $d$ and of the (smaller) $e$
are below $2^{-1}\log 2$ and the Hausdorff series $e\ast d$ is defined. The
difference $e\ast d-d$ is a sum of Lie monomials in $e$ and $d$, each of
which contains at least one occurrence of $e$; we claim that each such
monomial $\mathfrak{m}$, of total degree $n\geq 1$, satisfies%
\begin{equation}
\mathfrak{m}\in \mathrm{inn}\left( A\right) \quad \text{and}\quad \left\Vert 
\mathfrak{m}\right\Vert _{\mathrm{inn}\left( A\right) }\leq \left\Vert
e\right\Vert _{\mathrm{inn}\left( A\right) }\rho ^{n-1}\text{, where }\rho
:=\max \left\{ \left\Vert e\right\Vert ,\left\Vert d\right\Vert \right\} 
\text{.}  \label{eq:monomial-intrinsic-bound}
\end{equation}%
Indeed, write $\mathfrak{m}$ in right-nested form and consider the \emph{%
innermost }occurrence of $e$ in it. That occurrence sits inside a bracket $%
\left[ u,e\right] $ or $\left[ e,u\right] $ with $u\in \left\{ e,d\right\} $%
; by \eqref{intrinsic-bracket-estimate} and the antisymmetry of the bracket,
that sub-bracket lies in $\mathrm{inn}\left( A\right) $ with intrinsic norm
at most $\left\Vert u\right\Vert \left\Vert e\right\Vert _{\mathrm{inn}%
\left( A\right) }$ (if $\mathfrak{m}=e$ itself, there is nothing to prove).
Applying \eqref{intrinsic-bracket-estimate} once for each of the remaining $%
n-2$ entries, all of which are bracketed against an element of $\mathrm{inn}%
\left( A\right) $ and are estimated in the \emph{ambient }norm, yields %
\eqref{monomial-intrinsic-bound}.

By \eqref{monomial-intrinsic-bound}, the series obtained from $e\ast d-d$ by
replacing each monomial by its intrinsic norm is dominated by $\left\Vert
e\right\Vert _{\mathrm{inn}\left( A\right) }$ times the majorant of the
Hausdorff series provided by the Fundamental Estimate \cite[Theorem 5.29]%
{bonfiglioli_topics_2012}, which converges because $\left( e,d\right) \in
\Delta _{A}$ \cite[Theorem 5.30]{bonfiglioli_topics_2012}. Hence $e\ast
d-d\in \mathrm{inn}\left( A\right) $ and%
\begin{equation*}
\left\Vert e\ast d-d\right\Vert _{\mathrm{inn}\left( A\right) }\leq
K_{R}\left\Vert e\right\Vert _{\mathrm{inn}\left( A\right) }
\end{equation*}%
for a constant $K_{R}$ depending only on $R$. This gives %
\eqref{intrinsic-BCH-modulus} with the linear modulus $\omega _{R}\left(
t\right) :=K_{R}t$.
\end{proof}

\begin{lemma}
\label{Lemma:transfinite-chart-reflection}There exists $R_{\ast }>0$ such
that the following holds simultaneously for every countable ordinal $\xi $
and every $0<R<R_{\ast }$. Suppose that $d_{n}\in \mathrm{inn}\left(
A\right) $, $\left\Vert d_{n}\right\Vert _{\mathrm{inn}\left( A\right)}\leq R
$, and that $\exp \left( d_{n}\right) $ converges to $a$ in the canonical
topology of $\mathrm{Ph}^{1+\xi }\mathrm{Aut}\left( A\right)$. Then there is
a unique $d$ in the fixed closed principal chart such that 
\begin{equation*}
a=\exp \left( d\right) \text{, }d\in \mathrm{Ph}^{1+\xi }\mathrm{aut}\left(
A\right)\text{, and }d_{n}\longrightarrow d\quad \text{in }\mathrm{Ph}%
^{1+\xi }\mathrm{aut}\left( A\right)\text{.}
\end{equation*}
\end{lemma}

\begin{proof}
Choose $R_{\ast }>0$ so small that the closed ambient ball of radius $%
2R_{\ast }$ is contained in a set $C_{s}$ as in Lemma \ref%
{Lemma:exp-homeomorphism}, and so that all the BCH products below belong to
a fixed polydisc on which Lemma \ref{Lemma:intrinsic-local-estimates}(2)
applies. We prove the assertion by transfinite induction, simultaneously for
every $R<R_{\ast }$.

Suppose first that $\xi =0$. By \eqref{ambient-intrinsic-derivation}, $%
\left\Vert d_{n}\right\Vert \leq 2R$. Lemma \ref{Lemma:exp-homeomorphism}%
(2)--(3) gives a unique $d\in \mathrm{Ph}^{1}\mathrm{aut}\left( A\right) $
in the fixed chart such that 
\begin{equation*}
a=\exp \left( d\right) \quad \text{and}\quad d_{n}\longrightarrow d
\end{equation*}%
in the operator norm, which is the canonical topology of $\mathrm{Ph}^{1}%
\mathrm{aut}\left( A\right)$.

Suppose that the assertion holds at $\gamma $, and that 
\begin{equation*}
\exp \left( d_{n}\right) \longrightarrow a\quad \text{in }\mathrm{Ph}%
^{1+\gamma +1}\mathrm{Aut}\left( A\right) \text{.}
\end{equation*}%
The inclusion $\mathrm{Ph}^{1+\gamma +1}\mathrm{Aut}\left(
A\right)\rightarrow \mathrm{Ph}^{1+\gamma }\mathrm{Aut}\left( A\right)$ is
continuous, so the induction hypothesis supplies the unique principal-chart
logarithm $d\in \mathrm{Ph}^{1+\gamma }\mathrm{aut}\left( A\right)$ such
that 
\begin{equation*}
a=\exp \left( d\right) \quad \text{and}\quad d_{n}\longrightarrow d \quad 
\text{in }\mathrm{Ph}^{1+\gamma }\mathrm{aut}\left( A\right)\text{.}
\end{equation*}%
After passing to a tail, choose $\varepsilon _{n}\downarrow 0$ such that 
\begin{equation*}
d_{\gamma +1}^{\mathrm{Ph}^{1}\mathrm{Out}\left( A\right) }\left( \exp
\left( d_{n}\right) ,\exp \left( d\right) \right) <\varepsilon _{n}<1\text{.}
\end{equation*}%
The successor-stage metric formula from Lemma \ref{Lemma:metric} provides,
for every $n$, a sequence $\left( w_{n,k}\right) _{k}$ in $\mathrm{Inn}%
\left( A\right)$ such that 
\begin{equation}
q\left( w_{n,k},1\right) <\varepsilon _{n}\quad \text{and}\quad
w_{n,k}\longrightarrow \exp \left( d\right) \exp \left( -d_{n}\right) \quad 
\text{in }\mathrm{Ph}^{1+\gamma }\mathrm{Aut}\left( A\right)\text{.}
\label{eq:small-automorphism-corrections}
\end{equation}%
By Lemma \ref{Lemma:intrinsic-local-estimates}(1), after discarding finitely
many $n$, we can write 
\begin{equation*}
w_{n,k}=\exp \left( e_{n,k}\right) \text{, }e_{n,k}\in \mathrm{inn}\left(
A\right)\text{, and }\left\Vert e_{n,k}\right\Vert _{\mathrm{inn}\left(
A\right)}\leq C\varepsilon _{n}\text{.}
\end{equation*}%
Define $c_{n,k}:=e_{n,k}\ast d_{n}$. Then 
\begin{equation*}
\exp \left( c_{n,k}\right) =w_{n,k}\exp \left( d_{n}\right) \longrightarrow
\exp \left( d\right) \quad \text{in }\mathrm{Ph}^{1+\gamma }\mathrm{Aut}%
\left( A\right)\text{,}
\end{equation*}%
and Lemma \ref{Lemma:intrinsic-local-estimates}(2) gives 
\begin{equation}
c_{n,k}-d_{n}\in \mathrm{inn}\left( A\right)\quad \text{and}\quad \left\Vert
c_{n,k}-d_{n}\right\Vert _{\mathrm{inn}\left( A\right)}\leq r_{n}:=\omega
_{R}\left( C\varepsilon _{n}\right) ,  \label{eq:small-inner-BCH-corrections}
\end{equation}%
where $r_{n}\rightarrow 0$, uniformly in $k$.

Fix $R^{\prime }$ with $R<R^{\prime }<R_{\ast }$. For all sufficiently large 
$n$, $R+r_{n}<R^{\prime }$. For each such fixed $n$, the sequence $\left(
c_{n,k}\right) _{k}$ lies in the intrinsic $\mathrm{inn}\left( A\right)$%
-ball of radius $R^{\prime }$. The induction hypothesis at stage $\gamma $,
applied to this sequence, gives 
\begin{equation*}
c_{n,k}\longrightarrow d\quad \text{in }\mathrm{Ph}^{1+\gamma }\mathrm{aut}%
\left( A\right)\text{.}
\end{equation*}%
Consequently, 
\begin{equation}
d-d_{n}\in \overline{\left\{ e\in \mathrm{inn}\left( A\right):\left\Vert
e\right\Vert _{\mathrm{inn}\left( A\right)} \leq r_{n}\right\}}^{\mathrm{Ph}%
^{1+\gamma }\mathrm{aut}\left( A\right)}\text{.}
\label{eq:additive-inner-correction}
\end{equation}%
The successor-stage membership criterion and metric formula now imply 
\begin{equation*}
d\in \mathrm{Ph}^{1+\gamma +1}\mathrm{aut}\left( A\right)\quad \text{and}%
\quad d_{\gamma +1}^{\mathrm{Ph}^{1}\mathrm{out}\left( A\right) }\left(
d_{n},d\right) \leq r_{n}\longrightarrow 0\text{.}
\end{equation*}%
This proves the successor step.

Finally, let $\lambda $ be a nonzero countable limit ordinal, and choose
successor ordinals $\xi _{j}\uparrow \lambda $. Convergence in $\mathrm{Ph}%
^{1+\lambda }\mathrm{Aut}\left( A\right)$ implies convergence in every $%
\mathrm{Ph}^{1+\xi _{j}}\mathrm{Aut}\left( A\right)$. The induction
hypothesis and uniqueness in the fixed principal chart yield a single $d$
such that 
\begin{equation*}
a=\exp \left( d\right) \text{, }d\in \bigcap_{j}\mathrm{Ph}^{1+\xi _{j}}%
\mathrm{aut}\left( A\right)=\mathrm{Ph}^{1+\lambda }\mathrm{aut}\left(
A\right) \text{, and }d_{n}\longrightarrow d\quad \text{in every }\mathrm{Ph}%
^{1+\xi _{j}}\mathrm{aut}\left( A\right)\text{.}
\end{equation*}%
The projective-limit formula for $d_{\lambda }^{\mathrm{Ph}^{1}\mathrm{out}%
\left( A\right) }$ and dominated convergence give $d_{n}\rightarrow d$ in $%
\mathrm{Ph}^{1+\lambda }\mathrm{aut}\left( A\right)$. This completes the
induction.
\end{proof}

\begin{proposition}
\label{Proposition:closed-terminal-exponential-ball}Let $\xi <\omega _{1}$
and suppose that $\mathrm{Ph}^{1+\xi }\mathrm{aut}\left( A\right)=\mathrm{inn%
}\left( A\right)$. For every sufficiently small $R>0$, the set 
\begin{equation*}
F_{R}:=\left\{ \exp \left( d\right) :d\in \mathrm{inn}\left( A\right)\text{
and }\left\Vert d\right\Vert _{\mathrm{inn}\left( A\right)}\leq R\right\}
\end{equation*}%
is closed in $\mathrm{Ph}^{1+\xi }\mathrm{Aut}\left( A\right)$.
Consequently, 
\begin{equation*}
\mathrm{Inn}\left( A\right)\in \boldsymbol{\Sigma }_{2}^{0}\left( \mathrm{Ph}%
^{1+\xi }\mathrm{Aut}\left( A\right)\right) \quad \text{and}\quad \ell _{%
\mathrm{Ph}}\left( \mathrm{Ph}^{1}\mathrm{Out}\left( A\right)\right) \leq
\xi +1/2\text{.}
\end{equation*}
\end{proposition}

\begin{proof}
Suppose that $\exp \left( d_{n}\right) \in F_{R}$ converges to $a$ in $%
\mathrm{Ph}^{1+\xi }\mathrm{Aut}\left( A\right)$. Lemma \ref%
{Lemma:transfinite-chart-reflection} gives 
\begin{equation*}
a=\exp \left( d\right) \text{, }d\in \mathrm{Ph}^{1+\xi }\mathrm{aut}\left(
A\right)=\mathrm{inn}\left( A\right)\text{, and }d_{n}\longrightarrow d\quad 
\text{in }\mathrm{Ph}^{1+\xi }\mathrm{aut}\left( A\right)\text{.}
\end{equation*}%
The canonical topology on the underlying group $\mathrm{Ph}^{1+\xi }\mathrm{%
aut}\left( A\right)=\mathrm{inn}\left( A\right)$ and the intrinsic Banach
topology of $\mathrm{inn}\left( A\right)$ are Polish group topologies making
the inclusion into $\mathrm{Ph}^{1}\mathrm{aut}\left( A\right)$ continuous.
They therefore coincide by the uniqueness of the Polish group topology
witnessing that $\mathrm{inn}\left( A\right)$ is a Polish subgroup of $%
\mathrm{Ph}^{1}\mathrm{aut}\left( A\right)$; see \cite[Theorems 9.10 and 15.1%
]{kechris_classical_1995}. Hence 
\begin{equation*}
\left\Vert d_{n}-d\right\Vert _{\mathrm{inn}\left( A\right)}\longrightarrow 0
\end{equation*}%
and $\left\Vert d\right\Vert _{\mathrm{inn}\left( A\right)}\leq R$. Thus $%
a\in F_{R}$, proving that $F_{R}$ is closed in $\mathrm{Ph}^{1+\xi }\mathrm{%
Aut}\left( A\right)$.

By Lemma \ref{Lemma:exp-homeomorphism}(4), for sufficiently small $r<R$, the
set 
\begin{equation*}
\left\{ \exp \left( d\right) :d\in \mathrm{inn}\left( A\right)\text{ and }%
\left\Vert d\right\Vert _{\mathrm{inn}\left( A\right)}<r\right\}
\end{equation*}%
contains an identity neighborhood $V$ in the intrinsic Polish group $\mathrm{%
Inn}\left( A\right)$. Since $\mathrm{Inn}\left( A\right)$ is separable,
there are $g_{j}\in \mathrm{Inn}\left( A\right)$ such that 
\begin{equation*}
\mathrm{Inn}\left( A\right)=\bigcup_{j}g_{j}V\subseteq
\bigcup_{j}g_{j}F_{R}\subseteq \mathrm{Inn}\left( A\right)\text{.}
\end{equation*}%
Translations preserve the canonical topology, so the preceding cover
expresses $\mathrm{Inn}\left( A\right) $ as a countable union of closed
subsets of the ambient group. This proves the first conclusion, and the
second follows from the definition of phantom length.
\end{proof}

\begin{proposition}
\label{Proposition:half-stage-equivalence}For every countable ordinal $\xi $%
, 
\begin{equation*}
\mathrm{inn}\left( A\right)\in \boldsymbol{\Sigma }_{2}^{0}\left( \mathrm{Ph}%
^{1+\xi }\mathrm{aut}\left( A\right)\right) \quad \Longleftrightarrow \quad 
\mathrm{Inn}\left( A\right)\in \boldsymbol{\Sigma }_{2}^{0}\left( \mathrm{Ph}%
^{1+\xi }\mathrm{Aut}\left( A\right)\right) \text{.}
\end{equation*}
\end{proposition}

\begin{proof}
The reverse implication is Lemma \ref{Lemma:comparison-sigma}, applied at
stage $\xi $.

For the forward implication, suppose that 
\begin{equation*}
\mathrm{inn}\left( A\right)\in \boldsymbol{\Sigma }_{2}^{0}\left( \mathrm{Ph}%
^{1+\xi }\mathrm{aut}\left( A\right)\right) \text{.}
\end{equation*}%
By the definition of phantom length and the recursive identity in Lemma \ref%
{Lemma:metric}, this implies $\mathrm{Ph}^{1+\xi +1}\mathrm{aut}\left(
A\right)=\mathrm{inn}\left( A\right)$. The canonical topology on the
underlying group $\mathrm{Ph}^{1+\xi +1}\mathrm{aut}\left( A\right)=\mathrm{%
inn}\left( A\right)$ and the intrinsic Banach topology of $\mathrm{inn}%
\left( A\right)$ are Polish group topologies making the inclusion into $%
\mathrm{Ph}^{1}\mathrm{aut}\left( A\right)$ continuous. They therefore
coincide by the uniqueness of the Polish group topology witnessing that $%
\mathrm{inn}\left( A\right)$ is a Polish subgroup of $\mathrm{Ph}^{1}\mathrm{%
aut}\left( A\right)$; see \cite[Theorems 9.10 and 15.1]%
{kechris_classical_1995}.

Let $d_{\xi +1}^{\mathrm{Ph}^{1}\mathrm{out}\left( A\right) }$ be the
canonical metric on $\mathrm{Ph}^{1+\xi +1}\mathrm{aut}\left( A\right)=%
\mathrm{inn}\left( A\right)$. Choose $0<s<1$ so small that the open
canonical ball 
\begin{equation*}
B_{s}:=\left\{ d\in \mathrm{inn}\left( A\right):d_{\xi +1}^{\mathrm{Ph}^{1}%
\mathrm{out}\left( A\right) }\left( d,0\right) <s\right\}
\end{equation*}%
is contained in an intrinsically bounded subset of the fixed principal chart
on which Lemma \ref{Lemma:transfinite-chart-reflection} applies. This is
possible because the canonical and intrinsic topologies on $\mathrm{inn}%
\left( A\right)$ coincide. The set $B_{s}$ also contains an intrinsic
identity neighborhood in $\mathrm{inn}\left( A\right)$.

The successor-stage metric formula from Lemma \ref{Lemma:metric}, together
with the two implications concerning the infimum recorded immediately before
that lemma, gives 
\begin{equation}
B_{s}=\mathrm{inn}\left( A\right)\cap \bigcup_{\substack{ r\in \mathbb{Q} 
\\ 0<r<s}}\overline{\left\{ e\in \mathrm{inn}\left( A\right): \left\Vert
e\right\Vert _{\mathrm{inn}\left( A\right)}\leq r\right\}}^{\mathrm{Ph}%
^{1+\xi }\mathrm{aut}\left( A\right)}\text{.}
\label{eq:canonical-successor-ball}
\end{equation}%
Indeed, if the canonical distance of $d$ from zero is less than $s$, choose
a rational $r$ strictly between them; then $r$ belongs to the defining set
of radii. Conversely, membership in the closure for some $r<s$ makes the
canonical distance at most $r$. Here we also use $s<1$ and the fixed
truncation of the intrinsic metric. Since $\mathrm{inn}\left( A\right)$ is $%
\boldsymbol{\Sigma }_{2}^{0}$ in $\mathrm{Ph}^{1+\xi }\mathrm{aut}\left(
A\right)$, \eqref{canonical-successor-ball} expresses $B_s$ as a countable
union of closed subsets of $\mathrm{Ph}^{1+\xi }\mathrm{aut}\left( A\right)$%
. Write 
\begin{equation*}
B_{s}=\bigcup_{m\in \mathbb{N}}B_{m}
\end{equation*}%
with every $B_{m}$ closed in $\mathrm{Ph}^{1+\xi }\mathrm{aut}\left(
A\right) $ and contained in $B_s$.

We claim that $\exp \left( B_{m}\right) $ is closed in $\mathrm{Ph}^{1+\xi }%
\mathrm{Aut}\left( A\right)$ for every $m$. Indeed, suppose that $d_{n}\in
B_{m}$ and 
\begin{equation*}
\exp \left( d_{n}\right) \longrightarrow a\quad \text{in }\mathrm{Ph}^{1+\xi
}\mathrm{Aut}\left( A\right) \text{.}
\end{equation*}%
The sequence $\left( d_{n}\right) $ is uniformly bounded in the intrinsic
norm and lies in the fixed principal chart. Lemma \ref%
{Lemma:transfinite-chart-reflection} supplies $d\in \mathrm{Ph}^{1+\xi }%
\mathrm{aut}\left( A\right)$ such that 
\begin{equation*}
a=\exp \left( d\right) \quad \text{and}\quad d_{n}\longrightarrow d\quad 
\text{in }\mathrm{Ph}^{1+\xi }\mathrm{aut}\left( A\right)\text{.}
\end{equation*}%
Since $B_{m}$ is closed in $\mathrm{Ph}^{1+\xi }\mathrm{aut}\left( A\right)$%
, we have $d\in B_{m}$, proving the claim. Thus, $\exp \left( B_{s}\right) $
is $\boldsymbol{\Sigma }_{2}^{0}$ in $\mathrm{Ph}^{1+\xi }\mathrm{Aut}\left(
A\right)$.

Choose $r>0$ such that the intrinsic ball $O_{r}$ from Lemma \ref%
{Lemma:exp-homeomorphism}(4) is contained in $B_{s}$. Then $\exp \left(
O_{r}\right) \subseteq \exp \left( B_{s}\right) $ contains an identity
neighborhood in the intrinsic Polish group $\mathrm{Inn}\left( A\right)$.
Since $\mathrm{Inn}\left( A\right)$ is separable, countably many $\mathrm{Inn%
}\left( A\right)$-translates of $\exp \left( B_{s}\right) $ cover $\mathrm{%
Inn}\left( A\right)$. Translations are homeomorphisms of $\mathrm{Ph}^{1+\xi
}\mathrm{Aut}\left( A\right)$, and hence 
\begin{equation*}
\mathrm{Inn}\left( A\right)\in \boldsymbol{\Sigma }_{2}^{0}\left( \mathrm{Ph}%
^{1+\xi }\mathrm{Aut}\left( A\right)\right) \text{.}
\end{equation*}%
This concludes the proof.
\end{proof}

\begin{corollary}
\label{Corollary:half-level-equivalence}Let $A$ be a separable unital
C*-algebra. For every countable ordinal $\alpha $, 
\begin{equation*}
\ell _{\mathrm{aut}}\left( A\right) \leq \alpha +1/2 \quad
\Longleftrightarrow \quad \ell _{\mathrm{Aut}}\left( A\right) \leq \alpha
+1/2\text{.}
\end{equation*}
\end{corollary}

\begin{proof}
The reverse implication follows from Proposition \ref%
{Proposition:from-derivations-to-automorphisms}. For the forward
implication, suppose first that $\alpha =0$. By Definition \ref%
{Definition:derivation-length}, the derivation length of $A$ is either $0$
or of the form $1+\beta \geq 1$; the hypothesis therefore gives $\ell _{%
\mathrm{aut}}\left( A\right) =0$, and Theorem \ref{Theorem:main-unital}
gives $\ell _{\mathrm{Aut}}\left( A\right) \leq 1/2$.

Suppose now that $\alpha >0$. The same argument applies if $\ell _{\mathrm{%
aut}}\left( A\right) =0$, so assume otherwise; by Definition \ref%
{Definition:derivation-length} again, $\ell _{\mathrm{aut}}\left( A\right)
\geq 1$, so that the shift identity~(\ref{eq:phantom-length-shift}) applies.
Choose a countable ordinal $\xi $ such that $\alpha =1+\xi $ (take the
predecessor when $\alpha $ is finite, and take $\xi =\alpha $ when $\alpha $
is infinite). By Definition \ref{Definition:derivation-length} and the shift
identity~(\ref{eq:phantom-length-shift}), 
\begin{equation*}
\ell _{\mathrm{Ph}}\left( \mathrm{Ph}^{1}\mathrm{out}\left( A\right)\right)
\leq \xi +1/2\text{.}
\end{equation*}%
Thus, 
\begin{equation*}
\mathrm{inn}\left( A\right)\in \boldsymbol{\Sigma }_{2}^{0}\left( \mathrm{Ph}%
^{1+\xi }\mathrm{aut}\left( A\right)\right) \text{.}
\end{equation*}%
Proposition \ref{Proposition:half-stage-equivalence} gives 
\begin{equation*}
\mathrm{Inn}\left( A\right)\in \boldsymbol{\Sigma }_{2}^{0}\left( \mathrm{Ph}%
^{1+\xi }\mathrm{Aut}\left( A\right)\right) \text{,}
\end{equation*}%
and hence 
\begin{equation*}
\ell _{\mathrm{Ph}}\left( \mathrm{Ph}^{1}\mathrm{Out}\left( A\right)\right)
\leq \xi +1/2\text{.}
\end{equation*}%
Proposition \ref{Proposition:from-derivations-to-automorphisms} ensures that 
$\ell _{\mathrm{Aut}}\left( A\right) \geq 1$, so another application of the
shift identity yields 
\begin{equation*}
\ell _{\mathrm{Aut}}\left( A\right) \leq 1+\left( \xi +1/2\right) =\alpha
+1/2\text{.}
\end{equation*}%
This proves the forward implication.
\end{proof}

Following \cite[Section 4]{lupini_complexity_2025}, for a Polish subgroup $H$
of a Polish group $G$, write $\mathrm{rk}_{\mathrm{S}}(H,G)$ for the \emph{%
Solecki rank} of $H$ in $G$, namely the least countable ordinal $\xi $ such
that $\mathrm{Ph}^{\xi }(G/H)=\{\ast\}$. We define the \emph{Solecki rank}
of a homogeneous space with a Polish cover $G/N$ to be the Solecki rank of $%
N $ in $G$, and write%
\begin{equation*}
\mathrm{rk}_{\mathrm{S}}(G/N):=\mathrm{rk}_{\mathrm{S}}(N,G) \text{.}
\end{equation*}

\begin{corollary}
\label{Corollary:Solecki-ranks}For every countable ordinal $\xi $, 
\begin{equation*}
\mathrm{Ph}^{1+\xi }\mathrm{Aut}\left( A\right)=\mathrm{Inn}\left(
A\right)\quad \Longrightarrow \quad \mathrm{Ph}^{1+\xi }\mathrm{aut}\left(
A\right)=\mathrm{inn}\left( A\right) \quad \Longrightarrow \quad \mathrm{Ph}%
^{1+\xi +1}\mathrm{Aut}\left( A\right)=\mathrm{Inn}\left( A\right)\text{.}
\end{equation*}%
Consequently, 
\begin{equation*}
\mathrm{rk}_{\mathrm{S}}\left( \mathrm{Ph}^{1}\mathrm{out}(A)\right) \leq 
\mathrm{rk}_{\mathrm{S}}\left( \mathrm{Ph}^{1}\mathrm{Out}(A)\right) \leq 
\mathrm{rk}_{\mathrm{S}}\left( \mathrm{Ph}^{1}\mathrm{out}(A)\right) +1\text{.}
\end{equation*}
\end{corollary}

\begin{proof}
If $\mathrm{Ph}^{1+\xi }\mathrm{Aut}\left( A\right)=\mathrm{Inn}\left(
A\right)$, then the automorphism length is at most $1+\xi $. Lemma \ref%
{Lemma:comparison-ordinal} gives that the derivation length is at most $%
1+\xi $, and hence $\mathrm{Ph}^{1+\xi }\mathrm{aut}\left( A\right)=\mathrm{%
inn}\left( A\right)$. If $\mathrm{Ph}^{1+\xi }\mathrm{aut}\left( A\right)=%
\mathrm{inn}\left( A\right)$, Proposition \ref%
{Proposition:closed-terminal-exponential-ball} gives 
\begin{equation*}
\mathrm{Inn}\left( A\right)\in \boldsymbol{\Sigma }_{2}^{0}\left( \mathrm{Ph}%
^{1+\xi }\mathrm{Aut}\left( A\right)\right) \text{,}
\end{equation*}%
so the recursive identity in Lemma \ref{Lemma:metric} gives $\mathrm{Ph}%
^{1+\xi +1}\mathrm{Aut}\left( A\right)=\mathrm{Inn}\left( A\right)$. The
inequalities follow by taking the least such stages.
\end{proof}

\begin{lemma}
\label{Lemma:complexity}Suppose that $A$ is a separable unital C*-algebra,
and $2\leq \eta <\omega _{1}$.

\begin{enumerate}
\item If 
\begin{equation*}
\mathrm{inn}\left( A\right) \in \boldsymbol{\Sigma }_{\eta }^{0}(\mathrm{Ph}%
^{1}\mathrm{aut}\left( A\right) )\text{,}
\end{equation*}%
then 
\begin{equation*}
\mathrm{Inn}\left( A\right) \in \boldsymbol{\Sigma }_{\eta }^{0}(\mathrm{Ph}%
^{1}\mathrm{Aut}\left( A\right) )\text{.}
\end{equation*}

\item If 
\begin{equation*}
\mathrm{inn}\left( A\right) \in D(\boldsymbol{\Pi }_{\eta }^{0})(\mathrm{Ph}%
^{1}\mathrm{aut}\left( A\right) ),
\end{equation*}%
then 
\begin{equation*}
\mathrm{Inn}\left( A\right) \in \boldsymbol{\Sigma }_{\eta +1}^{0}(\mathrm{Ph%
}^{1}\mathrm{Aut}\left( A\right) ).
\end{equation*}
\end{enumerate}
\end{lemma}

\begin{proof}
Fix $s\in \left( 0,r_{0}\right) $, and let $C_{s}$ and $O_{s}$ be defined as
in Lemma \ref{Lemma:exp-homeomorphism}. Then $s<\log 2$ and $C_{s}\subseteq
W $. Define%
\begin{equation*}
O:=\mathrm{inn}\left( A\right) \cap C_{s}\text{.}
\end{equation*}%
Then $\varphi $ and $\exp $ agree on $C_{s}$, and by Lemma \ref%
{Lemma:exp-homeomorphism}(1) and (4),%
\begin{equation*}
\varphi \left( O_{s}\right) \subseteq \varphi \left( O\right) \subseteq 
\mathrm{Inn}\left( A\right) \text{,}
\end{equation*}%
where $\varphi \left( O_{s}\right) $ is an identity neighborhood in $\mathrm{%
Inn}\left( A\right) $. As $\mathrm{Inn}\left( A\right) $ is separable, it
follows that $\mathrm{Inn}\left( A\right) $ is the union of countably many
left translates of $\varphi \left( O\right) $ by elements of $\mathrm{Inn}%
\left( A\right) $, and such translations are homeomorphisms of $\mathrm{Ph}%
^{1}\mathrm{\mathrm{Aut}}\left( A\right) $.

(1) Since $C_{s}$ is closed in $\mathrm{Ph}^{1}\mathrm{aut}\left( A\right) $%
, we have that%
\begin{equation*}
O\in \boldsymbol{\Sigma }_{\eta }^{0}\left( \mathrm{Ph}^{1}\mathrm{aut}%
\left( A\right) \right)
\end{equation*}%
and, in particular, $O$ belongs to $\boldsymbol{\Sigma }_{\eta }^{0}$
relatively to $C_{s}$. By Lemma \ref{Lemma:exp-homeomorphism}(2), $\varphi
|_{C_{s}}$ is a homeomorphism from $C_{s}$ onto $\varphi \left( C_{s}\right) 
$, and hence $\varphi \left( O\right) $ belongs to $\boldsymbol{\Sigma }%
_{\eta }^{0}$ relatively to $\varphi \left( C_{s}\right) $. Since $\varphi
\left( C_{s}\right) $ is closed in $\mathrm{Ph}^{1}\mathrm{\mathrm{Aut}}%
\left( A\right) $ by Lemma \ref{Lemma:exp-homeomorphism}(3), this implies
that%
\begin{equation*}
\varphi \left( O\right) \in \boldsymbol{\Sigma }_{\eta }^{0}\left( \mathrm{Ph%
}^{1}\mathrm{Aut}\left( A\right) \right) \text{.}
\end{equation*}%
Since $\mathrm{Inn}\left( A\right)$ is the union of countably many
translates of $\varphi \left( O\right) $, we conclude that 
\begin{equation*}
\mathrm{Inn}\left( A\right) \in \boldsymbol{\Sigma }_{\eta }^{0}\left( 
\mathrm{Ph}^{1}\mathrm{Aut}\left( A\right) \right) \text{.}
\end{equation*}

(2) Adopting the notation from (1), the same argument gives 
\begin{equation*}
O\in D(\boldsymbol{\Pi }_{\eta }^{0})(\mathrm{Ph}^{1}\mathrm{aut}\left(
A\right) )\quad \text{and}\quad \varphi \left( O\right) \in D(\boldsymbol{%
\Pi }_{\eta }^{0})(\mathrm{Ph}^{1}\mathrm{Aut}\left( A\right) )\text{,}
\end{equation*}%
where we use that $D(\boldsymbol{\Pi }_{\eta }^{0})$, like $\boldsymbol{%
\Sigma }_{\eta }^{0}$, is closed under intersections with closed sets. Since 
$\mathrm{Inn}\left( A\right) $ is the union of countably many translates of $%
\varphi \left( O\right) $, and $D(\boldsymbol{\Pi }_{\eta }^{0})\subseteq 
\boldsymbol{\Sigma }_{\eta +1}^{0}$, we conclude%
\begin{equation*}
\mathrm{Inn}\left( A\right) \in \boldsymbol{\Sigma }_{\eta +1}^{0}\left( 
\mathrm{Ph}^{1}\mathrm{Aut}\left( A\right) \right) \text{.}
\end{equation*}
\end{proof}

\begin{proposition}
\label{Proposition:from-automorphisms-to-derivations}Let $A$ be a separable
unital C*-algebra, and let $\alpha <\omega _{1}$ be a countable ordinal.

\begin{enumerate}
\item $\ell _{\mathrm{aut}}\left( A\right) \leq \ell _{\mathrm{\mathrm{Aut}}%
}\left( A\right) $;

\item $\ell _{\mathrm{aut}}\left( A\right) \leq \alpha +1/2$ if and only if $%
\ell _{\mathrm{Aut}}\left( A\right) \leq \alpha +1/2$;

\item $\ell _{\mathrm{aut}}\left( A\right) \leq \alpha +1/2+\varepsilon $ if
and only if $\ell _{\mathrm{Aut}}\left( A\right) \leq \alpha
+1/2+\varepsilon $;

\item If $\ell _{\mathrm{aut}}\left( A\right) =\alpha $, then 
\begin{equation*}
\ell _{\mathrm{Aut}}\left( A\right) \in \left\{ \alpha ,\alpha +1/2\right\} ;
\end{equation*}

\item If $\ell _{\mathrm{aut}}\left( A\right) =\alpha +1/2$, then $\ell _{%
\mathrm{Aut}}\left( A\right) =\alpha +1/2$;

\item $\ell _{\mathrm{aut}}\left( A\right) =\alpha +1/2+\varepsilon $ if and
only if $\ell _{\mathrm{Aut}}\left( A\right) =\alpha +1/2+\varepsilon $.
\end{enumerate}
\end{proposition}

\begin{proof}
We will also apply Theorem \ref{Theorem:main-unital}.

(1) This is by Proposition \ref%
{Proposition:from-derivations-to-automorphisms}.

(2) This is Corollary \ref{Corollary:half-level-equivalence}.

(3) The reverse implication follows from (1). For the forward implication,
suppose first that $\alpha =0$. Definition \ref{Definition:derivation-length}
shows that the hypothesis forces $\ell _{\mathrm{aut}}\left( A\right) =0$,
so Theorem \ref{Theorem:main-unital} gives $\ell _{\mathrm{Aut}}\left(
A\right) \leq 1/2$. The same argument applies whenever $\ell _{\mathrm{aut}%
}\left( A\right) =0$, so assume from now on that $\alpha >0$ and $\ell _{%
\mathrm{aut}}\left( A\right) >0$.

Suppose that $\alpha $ is a limit ordinal. The shift identity and Theorem %
\ref{Theorem:lupini-complexity} give 
\begin{equation*}
\mathrm{inn}\left( A\right) \in \Gamma _{\alpha +1/2+\varepsilon }\left( 
\mathrm{Ph}^{1}\mathrm{aut}\left( A\right) \right) =D(\boldsymbol{\Pi }%
_{\alpha +1}^{0})\left( \mathrm{Ph}^{1}\mathrm{aut}\left( A\right) \right) 
\text{.}
\end{equation*}%
Lemma \ref{Lemma:complexity} and \cite[Theorem 6.1]{lupini_complexity_2025}
therefore give 
\begin{equation*}
\mathrm{Inn}\left( A\right) \in D(\boldsymbol{\Pi }_{\alpha +1}^{0})\left( 
\mathrm{Ph}^{1}\mathrm{Aut}\left( A\right) \right) =\Gamma _{\alpha
+1/2+\varepsilon }\left( \mathrm{Ph}^{1}\mathrm{Aut}\left( A\right) \right) 
\text{.}
\end{equation*}%
Since $\alpha $ is a limit ordinal, $\alpha +1/2+\varepsilon $ belongs to $%
\omega _{1}^{\mathbf{\Delta }}$ and is not of the form $\sigma +1/2$ for a
successor ordinal $\sigma $, so the converse implication of Theorem \ref%
{Theorem:lupini-complexity} applies; followed by the shift identity and $%
1+\alpha =\alpha $, it gives 
\begin{equation*}
\ell _{\mathrm{Aut}}\left( A\right) \leq \alpha +1/2+\varepsilon \text{.}
\end{equation*}

It remains to consider the case where $\alpha $ is a successor ordinal.
Choose $\gamma $ such that $\alpha =1+\gamma $ (take the predecessor when $%
\alpha $ is finite, and take $\gamma =\alpha $ when $\alpha $ is infinite).
Thus, $\gamma $ is either zero or a successor ordinal. The phantom length of 
$\mathrm{Ph}^{1}\mathrm{out}\left( A\right) $ is at most $\gamma
+1/2+\varepsilon $. Since $\Gamma _{\gamma +1/2+\varepsilon }=D(\boldsymbol{%
\Pi }_{1+\gamma +1}^{0})$ for every $\gamma <\omega _{1}$, the forward
implication of Theorem \ref{Theorem:lupini-complexity}---which, we stress,
applies to every element of $\omega _{1}^{\mathrm{Ph}}$, and not only to the
elements of $\omega _{1}^{\mathbf{\Delta }}$---gives%
\begin{equation*}
\mathrm{inn}\left( A\right) \in \Gamma _{\gamma +1/2+\varepsilon }\left( 
\mathrm{Ph}^{1}\mathrm{aut}\left( A\right) \right) =D(\boldsymbol{\Pi }%
_{1+\gamma +1}^{0})\left( \mathrm{Ph}^{1}\mathrm{aut}\left( A\right) \right) 
\text{.}
\end{equation*}%
Then by Lemma \ref{Lemma:complexity}, 
\begin{equation*}
\mathrm{Inn}\left( A\right) \in \boldsymbol{\Sigma }_{1+\gamma +2}^{0}(%
\mathrm{Ph}^{1}\mathrm{Aut}\left( A\right) )\text{.}
\end{equation*}%
By \cite[Theorem 6.1]{lupini_complexity_2025}, this implies%
\begin{equation*}
\mathrm{Inn}\left( A\right) \in D(\boldsymbol{\Pi }_{1+\gamma +1}^{0})\left( 
\mathrm{Ph}^{1}\mathrm{Aut}\left( A\right) \right) =\Gamma _{\gamma
+1/2+\varepsilon }\left( \mathrm{Ph}^{1}\mathrm{Aut}\left( A\right) \right) 
\text{.}
\end{equation*}%
If $\gamma =0$, then $\gamma +1/2+\varepsilon $ belongs to $\omega _{1}^{%
\mathbf{\Delta }}$ and is not of the form $\sigma +1/2$ for a successor
ordinal $\sigma $, so the converse implication of Theorem \ref%
{Theorem:lupini-complexity} applies verbatim. If $\gamma >0$, then $\gamma $
is a successor ordinal and $\Gamma _{\gamma +1/2+\varepsilon }=\Gamma
_{\gamma +1/2}$, so that $\mathrm{Inn}\left( A\right) \in \Gamma _{\gamma
+1/2}\left( \mathrm{Ph}^{1}\mathrm{\mathrm{Aut}}\left( A\right) \right) $,
and the exceptional case of that converse implication applies at $\gamma
+1/2 $. In either case one obtains 
\begin{equation*}
\ell _{\mathrm{Ph}}\left( \mathrm{Ph}^{1}\mathrm{Out}\left( A\right) \right)
\leq \gamma +1/2+\varepsilon \text{,}
\end{equation*}%
and the shift identity yields 
\begin{equation*}
\ell _{\mathrm{Aut}}\left( A\right) \leq 1+\gamma +1/2+\varepsilon =\alpha
+1/2+\varepsilon \text{.}
\end{equation*}

(4) If $\ell _{\mathrm{aut}}\left( A\right) =\alpha $, then (1) and (2) give 
\begin{equation*}
\alpha \leq \ell _{\mathrm{Aut}}\left( A\right) \leq \alpha +1/2\text{.}
\end{equation*}%
There is no value of the formal phantom scale strictly between $\alpha $ and 
$\alpha +1/2$.

(5) If $\ell _{\mathrm{aut}}\left( A\right) =\alpha +1/2$, then (1) and (2)
give 
\begin{equation*}
\alpha +1/2\leq \ell _{\mathrm{Aut}}\left( A\right) \leq \alpha +1/2\text{.}
\end{equation*}

(6) The forward implication follows from (1) and (3). Conversely, suppose
that 
\begin{equation*}
\ell _{\mathrm{Aut}}\left( A\right) =\alpha +1/2+\varepsilon \text{.}
\end{equation*}%
By (1), the derivation length is at most $\alpha +1/2+\varepsilon $. If it
were at most $\alpha +1/2$, then (2) would force the automorphism length to
be at most $\alpha +1/2$, a contradiction. Since there is no value of the
formal phantom scale strictly between $\alpha +1/2$ and $\alpha
+1/2+\varepsilon $, the derivation length is $\alpha +1/2+\varepsilon $.
\end{proof}

\begin{theorem}
\label{Theorem:equality}Let $A$ be a separable unital C*-algebra. Then
either 
\begin{equation*}
\ell _{\mathrm{Aut}}\left( A\right) =\ell _{\mathrm{aut}}\left( A\right)
\end{equation*}%
or there is a countable ordinal $\alpha \neq 1$ such that 
\begin{equation*}
\ell _{\mathrm{aut}}\left( A\right) =\alpha \quad \text{and}\quad \ell _{%
\mathrm{Aut}}\left( A\right) =\alpha +1/2\text{.}
\end{equation*}
\end{theorem}

\begin{proof}
Every value in the formal phantom scale is an ordinal, an element of the
form $\alpha +1/2$, or an element of the form $\alpha +1/2+\varepsilon $.
The classification follows from Proposition \ref%
{Proposition:from-automorphisms-to-derivations}(4)--(6).

If $\ell _{\mathrm{aut}}\left( A\right) =1$, then Proposition \ref%
{Proposition:KLR} gives 
\begin{equation*}
\ell _{\mathrm{Aut}}\left( A\right) \leq 1\text{,}
\end{equation*}%
so one must have $\alpha \neq 1$ in the second option. Likewise, by \cite[%
Proposition 10.3 and Corollary 10.4]{lupini_complexity_2025}, $\ell _{%
\mathrm{aut}}\left( A\right) $ cannot be a limit ordinal.
\end{proof}

\begin{corollary}
\label{Corollary:non-limit}If $A$ is a separable unital C*-algebra, then
neither the derivation length nor the automorphism length of $A$ is a limit
ordinal.
\end{corollary}

\begin{proof}
If every derivation of $A$ is inner, then its derivation length is zero.
Otherwise, Definition \ref{Definition:derivation-length} gives 
\begin{equation*}
\ell _{\mathrm{aut}}\left( A\right) =1+\ell _{\mathrm{Ph}}\left( \mathrm{Ph}%
^{1}\mathrm{out}\left( A\right) \right)\text{.}
\end{equation*}%
The phantom Banach space $\mathrm{Ph}^{1}\mathrm{out}\left( A\right) $ has
Borel length that is not a limit ordinal by \cite[Proposition 10.3 and
Corollary 10.4]{lupini_complexity_2025}. If the derivation length were a
nonzero limit ordinal $\lambda $, then the phantom length in the display
would also be $\lambda $, and Theorem \ref{Theorem:lupini-complexity} would
force its Borel length to be $\lambda $, a contradiction. Thus, the
derivation length of $A$ is not a limit ordinal. By Theorem \ref%
{Theorem:equality}, an ordinal-valued automorphism length must equal the
derivation length. Hence the automorphism length is not a limit ordinal
either.
\end{proof}

%\appendix
%    Include appendix "chapters" here.
%\include{}

%    Bibliography styles amsplain or harvard are also acceptable.
\bibliographystyle{amsplain}
\bibliography{bibliography}
%    See note above about multiple indexes.

%\appendix
%    Include appendix "chapters" here.
%\include{}

\end{document}